\documentclass[10pt]{article} 

\usepackage{graphicx}
\usepackage{amsfonts}
\usepackage{amsmath}
\usepackage{amscd}
\usepackage{amssymb}
\usepackage{verbatim}
\usepackage{graphpap}
\usepackage{exscale,relsize}

\newenvironment{proof}{\noindent\emph{Proof}.}{$\square$\smallskip}
\newenvironment{proof*}{\noindent\emph{Proof}}{$\square$\smallskip}
\newenvironment{proofnosquare}{\noindent\emph{Proof}.}{\smallskip}
\newtheorem{theorem}{Theorem}[section]

\newtheorem{Definition}[theorem]{Definition}
\newtheorem{lemma}[theorem]{Lemma}
\newtheorem{Example}[theorem]{Example}
\newtheorem{corollary}[theorem]{Corollary}
\newtheorem{Remark}[theorem]{Remark}
\newtheorem{proposition}[theorem]{Proposition}

\newtheorem{fact}[theorem]{Fact}
\newtheorem{claim}[theorem]{Claim}

\newenvironment{definition}{\begin{Definition}\normalfont}{\end{Definition}}
\newenvironment{example}{\begin{Example}\normalfont}{\end{Example}}
\newenvironment{remark}{\begin{Remark}\normalfont}{\end{Remark}}

\newcommand{\br}{\ensuremath{\mathbb{R}}} 
\newcommand{\bz}{\ensuremath{\mathbb{Z}}} 
\newcommand{\bn}{\ensuremath{\mathbb{N}}} 
\newcommand{\bq}{\ensuremath{\mathbb{Q}}} 
\newcommand{\bc}{\ensuremath{\mathbb{C}}} 
\newcommand{\bm}{\ensuremath{\mathbb{M}}} 
\newcommand{\id}{\ensuremath{\mathrm{id}}} 
\newcommand{\co}{\colon\thinspace} 

\newcommand{\ligx}{\ensuremath{\mathcal{G}_{LI}(X)}} 
\newcommand{\ligy}{\ensuremath{\mathcal{G}_{LI}(Y)}} 
\newcommand{\lsgx}{\ensuremath{\mathcal{G}_{LS}(X)}} 
\newcommand{\lsgy}{\ensuremath{\mathcal{G}_{LS}(Y)}} 
\newcommand{\symtv}{\ensuremath{Sym_\infty(T,v)}} 
\newcommand{\symsw}{\ensuremath{Sym_\infty(S,w)}} 
\newcommand{\plix}{\ensuremath{\mathcal{P}_{LI}(X)}} 
\newcommand{\eligx}{\ensuremath{{\mathcal{G}}^{\epsilon}_{LI}(X)}} 
\newcommand{\eiligx}{\ensuremath{{\mathcal{G}}^{\epsilon_i}_{LI}(X)}} 
\newcommand{\eioligx}{\ensuremath{{\mathcal{G}}^{\epsilon_{i+1}}_{LI}(X)}} 
\newcommand{\pgd}{\ensuremath{{\mathcal{PG}}({D},v_0)}} 
\newcommand{\pgb}{\ensuremath{{\mathcal{PG}}({B(T,v)})}} 
\newcommand{\pg}{\ensuremath{{\mathcal{PG}}}} 
\newcommand{\pb}{\ensuremath{{\mathcal{P}}({B(T,v)})}} 
\newcommand{\gammax}{\ensuremath{\mathcal{G}_{\Gamma}(X)}} 
\newcommand{\gammaxh}{\ensuremath{\mathcal{G}_{h^{-1}\Gamma h}(X)}} 
\newcommand{\gammay}{\ensuremath{\mathcal{G}_{\Gamma}(Y)}} 
\newcommand{\g}{\ensuremath{\mathcal{G}}} 
\newcommand{\Cay}{\ensuremath{{\rm{Cay}}}} 
                            
\title{Trees, Ultrametrics, and \\Noncommutative Geometry}  
\author{Bruce Hughes\\ Vanderbilt University}

\begin{document}

\maketitle

\begin{abstract} Noncommutative geometry is used
to study the local geometry of ultrametric spaces and the geometry of 
trees at infinity. 
Connes's example of the noncommutative space of Penrose tilings
is interpreted as a non-Hausdorff orbit space of a compact,
ultrametric space under the action of its local isometry group.
This is generalized to compact, locally rigid, ultrametric spaces.
The local isometry types and the local similarity types in those
spaces can be analyzed using groupoid $C^*$-algebras.

The concept of a locally rigid action of a countable group $\Gamma$
on a compact, ultrametric space by local similarities is introduced. It is 
proved that there is a faithful unitary representation of $\Gamma$ into
the germ groupoid $C^*$-algebra of the action.
The prototypical example is the standard action of Thompson's group $V$
on the ultrametric Cantor set. In this case, the $C^*$-algebra is 
the Cuntz algebra $\mathcal O_2$
and representations originally due to
Birget and Nekrashevych are recovered.

End spaces of trees are sources of ultrametric spaces. Some connections
are made between locally rigid, ultrametric spaces and a concept in 
the theory of tree lattices
of Bass and Lubotzky.
\end{abstract}

{\footnotesize
\thanks{Supported in part by NSF Grant 
DMS--0504176 \newline
{\hspace*{15pt}}2000 Math.\ Subject Class. Primary 46L85,
54E45, 20F65; 
Secondary 22A22, 05C05, 19K14
}
}

{\footnotesize
\tableofcontents}


\section{Introduction}      

This paper applies techniques from noncommutative geometry, 
in particular, Renault's 
group\-oid $C^*$-algebras \cite{Ren}, to study various aspects of the local
geometry of ultrametric spaces. Since the geometry of infinite 
trees at infinity and the local geometry 
of ultrametric spaces are related by the end space functor, 
results about the large-scale geometry of
infinite trees are also obtained.

There are two main motivating examples for this paper.
First, Connes \cite{Con} used the space of Penrose
tilings as an illustration of a  noncommutative space.
We observe here that the space of Penrose tiles can be interpreted
as a certain compact ultrametric space modulo local isometry type.
We then address the problem of finding other compact ultrametric 
spaces for which the noncommutative geometric point of view can 
be used to study local isometry types.
To solve this problem, locally rigid ultrametric spaces are 
introduced.
The results are 
described in Section~\ref{subsec:LI} below.

In addition to the equivalence relation of local isometry type,
we also analyze a closely related equivalence relation on the 
points of a locally rigid ultrametric space, namely local similarity type.
Those results are 
summarized in Section~\ref{subsec:LS}.

The second main motivating example
is Birget's faithful unitary representation of Thompson's group $V$
into the Cuntz algebra ${\cal O}_2$
\cite{Bir}, also obtained independently by 
Nekrashevych \cite{Nek}.
In this paper we derive such a representation from a more general
result establishing a faithful unitary representation of any countable
group $\Gamma$ acting locally rigidly by local similarities
on a compact ultrametric space $X$ into the $C^*$-algebra of a 
groupoid associated to the action of $\Gamma$ on $X$.
See Section~\ref{subsec:FUR}
for more details.

The concept of a rigid tree appears in the Bass-Lubotzky theory of
tree lattices. In Section~\ref{subsec:tree-lattices} it is observed that this 
condition is equivalent to the end space of the tree being a locally
rigid, ultrametric space.

\subsection{Local isometries}
\label{subsec:LI}
The first main motivating example for this paper is the description by Connes \cite{Con} of the space of Penrose tilings. 
Consider the space
$X$ of infinite sequences of $0$'s and $1$'s, where any $1$ must be followed by 
$0$; the topology on $X$ comes from considering $X$ as a subspace of the countable
product of the discrete space $\{ 0,1\}$: 
$$X = \left\{ (x_0,x_1,x_2,\dots) \in  \prod_0^\infty\{ 0,1\} ~|~ \text{$x_i=1$ implies $x_{i+1}=0$} \right\}.$$

It is known that 
the set of Penrose tilings is parametrized as
a quotient space of $X$ with respect to the equivalence relation $R$
of {\it tail equivalence}
(where 
$x=(x_i)_{i=0}^\infty$ and $y=(y_i)_{i=0}^\infty$ are tail equivalent if and only if
there exists $N\geq 0$ such that $x_i=y_i$ for all $i\geq N$)
(see Gr\"unbaum and Shephard \cite{GruS}).
The problem is that $X/R$ is not Hausdorff and, hence, cannot be studied 
by ordinary topological methods.
Nevertheless, Connes shows how to associate to $X/R$ 
a natural noncommutative $C^*$-algebra; that is to say, 
$X/R$ can be viewed as 
a {\it noncommutative 
topological space} (whereas the Gelfand-Naimark Theorem
shows that commutative $C^*$-algebras correspond to locally compact
Hausdorff spaces).
In fact, the $C^*$-algebra constructed
by Connes for $X/R$ is AF (that is, {\it approximately 
finite}, or the norm closure of a direct limit of finite dimensional matrix algebras over
$\bc$).

\begin{figure}[htbp]

\begin{picture}(300,150)(-50,-10)

\put(110,130){$0$}\put(190,130){$1$}
\put(40,70){$0$}\put(100,70){$1$}\put(200,70){$0$}
\put(10,35){$0$}\put(50,35){$1$}\put(90,35){$0$}
\put(170,35){$0$}\put(210,35){$1$}
\put(-2,10){$0$}\put(18,10){$1$}\put(38,10){$0$}
\put(78,10){$0$}\put(98,10){$1$}
\put(158,10){$0$}\put(178,10){$1$}\put(198,10){$0$}


\put(0,0){\line(1,2){10}}
\put(40,0){\line(1,2){10}}
\put(80,0){\line(1,2){10}}
\put(160,0){\line(1,2){10}}
\put(200,0){\line(1,2){10}}

\put(20,0){\line(-1,2){10}}
\put(100,0){\line(-1,2){10}}
\put(180,0){\line(-1,2){10}}

\put(10,20){\line(2,3){20}}
\put(90,20){\line(2,3){20}}
\put(170,20){\line(2,3){20}}

\put(50,20){\line(-2,3){20}}
\put(210,20){\line(-2,3){20}}

\put(30,50){\line(1,1){40}}
\put(190,50){\line(1,1){40}}

\put(110,50){\line(-1,1){40}}

\put(70,90){\line(4,3){80}}

\put(230,90){\line(-4,3){80}}

\put(150,150){\circle*{5}}
\put(70,90){\circle*{5}}\put(230,90){\circle*{5}}
\put(30,50){\circle*{5}}\put(110,50){\circle*{5}}\put(190,50){\circle*{5}}
\put(10,20){\circle*{5}}\put(50,20){\circle*{5}}\put(90,20){\circle*{5}}
\put(170,20){\circle*{5}}\put(210,20){\circle*{5}}

\put(0,0){\vector(-1,-2){0}}\put(40,0){\vector(-1,-2){0}}\put(80,0){\vector(-1,-2){0}}
\put(160,0){\vector(-1,-2){0}}\put(200,0){\vector(-1,-2){0}}

\put(20,0){\vector(1,-2){0}}
\put(100,0){\vector(1,-2){0}}
\put(180,0){\vector(1,-2){0}}

\multiput(30,-12)(160,0){2}{\vdots}
\put(90,-12){\vdots}

\end{picture}
\caption{The Fibonacci tree}
\label{fig:fib-tree}
\end{figure}
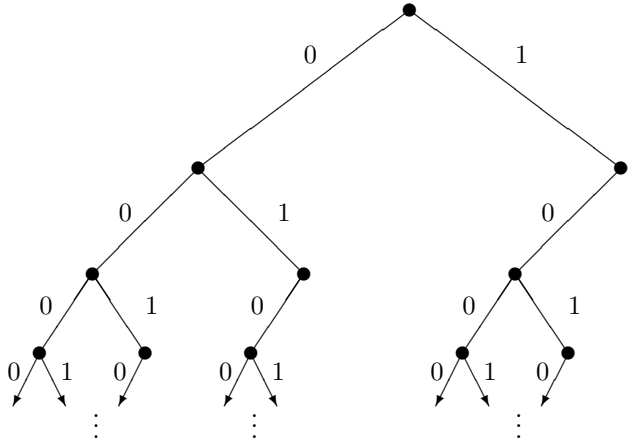

We call $(X,d)$ the {\it Fibonacci space} because it is the end 
space of the Fibonacci tree (see Figure~\ref{fig:fib-tree})
as will be explained below.

\bigskip

%


The point of view of this paper begins with the simple 
observation that the equivalence relation $R$ on $X$ has a
geometric interpretation when $X$ is endowed with the natural metric $d(x,y)=e^{-n}$, 
where $n=\inf\{ i\geq 0 ~|~ x_i\not= y_i\}$.
Namely, $xRy$ if and only if $X$ has the same {\it local isometry type}
at $x$ and $y$ (that is, there exists
$\epsilon >0$ and an isometry $h\co B(x,\epsilon)\to B(y,\epsilon)$ with $hx=y$).

One may then ask, what is a natural class of metric spaces whose local isometry types are able to be studied
by noncommutative geometric methods?
It is this question that we seek to answer in the first part of this paper.

A key feature of the metric $d$ on $X$ is that it is an 
{\it ultrametric}; that is, $d$ satisfies the strong triangle
inequality $d(x,y)\leq\max\{ d(x,z), d(z,y)\}$ for all $x,y,z\in X$.

A important geometric property of the Fibonacci space 
$X$ is that it is {\it rigid}; that is, $X$ has no isometries other than 
the identity. In particular,
local isometry types in $X$ are not reflected by global symmetries in $X$.
Furthermore, the rigidity property of $X$ is inherited
by balls, so that, in particular, $X$ is {\it locally rigid}.

In general, we find that it is compact ultrametric spaces satisfying the local 
rigidity property whose local isometry types can be studied using
noncommutative geometric methods.
The entry into these methods is through the theory of groupoids and their
$C^*$-algebras.

The equivalence relation $R$ on $X$ is an example of a groupoid. In general, one can define a 
groupoid $\ligx$ of local isometries on a metric space $X$. 
If the groupoid is sufficiently well-behaved, Renault \cite{Ren} showed
how to define the $C^*$-algebra of the groupoid (generalizing the 
$C^*$-algebra of a group).
We show that compact, locally rigid, ultrametric spaces have groupoids
of local isometries to which Renault's theory can be applied.

The results we are able to obtain in this general situation are summarized in the following theorem.

\begin{theorem} 
\label{theorem:first-main}
If $X$ is a compact, locally rigid, ultrametric space and $\ligx$ is the groupoid of local isometries on $X$, then
\begin{enumerate}
\item
\label{first-main-one}
$\ligx$ is a locally compact, Hausdorff, 
second countable, \'etale groupoid;
\item 
\label{first-main-two}
the groupoid $C^*$-algebra $C^*\ligx$ is a unital AF $C^*$-algebra;
\item 
\label{first-main-five}
the topological groupoid $\ligx$, 
the unital groupoid $C^*$-algebra $C^*\ligx$, and the unital, partially
ordered abelian group $K_0C^*\ligx$ are each invariants
of $X$ up to micro-scale equivalence of $X$;
\item 
\label{first-main-three}
there exists a 
Bratteli diagram $B(X)$ such that $\ligx$ is the path groupoid of $B(X)$;
\item
\label{first-main-four}
$K_0C^*\ligx$ as a unital partially ordered abelian group
is isomorphic to the symmetry at infinity group
$\symtv$ of any  rooted, geodesically complete,
locally finite simplicial tree $(T,v)$ whose end space is isometric to $X$.
\end{enumerate}
\end{theorem}

To measure local isometry types in $X$, it might seem more natural to focus on the group $LI(X)$ of local isometries
from $X$ to itself, rather than 
the groupoid $\ligx$.
One of the main purposes of this paper is to show that the groupoid approach can be used to study local
isometries on a compact, locally rigid ultrametric space.
In this case, the groupoid $\ligx$ is an effective replacement of the quotient space $X/LI(X)$, which, in general, need
not be Hausdorff. This is quite different from what happens for the isometry group $Isom(X)$, in which case
$X/Isom(X)$ is a perfectly well-behaved space.

As with the Fibonacci space and tree, there is in general a well-known
relationship between ultrametric spaces and trees. In fact, compact 
ultrametric spaces are exactly the end spaces of rooted, proper
$\br$-trees.
Under this correspondence, local isometries of the end space come from 
{\it uniform isometries at infinity}
of the rooted tree. These are isometries between complements of open balls
in the tree centered at the root. 
It follows that Theorem~\ref{theorem:first-main} provides an invariant
of rooted, geodesically complete, proper $\br$-trees, 
with locally rigid end spaces,
up to uniform isometry at infinity. 
See Corollary~\ref{unif-iso-inf-inv}
for more details.

Elliott \cite{Ell} proved that a unital AF
$C^*$-algebra is determined up to 
isomorphism by its $K_0$ group (as a unital partially ordered abelian group).
Consequently, in light of
Theorem~\ref{theorem:first-main}(\ref{first-main-two}), 
knowing $K_0C^*\ligx$
becomes important.
For any geodesically complete,
locally finite simplicial tree $T$ with root $v$,  
a group $Sym_\infty(T,v)$ is introduced
in Section~\ref{section:TheSymmetryAtInfinityGroup}, called {\it the
group of symmetries at infinity of $T$}. 
It is a unital partially ordered abelian group, constructed as a direct limit
of a sequence of finitely generated free abelian groups. The direct sequence is
elementary to construct from a diagram of the tree. 
Item (\ref{first-main-four}) in Theorem~\ref{theorem:first-main} is established by showing
that $Sym_\infty(T,v)$ is  isomorphic to $K_0C^*{\g_{LI}}(end(T,v))$
as a unital partially
ordered abelian group. 

\begin{figure}[htbp]


\begin{picture}(300,160)(-25,-10)

\put(110,130){$0$}\put(190,130){$1$}
\put(40,70){$0$}\put(100,70){$1$}\put(200,70){$0$}\put(260,70){$1$}
\put(10,35){$0$}\put(50,35){$1$}\put(90,35){$0$}\put(130,35){$1$}\put(170,35){$0$}\put(210,35){$1$}\put(250,35){$0$}\put(290,35){$1$}
\put(-2,10){$0$}\put(18,10){$1$}\put(38,10){$0$}\put(58,10){$1$}\put(78,10){$0$}\put(98,10){$1$}\put(118,10){$0$}\put(138,10){$1$}
\put(158,10){$0$}\put(178,10){$1$}\put(198,10){$0$}\put(218,10){$1$}\put(238,10){$0$}\put(258,10){$1$}\put(278,10){$0$}\put(298,10){$1$}


\put(0,0){\line(1,2){10}}
\put(40,0){\line(1,2){10}}
\put(80,0){\line(1,2){10}}
\put(120,0){\line(1,2){10}}
\put(160,0){\line(1,2){10}}
\put(200,0){\line(1,2){10}}
\put(240,0){\line(1,2){10}}
\put(280,0){\line(1,2){10}}

\put(20,0){\line(-1,2){10}}
\put(60,0){\line(-1,2){10}}
\put(100,0){\line(-1,2){10}}
\put(140,0){\line(-1,2){10}}
\put(180,0){\line(-1,2){10}}
\put(220,0){\line(-1,2){10}}
\put(260,0){\line(-1,2){10}}
\put(300,0){\line(-1,2){10}}

\put(10,20){\line(2,3){20}}
\put(90,20){\line(2,3){20}}
\put(170,20){\line(2,3){20}}
\put(250,20){\line(2,3){20}}

\put(50,20){\line(-2,3){20}}
\put(130,20){\line(-2,3){20}}
\put(210,20){\line(-2,3){20}}
\put(290,20){\line(-2,3){20}}

\put(30,50){\line(1,1){40}}
\put(190,50){\line(1,1){40}}

\put(110,50){\line(-1,1){40}}
\put(270,50){\line(-1,1){40}}

\put(70,90){\line(4,3){80}}

\put(230,90){\line(-4,3){80}}

\put(150,150){\circle*{5}}
\put(70,90){\circle*{5}}\put(230,90){\circle*{5}}
\put(30,50){\circle*{5}}\put(110,50){\circle*{5}}\put(190,50){\circle*{5}}\put(270,50){\circle*{5}}
\put(10,20){\circle*{5}}\put(50,20){\circle*{5}}\put(90,20){\circle*{5}}\put(130,20){\circle*{5}}
\put(170,20){\circle*{5}}\put(210,20){\circle*{5}}\put(250,20){\circle*{5}}\put(290,20){\circle*{5}}

\put(0,0){\vector(-1,-2){0}}\put(40,0){\vector(-1,-2){0}}\put(80,0){\vector(-1,-2){0}}\put(120,0){\vector(-1,-2){0}}
\put(160,0){\vector(-1,-2){0}}\put(200,0){\vector(-1,-2){0}}\put(240,0){\vector(-1,-2){0}}\put(280,0){\vector(-1,-2){0}}

\put(20,0){\vector(1,-2){0}}\put(60,0){\vector(1,-2){0}}\put(100,0){\vector(1,-2){0}}\put(140,0){\vector(1,-2){0}}
\put(180,0){\vector(1,-2){0}}\put(220,0){\vector(1,-2){0}}\put(260,0){\vector(1,-2){0}}\put(300,0){\vector(1,-2){0}}

\multiput(30,-12)(80,0){4}{\vdots}

\end{picture}
%
\caption{The Cantor tree}
\label{fig:cantor-tree}
\end{figure}
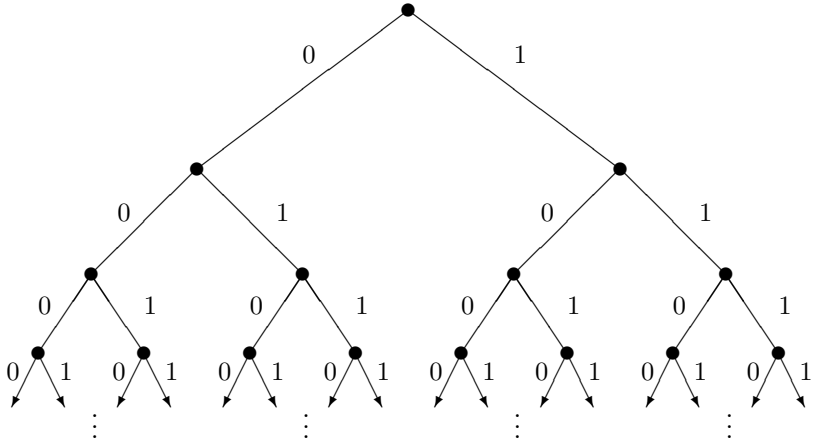

Theorem~\ref{theorem:first-main} only applies in the locally rigid case;
however, there are many compact ultrametric spaces that are not locally rigid.
An important example is the end space of the Cantor tree in 
Figure~\ref{fig:cantor-tree} (see Examples~\ref{ex:Cantor Tree}).

%
To a certain extent noncommutative geometric methods can still be applied
to such spaces, and this is discussed in Section~\ref{subsec:FUR}
below. But first we turn to the second main result of this paper.

\subsection{Local similarities}
\label{subsec:LS}

In addition to studying local isometry types in compact ultrametric spaces, one may
relax this relation and instead study {\it local similarity types.}
When considering end spaces of rooted trees, just as local isometries
of the end space correspond to uniform isometries at infinity of
the rooted tree, local similarities of 
the end space correspond to
(not necessarily uniform) {\it isometries at infinity} of the rooted tree.
These are isometries between complements of the interiors of finite
subtrees of the tree containing the root.

In the case that $X$ is the Fibonacci space, 
then $x,y\in X$ are {\it tail equivalent with lag} (i.e., there exists $m,n\geq 0$
such that $x_{m+j}=y_{n+j}$ for all $j\geq 0$) if and only if 
$X$ has the same local similarity type at $x$ and $y$ (i.e., there exist $\epsilon, \lambda >0$ and a $\lambda$-similarity
$h\co B(x,\epsilon)\to B(y,\lambda\epsilon)$ with $hx=y$) (see \cite{Hug}).

In analogy with the groupoid of local isometries, one can define a groupoid
$\lsgx$ of local similarities on a metric space $X$.
We show that compact, locally rigid, ultrametric spaces have groupoids of local
similarities to which Renault's theory can be applied. 
The results in this general situation are summarized in the following 
theorem.

\begin{theorem} 
\label{theorem:ls}
If $X$ is a compact, locally rigid, ultrametric space and $\lsgx$ is the groupoid of local similarities on $X$, then
\begin{enumerate}
\item $\lsgx$ is a locally compact, Hausdorff, second countable, \'etale groupoid,
and 
\item $\lsgx$ is invariant up to local similarity of $X$.
\end{enumerate}
\end{theorem}

The contrast between local isometry types (Theorem~\ref{theorem:first-main})
and local similarity types (Theorem~\ref{theorem:ls}) is already
foreshadowed in the $C^*$-algebra literature. For example, Mingo \cite{Min} studied 
$C^*$-algebras of spaces of sequences modulo tail equivalence (generalizing Connes \cite{Con}).
On the other hand, Kumjian, Pask, Raeburn and Renault \cite{KPRR} studied Cuntz-Krieger
$C^*$-algebras of spaces of sequences modulo tail equivalence with lag. 
Roughly,
tail equivalence corresponds to local isometry type, 
whereas tail equivalence with
lag corresponds to local similarity type. For the end space of the
Fibonacci tree (and some other trees), the analogy  is exact.


\subsection{Faithful unitary representations}
\label{subsec:FUR}

The second main motivating example for this paper is Birget's faithful 
unitary representation of Thompson's group $V$
into the Cuntz algebra ${\cal O}_2$
\cite{Bir}.
Such a representation  was obtained independently by 
Nekrashevych \cite{Nek}.
This is related to this paper in the following two ways:
\begin{enumerate}
\item Thompson's group $V$ is a subgroup of the group of
local similarities on the end space $Y$ of the Cantor tree. 
(References for this are given in Section~\ref{subsec:Thompson}.)
\item Renault \cite{Ren} defined a groupoid $O_2$, called the {\it Cuntz groupoid}, 
and showed that the $C^*$-algebra of $O_2$ is the Cuntz algebra
${\cal O}_2$. The groupoid $O_2$ is easily seen to be a groupoid of local similarities on the end space of the 
Cantor tree.
\end{enumerate}

However, as pointed out above,
the end space of the Cantor tree is not locally rigid; therefore, the point of view as developed in 
Theorems~\ref{theorem:first-main} and \ref{theorem:ls} does not apply directly.

The key observation needed to overcome the lack of local rigidity is that in items (1) and (2) above,
not {\it all} local similarities of the end space $Y$ of the Cantor tree are needed---just 
those that {\it locally preserve the natural total order}. 
The group $\Gamma$ of local similarities of $Y$ that are locally order preserving has an important property:
it {\it acts locally rigidly on $Y$.} (See 
Section~\ref{subsection:LocallyRigidActions} for the definition and
Section~\ref{subsec:Thompson} for the fact that this action is 
locally rigid.)
This is the key property shared with the full group $LS(X)$ of all local 
similarities on a compact, locally rigid ultrametric space $X$. 

Thus, we are led to consider subgroups $\Gamma$ of the group $LS(X)$
of local similarities of an arbitrary compact ultrametric space $X$ 
that act locally rigidly on $X$.

In analogy with the groupoids of local isometries and local similarities, 
one can define a subgroupoid $\gammax$ of
$\lsgx$ whenever $\Gamma$ is a subgroup of the group of 
local similarities on a metric space $X$.
In the case that $X$ is a compact ultrametric space and
the action of $\Gamma$ on $X$ is locally rigid,
we show that  Renault's theory can be applied. 
The results in this general situation are summarized in the following 
theorem.

\begin{theorem}
\label{theorem:second-main}
If $X$ is a compact, ultrametric space and $\Gamma$ is a countable group acting locally rigidly on $X$
by local similarities, then 
\begin{enumerate}
\item the germ groupoid $\gammax$ of $\Gamma$ on $X$ 
is a locally compact, Hausdorff, second countable, \'etale groupoid;
\item if $h\in LS(X)$, then $h^{-1}\Gamma h$ also acts locally rigidly on $X$ by local similarities and $\gammax$ and 
$\gammaxh$
are isomorphic topological groupoids; 
\item there is a faithful unitary representation of $\Gamma$ into $C^*\gammax$.
\end{enumerate}
\end{theorem}

\begin{example}
If $Y$ is the end space of the Cantor tree and $\Gamma$ is the subgroup of $LS(Y)$ consisting of locally
order preserving local similarities, then $\Gamma=V$, Thompson's group, and $\Gamma$ acts locally rigidly on $Y$.
In this case, $\gammay=O_2$, the Cuntz groupoid.
Since $C^*\gammay={\cal O}_2$, the Cuntz algebra by Renault \cite{Ren}, 
the representation of Birget \cite{Bir} and Nekrashevych \cite{Nek}
is a special 
case of Theorem~\ref{theorem:second-main}.

More generally, Birget and Nekrashevych obtained faithful unitary 
representations of the Higman--Thompson groups $G_{n,1}$ into the Cuntz
algebra ${\cal O}_n$ for all $n\geq 2$ ($G_{2,1}=V)$. Such representations
also follow from the results in this 
paper---see Section~\ref{subsec:Thompson}.
\end{example}

If $X$ is a compact, locally rigid ultrametric space, then $\Gamma=LI(X)$ acts locally rigidly on $X$ and
$\ligx=\gammax$. Thus, Theorems~\ref{theorem:first-main} and \ref{theorem:second-main} imply the following result.

\begin{corollary} The local isometry group $LI(X)$ of a compact, locally rigid ultrametric
space $X$ has a faithful unitary representation into the AF $C^*$-algebra $C^*\ligx$.
\end{corollary}

It should be mentioned that Berestovskii \cite{Ber} proved that the
isometry group of any 
compact ultrametric space has
a faithful representation into the orthogonal group of a 
separable, real  Hilbert space.

\subsection{Guide}

We indicate where the proofs of the theorems in the introduction
may be found in the body of the paper.
The proof of Theorem~\ref{theorem:first-main} can be found as follows.
Item (\ref{first-main-one}) is 
proved in Section~\ref{subsection:LocallyRigidUltrametricSpaces}; 
item (\ref{first-main-two}) is proved in 
Section~\ref{section:TheApproximatingGroupoids}; 
items (\ref{first-main-five}) and (\ref{first-main-three}) are proved in 
Section~\ref{section:Perturbations};
item (\ref{first-main-four}) is given
by Corollary~\ref{sym-and-k0}.

The proof of Theorem~\ref{theorem:ls} can be found as follows.
Item (1) is given by Corollary~\ref{cor:ls-grpd-props}.
Item (2) is in Section~\ref{section:LocSimLocIsomGrpd}.

The proof of Theorem~\ref{theorem:second-main} can be found as follows.
Item (1) is 
given in Section~\ref{section:LR,LRA,H} 
(see Theorem~\ref{theorem:second-main-one-proof}).
Item (2) is proved in Section~\ref{section:LocSimLocIsomGrpd}
(Proposition~\ref{prop:ConjugateGroupoids})
and Section~\ref{section:LR,LRA,H} (Proposition~\ref{prop:ConjugateActions}).
Item (3)
is in Section~\ref{section:FUR}.

For the theory of $C^*$-algebras of groupoids, see  Muhly \cite{Muh}, Paterson \cite{Pat}, and
Renault \cite{Ren}.
For some general background on ultrametric spaces, see
Khrennikov \cite{Khr04} and
Robert \cite{Rob}.
In addition to other references in the body of the paper, see
Nikolaev \cite{Nik} for an exposition of AF $C^*$-algebras and their
$K$-theory.

Throughout this paper we use the notation
$\bn = \{ 1,2,3,\dots\}$ for the natural numbers,
$\bz$ for the integers, $\bz_+$ for the nonnegative integers,
$\br$ for the real numbers, and $\bc$ for the complex numbers.

\bigskip
\noindent
{\bf Acknowledgements.} I have benefited from conversations with
Berndt Brenken, Jon Brown, Alex Kumjian, Igor Nikolaev, Mark Sapir,
Jack Spielberg, Andreas Thom, and Guoliang Yu.


\section{Local isometry and local similarity groups}
\label{section:LocSimLocIsom}

If $(X,d)$ is a metric space, $x\in X$ and $\epsilon >0$, then we use the notation
$B(x,\epsilon)=\{ y\in X ~|~ d(x,y)<\epsilon\}$ for the {\it open ball about $x$ of radius
$\epsilon$}, and $\bar B(x,\epsilon)=\{ y\in X ~|~ d(x,y)\leq\epsilon\}$ 
for the {\it closed ball about $x$ of radius
$\epsilon$}.

\begin{definition} Let $(X,d_X)$ and $(Y,d_Y)$ be metric spaces.
A homeomorphism $h\co X\to Y$ is
\begin{enumerate}
\item an {\it isometry} if $d_Y(hx,hy)=d_X(x,y)$ for all $x,y\in X$.
\item a {\it similarity} if there exists $\lambda >0$ such that $d_Y(hx,hy)=\lambda d_X(x,y)$ for all $x,y\in X$.
In this case, $h$ is a {\it $\lambda$-similarity}, $\lambda$ is the {\it similarity modulus} of $h$, and we
write $sim(h)=\lambda$.
\item a {\it local isometry} if for every $x\in X$ there exists $\epsilon >0$ such that
$h$ restricts to an isometry $h|\co B(x,\epsilon)\to B(hx,\epsilon)$.
\item a {\it local similarity} if for every $x\in X$ there exist $\epsilon >0$ and $\lambda >0$ 
such that $h$ restricts to a $\lambda$-similarity $h|\co B(x,\epsilon)\to B(hx,\lambda\epsilon)$. In this case,
$\lambda$ is the {\it similarity modulus of $h$ at $x$} and we write $sim(h,x)=\lambda$.
\item a {\it uniform local similarity} if 
there exist $\epsilon >0$ and $\lambda > 0$
such that for every $x\in X$ the restriction $h|\co B(x,\epsilon)\to
B(hx,\lambda\epsilon)$ is a $\lambda$-similarity. In this case,
$sim(h,x)=\lambda$ for all $x\in X$.
\end{enumerate}
\end{definition}

An important point of this definition is that each of these maps is surjective.\footnote{The terminology conflicts slightly
with \cite{Hug}, where {\it similarity}, {\it local similarity},
and {\it  uniform local similarity} were modified by 
{\it equivalence} to indicate that they were homeomorphisms.}

For a local similarity $h$, the similarity modulus $sim(h,x)$ is uniquely determined by $h$ and $x$, except in the
case $x$ is an isolated point of $X$. In that case, we will always take $sim(h,x)=1$.

For a metric space $X$ (with a given metric),
$Homeo(X)$ is the group of homeomorphisms from $X$ to $X$,
$Isom(X)$ is the group of isometries from $X$ to $X$,
$LI(X)$ is the group of local isometries from $X$ to $X$, and
$LS(X)$ is the group of local similarities from $X$ to $X$.

Note that there are inclusions of subgroups
$$Isom(X)\leq LI(X) \leq LS(X) \leq Homeo(X).$$
These groups are given the compact-open topology unless otherwise specified.

If $X$ is a compact metric space and $h\co X\to X$ is a similarity
(respectively, a uniform local similarity), then $h$ is an isometry
(respectively, a local isometry) (for the second of these statements, 
see \cite{Hug}).


\section{Local  isometry and local similarity groupoids}
\label{section:LocSimLocIsomGrpd}

In this section we define the topological groupoids of local similarities 
and local isometries of a metric space.
Unfortunately, these groupoids are rarely Hausdorff or second countable (see
Examples~\ref{example:non-Hausdorff} and \ref{example:non-secondcountable})---two 
conditions needed for Renault's machinery \cite{Ren} 
to work.\footnote{The Hausdorff condition is relaxed a bit in Paterson's approach \cite{Pat}.}
We will eventually overcome this problem in Section~\ref{section:LR,LRA,H}
by either restricting to second countable, locally rigid ultrametric
spaces or to certain subgroupoids.

An alternative treatment of local isometry groupoids is in 
Bridson and Haefliger \cite[Part III, Chapter $\mathcal{G}$]{BrH};
however, beyond the basic definitions, their point of view is quite a bit different from the present paper 
(in particular, they do not discuss $C^*$-algebras).

Let $(X,d)$ be a metric space.

\begin{definition} Let $x_1,x_2\in X$. A {\it local similarity germ from $x_1$ to
$x_2$ in $X$\/} is an equivalence class $[g,x_1]$ represented by a 
$\lambda$-similarity
$g\co B(x_1,\epsilon)\to B(x_2,\lambda\epsilon)$
for some $\epsilon >0$ and $\lambda >0$ such that $gx_1=x_2$.
Another such similarity $g'\co B(x_1,\epsilon')\to B(x_2,\lambda'\epsilon')$
is {\it equivalent\/} to $g$ if
$g|B(x_1,\epsilon'')= g'|B(x_1,\epsilon'')$ for some $\epsilon''>0$ with $\epsilon''\leq\min\{\epsilon,\epsilon'\}$.
\end{definition}

If $[g,x]$ is a local similarity germ, then the modulus $sim(g,x)$ is independent of the choice of
representative for the equivalence class.

\begin{definition} The {\it local similarity groupoid $\lsgx$ of $X$\/} is the set of
all local similarity germs between pairs of points in $X$.
\end{definition}

The groupoid structures on $\lsgx$ are the obvious ones. Thus, the unit space is $X$ and the
domain $d\co \lsgx\to X$ and range $r\co \lsgx\to X$ maps are given by
$d([g,x])=x$ and $r([g,x])=gx$. 
If $[g_1,x_1]$ and $[g_2,x_2]$ are local similarity germs from $x_1$ to $x_2$ and from $x_2$ to
$x_3$, respectively, then the composition $[g_2,x_2][g_1,x_1]$ is the local similarity germ
from $x_1$ to $x_3$ defined by composing $g_1$ and $g_2$ after suitably restricting their
domains: $[g_2, x_2][g_1, x_1] = [g_2g_1, x_1]$. 
The inverse is $[g,x]^{-1}= [g^{-1}, gx]$.

The topology on $\lsgx$ is determined as follows.

\begin{definition}  For every germ
$[g,x]$ represented by a $\lambda$-similarity 
$g\co B(x,\epsilon)\to B(gx,\lambda\epsilon)$ and every 
$y\in B(x,\epsilon)$, 
there is a $\lambda$-similarity $g|\co B(y,\delta)\to B(gy,\lambda\delta)$ where $\delta=
\epsilon-d(x,y)$ representing a germ $[g,y]$. Let 
$$U(g,x,\epsilon)= \left\{ [g,y] ~|~ y\in B(x,\epsilon)\right\}\subseteq \lsgx.$$
The collection of all such $U(g,x,\epsilon)$ for $[g,x]\in\lsgx$ forms a basis for a topology on
$\lsgx$ called the {\it germ topology.\/}
\end{definition}

Note that $U(g_1,x_1,\epsilon_1)\cap U(g_2,x_2,\epsilon_2)=$
$$\bigcup\{ U(g,x,\epsilon) ~|~ B(x,\epsilon)\subseteq B(x_1,\epsilon_1)\cap
B(x_2,\epsilon_2) ~\textrm{and}~
g=g_1|B(x,\epsilon)=g_2|B(x,\epsilon)\}.
$$

Throughout the rest of this paper, $\lsgx$ will always be given the germ topology.

The following result gives a proof of Theorem~\ref{theorem:ls}(2) in the
Introduction.

\begin{proposition}
\label{prop:loc-sim}
A local similarity $h\co X\to Y$ of metric spaces
induces an isomorphism $h_*\co \lsgx\to\lsgy$ of topological
groupoids.
\end{proposition}

\begin{proof}
Let $[g,x]\in\lsgx$ be a local similarity germ. 
We may assume that $g$ is defined on $B(x,\epsilon)$ with $\epsilon> 0$
sufficiently small that there exists $\lambda_1>0$ such that
$h|\co B(x,\epsilon)\to B(hx,\lambda_1\epsilon)$
is a $\lambda_1$-similarity and there exists
$\lambda_2>0$ such that
$h|\co B(gx,\epsilon)\to B(hgx,\lambda_1\epsilon)$
is a $\lambda_2$-similarity.
Then $h_*[g,x]\co =[hgh^{-1},hx]$ can be seen to define an isomorphism
of topological groupoids.
\end{proof}

\begin{lemma}\label{lemma:etale} 
$\lsgx$ is an \'etale groupoid. That is, $r\co \lsgx\to X$ is a local
homeomorphism. In fact, the collection $\lsgx^{op}$ of open subsets $A$ of $\lsgx$ such that $d|A$
and $r|A$ are homeomorphisms onto open subsets of $X$ forms a basis for the germ topology
\footnote{Some authors take this as the definition of \'etale; others refer to it as
{\it $r$-discreteness}. Note that we are not insisting that our groupoids are locally compact,
Hausdorff or second countable.} on
$\lsgx$.
\end{lemma}

\begin{proof}
It suffices to observe that for each $\lambda$-similarity
$g\co B(x,\epsilon)\to B(gx,\lambda\epsilon)$,
$d|\co U(g,x,\epsilon)\to B(x,\epsilon)$ and $r|\co U(g,x,\epsilon)\to B(gx,\lambda\epsilon)$ are
homeomorphisms. It is clear that these are bijections; that they are homeomorphisms follows from the fact that
$U(g,y,\delta)\subseteq U(g,x,\epsilon)$ whenever $y\in B(x,\epsilon)$ and
$0 < \delta\leq \epsilon-d(x,y)$.
\end{proof}

\begin{lemma}\label{lemma:loccom} 
If $X$ is locally compact, then $\lsgx$ is locally compact.
\end{lemma}

\begin{proof}
From the proof of Lemma~\ref{lemma:etale}, $\lsgx$ has a basis of open sets homeomorphic to open balls of $X$.
\end{proof}
 
\begin{example} 
If $(X,d)$ is a discrete metric space, then $\lsgx$ and $X\times X$ 
are isomorphic as topological
groupoids when $X\times X$ is given the pair groupoid structure
(e.g. see \cite[page 747]{Weins}).
\end{example}

\begin{definition} An {\it open subgroupoid} of $\lsgx$ is a subset $\g$ of $\lsgx$ such that
\begin{enumerate}
\item $\g$ is open in $\lsgx$ (as topological spaces),
\item $\g$ is closed under composition,
\item $\g$ contains the unit space $X$.
\end{enumerate}
\end{definition}

\begin{definition} The {\it local isometry groupoid $\ligx$ of $X$} is the groupoid of all
local isometry germs between pairs of points of $X$; that is,
$$\ligx=\{ [g,x]\in\lsgx ~|~ sim(g,x)=1\}.$$
\end{definition}

\begin{definition} For a subgroup $\Gamma\leq LS(X)$, the {\it germ groupoid $\gammax$ of $\Gamma$ on $X$}
is the subgroupoid of $\lsgx$ given by
$$\gammax = \{ [g,x] \in\lsgx ~|~ g\in\Gamma\}.$$
\end{definition}

\begin{remark}\label{remark:opensubgrpoids}
It is clear that both $\ligx$ and $\gammax$ are open subgroupoids of $\lsgx$. In fact,
if $g\co B(x,\epsilon)\to B(gx,\epsilon)$ represents a $[g,x]\in\ligx$, then
$[g,x]\in U(g,x,\epsilon)\subseteq\ligx$.
Likewise, if $g\in\Gamma$, $x\in X$ and $\epsilon >0$,
then $[g,x]\in U(g,x,\epsilon)\subseteq\gammax$.
\end{remark}

\begin{remark}\label{remark:unit-space}
The unit space is naturally an open subspace of $\lsgx$ via the map
$\alpha\co X\to\lsgx, x\mapsto [\id, x]$, where 
$\id\co B(x,\epsilon)\to B(x,\epsilon)$ is the identity
for some $\epsilon > 0$. Obviously, $\alpha$ is injective. To see that it is 
continuous, let $U(g,y,\epsilon)$ be a basis element of $\lsgx$ and suppose $\alpha(x)\in U(g,y,\epsilon)$.
Then $[\id, x]\in U(g,y,\epsilon)$, so $x\in B(y,\epsilon)$ and $g=\id$ near $x$.
It follows that if $z$ is close enough to $x$, then $\alpha(z)\in U(g,y,\epsilon)$,
thereby verifying continuity of $\alpha$.
To see that $\alpha$ is an open map, note that
$\alpha(B(x,\epsilon))=\cup \{[\id,y] ~|~ y\in B(x,\epsilon)\}= U(\id, x, \epsilon)$.
In particular, $\alpha(X)$ is an open subset of $\lsgx$. It need not be the case that $\alpha(X)$ is
closed in $\ligx$, but see Remark~\ref{remark:locally-rigid-unit-space}
for an instance when it is. 
\end{remark}

\begin{example}\label{example:non-Hausdorff}
If $X$ is the end space of the Cantor tree $C$, then $\ligx$ is not Hausdorff.
On the other hand, if $Y$ is the end space of the Fibonacci tree $F$,
then $\ligy$ is Hausdorff.
See Theorem~\ref{theorem:Hausdorff-groupoid} and Example~\ref{example:locally-rigid}.
\end{example}

\begin{example}\label{example:non-secondcountable}
If $X=\{ x\in\br^2 ~|~ ||x||\leq 1\}$, the closed unit ball in $\br^2$ with the usual metric, then
$\ligx$ is not second countable. To see this,
for each $0\leq\theta < 2\pi$, let $g_\theta\co X\to X$ be counterclockwise rotation through angle $\theta$.
Then for every $x\in X$ and $\theta_1\not=\theta_2$, $[g_{\theta_1},x]\not= [g_{\theta_2},x]$. 
It follows that if ${\bf 0}$ is the origin in $\br^2$,
$\{ U(g_\theta,{\bf 0},1) ~|~ 0\leq\theta < 2\pi\}$ is an uncountable collection of mutually disjoint, nonempty open subsets
of $\ligx$. Hence, $\ligx$ is not second countable.
\end{example}

\begin{proposition}
\label{prop:ConjugateGroupoids}
If $X$ is a metric space with a subgroup $\Gamma\leq LS(X)$ and $h\in LS(X)$, then $\gammax$ and $\gammaxh$ are isomorphic topological
groupoids.
\end{proposition}

\begin{proof}
According to Proposition~\ref{prop:loc-sim},
$(h^{-1})_\ast\co\lsgx\to\lsgx$, $[g,x]\mapsto [h^{-1}gh, h^{-1}x]$, is an isomorphism of topological groups.
Clearly, $(h^{-1})_\ast$ takes the open subgroupoid $\gammax$ onto 
$\gammaxh$.
\end{proof}


\section{Ultrametric spaces}
\label{section:UltarametricSpaces}

In this section we recall the definition of ultrametric spaces and some of their well-known
properties. We then establish some elementary properties 
which have not appeared previously in the literature.
These properties  will be useful in studying 
local isometry and similarity groups and groupoids of ultrametric spaces.

\begin{definition}
\label{def:ultrametric}
If $(X,d)$ is a metric space and $d(x,y)\leq\max\{ d(x,z),d(z,y)\}$
for all $x,y,z\in X$, then $d$ is an {\it ultrametric\/} and $(X,d)$ is an {\it ultrametric
space\/}.
\end{definition}

The following proposition lists some well-known properties of ultrametric spaces. They are
readily verified.

\begin{proposition}[Elementary properties of ultrametric spaces]
\label{prop:ElemProp}
The
following properties hold in any ultrametric space $(X,d)$.
\begin{enumerate}
\item If two open balls (or two closed balls) in $X$ intersect, then one contains the other.
\item {\bf (Egocentricity)} Every point in an open (or closed) ball is a center of the ball.
\item Every open ball is closed, and every closed ball is open.
\item {\bf (ISB)} Every triangle in $X$ is isosceles with a short base. That is, if $x_1,x_2,x_3\in
X$, then there exists an $i$ such that $d(x_j,x_k)\leq d(x_i,x_j)=d(x_i,x_k)$ whenever 
$j\not= i\not= k$. {$\square$\smallskip}
\end{enumerate}
\end{proposition}

\begin{lemma}[Isometry Extension]
\label{lemma:IsometryExtension} 
Suppose $X$ is an ultrametric space, $x\in X$ and
$\epsilon >0$. If $h\co B(x,\epsilon)\to B(x,\epsilon)$ is an isometry, then 
$\tilde h\co X\to X$ defined by 
$$\tilde h= \left\{ \begin{array}{ll}
h & \textrm{on $B(x,\epsilon)$}\\
\textrm{inclusion} & \textrm{on $X\setminus B(x,\epsilon)$}\\
\end{array} \right.$$
is also an isometry.
\end{lemma}

\begin{proofnosquare} First observe that for all $y,z\in X$,
$y,z\in B(x,\epsilon)$ implies $d(y,z)<\epsilon$,
and $y\in B(x,\epsilon), z\notin B(x,\epsilon)$ implies $d(y,z)\geq\epsilon$.
[The first implication follows immediately from the ultrametric inequality. 
The second follows because $\epsilon\leq d(x,z)\leq\max\{ d(x,y),d(y,z)\}$, and
$d(x,y)<\epsilon$; thus, $\epsilon\leq d(y,z)$.]
Now to show that $\tilde h$ is an isometry, it suffices to let $x_1\in B(x,\epsilon)$,
$x_2\notin B(x,\epsilon)$ and show $d(x_1,x_2)=d(hx_1,x_2)$. For this note that on one hand,
$$d(x_1,x_2)\leq\max\{ d(x_1,hx_1),d(hx_1,x_2)\}=d(hx_1,x_2).$$
And on the other hand,
$$d(hx_1,x_2)\leq\max\{ d(hx_1,x_1),d(x_1,x_2)\}=d(x_1,x_2). ~~~\square$$
\end{proofnosquare}

\begin{remark} Lemma~\ref{lemma:IsometryExtension} need not hold for isometries $h\co B(x,\epsilon)\to B(y,\epsilon)$.
For example, the end space of the Fibonacci tree is rigid, but there are some local isometries
(see \cite{Hug}, Prop. 9.5).
\end{remark}

\begin{lemma}[Circular Equidistance]
\label{lemma:CircularEquidistance} 
If $(X,d)$ is an ultrametric space with points 
$w,x,y,z$ in $X$ such that $d(x,w)\not= d(x,y)= d(x,z)$, then $d(w,y)=d(w,z)$.
That is, if there exists a point $x\in X$ an equidistance $\ell$ to two points
$y,z$ then every other point $w$ whose distance from $x$ is different from $\ell$
is equidistant to $y$ and $z$. 
In yet other words, let $\ell>0$, $x\in X$ and consider the ``circle'' $C= \{ y\in X 
~|~ d(x,y)=\ell\}$. Then any point not on the circle $C$ is equidistant to any two points on
the circle $C$.
\end{lemma}

\begin{proof} Let 
$\ell=d(x,y)=d(x,z)$ and let $r=d(x,w)$. If $r<\ell$, then by the ISB property (Proposition~\ref{prop:ElemProp}), since
$d(x,y)>d(x,w)$, it must be the case that $d(y,w)=d(y,x)$. 
Likewise $d(z,w)=d(z,x)$. Hence, $d(z,w)=d(y,w)=\ell$.
If $r>\ell$, then by the ISB property, since $d(w,x)>d(y,x)$, it must be the case that
$d(w,y)=d(y,x)$. Likewise, $d(z,w)=d(z,x)$. Hence, $d(z,w)=d(y,w)=r$
\end{proof}

\begin{lemma}[Modification of Local Isometry]
\label{lemma:modification}
If $(X,d)$ is an ultrametric space,
$x\in X$, $\epsilon>0$ and $g\co B(x,\epsilon)\to B(x,\epsilon)$ is an isometry such that
$g(x)=x$ and $g$ is non-trivial arbitrarily close to $x$ (that is, for every $\delta>0$,
$\delta\leq\epsilon$, there exists $y\in B(x,\delta)$ such that $g(y)\not= y$), then
there exists an isometry $\tilde g\co B(x,\epsilon)\to B(x,\epsilon)$ such that $\tilde g(x)=x$,
$\tilde g$ is non-trivial arbitrarily close to $x$, and for every $\delta > 0$, $\delta\leq
\epsilon$, there exists $y\in B(x,\delta)$ and $\mu>0$ such that
$\tilde g|\co B(y,\mu)\to B(y,\mu)$ is the identity.
\end{lemma}

\begin{proof}
Choose a sequence $\{ x_i\}_{i=1}^\infty$ of distinct points in $B(x,\epsilon)$ converging to
$x$ such that
\begin{enumerate}
\item for every $i\in\bn$, $g(x_i)\not= x_i$, and
\item $d(x,x_1) > d(x,x_2) > d(x,x_3) > \cdots$.
\end{enumerate}
For each $i\in\bn$ let $C_i=\{ y\in X ~|~ d(x,y)=d(x,x_i)\}$ and note that $g(C_i)= C_i$.
Define $$\tilde g(x)= \left\{ \begin{array}{ll}
g(x) & \textrm{if $x\in\cup_{i=1}^\infty C_{2i}$}\\
x & \textrm{if $x\notin\cup_{i=1}^\infty C_{2i}$.}\\
\end{array} \right.$$
Note that the Circular Equidistance Lemma~\ref{lemma:CircularEquidistance} implies that $\tilde g\co B(x,\epsilon)\to B(x,\epsilon)$
is an isometry. The rest of the properties are straightforward to verify.
\end{proof}

\begin{lemma}[Local Isometry Extension]
\label{lemma:LIExtension}
If $(X,d)$ is an ultrametric space,
then $\ligx=\gammax$, where $\Gamma= LI(X)$.
\end{lemma}

\begin{proof}
Clearly $\gammax\subseteq\ligx$. 
Now let $g\co B(x,\epsilon)\to B(gx,\epsilon)$ be an isometry representing
$[g,x]\in\ligx$. Define $\tilde g\co X\to X$ by
$$\tilde g= \left\{ \begin{array}{ll}
g & \textrm{on $B(x,\epsilon)$}\\
g^{-1} & \textrm{on $B(gx,\epsilon)\setminus B(x,\epsilon)$}\\
\textrm{inclusion} & \textrm{on $X\setminus(B(x,\epsilon)\cup B(gx,\epsilon))$}\\
\end{array} \right.$$
It is easy to verify that $\tilde g$ is a local isometry
(recall that open balls are closed and $B(x,\epsilon)=B(gx,\epsilon)$ or
$B(x,\epsilon)\cap B(gx,\epsilon)=\emptyset$).
Thus,
$[g,x]=[\tilde g,x]\in\gammax$.
\end{proof}

A similar result does not hold for local similarities as the next example shows.

\begin{example} Let $X= \{ z_\infty, z_0, z_1, z_2,\dots\}$ with ultrametric $d$ given by
$$d(z_i,z_j)= e^{-\min\{ i,j\}} \text{ if } i\not= j.$$
Thus, $X$ is the end space of 
the Sturmian tree---see Example~\ref{ex:Sturmian Tree}.
Define $g\co X\to X$ by $gz_\infty=z_\infty$ and $gz_i=z_{i+1}$ for $i=0,1,2,\dots$.
Then $g|\co B(z_\infty,1)\to B(z_\infty, e^{-1})$ is an $e^{-1}$-similarity representing
$[g,z_\infty]\in\ligx$. However, there is no
local similarity $h\co X\to X$ with $[h,z_\infty]=[g,z_\infty]$.
Hence, $\gammax\subsetneqq\lsgx$,
where $\Gamma= LS(X)$.
\end{example}

The groupoids associated to a compact ultrametric space $X$ studied in this paper  are of 
two types. First, there is the full groupoid $\lsgx$ of local similarity germs on $X$.
Second, there are the groupoids of the form $\gammax$ where $\Gamma$ is a subgroup of
$LS(X)$.
By Lemma~\ref{lemma:LIExtension}, this second type includes $\ligx$. 
Moreover, by Remark~\ref{remark:opensubgrpoids},
the groupoids $\gammax$ are open subgroupoids of $\lsgx$.


\section{Recollections on trees and their ends}
\label{section:trees}

The material in this section is well-known; we collect it here for the convenience of the reader.
For more background on $\br$-trees, see Bestvina \cite{Bes}, Chiswell \cite{Chi}, and Morgan and Shalen \cite{MoS}.
For more information and references
on end spaces of $\br$-trees, see Hughes \cite{Hug} and
Mart\'inez-P\'erez and  Mor\'on
\cite{MPM}.


\subsection{Trees}
\label{subsection:trees}

An {\it $\br$-tree} is a metric space $(T,d)$ that is uniquely arcwise connected, and for any two 
points $x,y\in T$ the unique arc from $x$ to $y$, denoted $[x,y]$, is
isometric to the subinterval $[0, d(x,y)]$ of $\br$.

An $\br$-tree is {\it proper} if every closed metric ball in $T$ is compact.

As an example, let $T$ be a {\it locally finite simplicial tree}; that is, $T$ is the (geometric realization of)
a locally finite, one-dimensional, simply  connected, simplicial complex.
There is a natural unique metric $d$ on $T$ such that 
$(T,d)$ is an $\br$-tree, every edge is isometric
to the closed unit interval $[0,1]$, and the distance between distinct vertices $v_1, v_2$ is the minimum number
of edges in a sequence of edges $e_0, e_1,\dots, e_n$ with
$v_1\in e_0$, $v_2\in e_n$ and $e_i\cap e_{i+1}\not=\emptyset$ for $0\leq i\leq n-1$.
It follows that $(T,d)$ is a proper $\br$-tree.

Whenever we refer to a {\it locally finite simplicial tree} $T$, the metric $d$ on $T$ will be understood
to be the natural one just described.

Choose a {\it root} (i.e., a base vertex) $v\in T$. The rooted tree $(T,v)$ is {\it geodesically complete} if for every isometric
embedding $x\co [0,t]\to T$, $t> 0$, with $x(0)=v$, extends to an isometric embedding $\tilde x\co [0,\infty)\to T$. 
Such a map $\tilde f$ is a {\it geodesic ray in $T$ beginning at $v$}. In other words, $T$ is geodesically complete
if every vertex of $T$, except possibly the root, lies in at least two edges.

\subsection{Ends of trees}
\label{sec:ends}

The {\it end space} of a rooted $\br$-tree $(T,v)$ is given by
$$end(T,v) = \{ x\co [0,\infty)\to T ~|~ x(0)=v \text{ and } x \text{ is an isometric 
embedding}\}.$$
For $x,y\in end(T,v)$, define
$$d_e(x,y) = \left\{ \begin{array}{ll}
0 & \text{if $x=y$}\\
1/e^{t_0} &\text{if $x\not= y$ and $t_0=\sup\{ t\geq 0 ~|~ x(t)=y(t)\}$}. \end{array}\right.
$$

It follows that $(end(T,v),d_e)$ is a 
complete ultrametric space of diameter $\leq 1$.
The elements of $end(T,v)$ are called {\it ends} of $(T,v)$. 

\begin{proposition}
\label{prop:proper tree} Let $(T,v)$ be a geodesically complete, rooted $\br$-tree.
Then $T$ is proper if and only if $end(T,v)$ is compact.
\end{proposition}

\begin{proof}
First, assume $T$ is proper and show that $end(T,v)$ is totally bounded.
Let $\epsilon>0$ be given; to show that $end(T,v)$ can be covered by a finite number of
closed $\epsilon$-balls, we may assume $\epsilon<1$.
Let $r=-\ln\epsilon$.
Since $\overline B(v,r)$ is compact, so is $\partial\overline B(v,r)=\{ t\in T ~|~ d(t,v)=r\}$.
We claim that $\partial\overline B(v,r)$ is finite. On the contrary assume that there is an infinite
set $\{ t_i\}_{i=1}^\infty$ of distinct points in $\partial\overline B(v,r)$.
Choose $\{ x_i\}_{i=1}^\infty\subseteq end(T,v)$ such that $x_i(r)=t_i$ for all $i\geq 1$.
Then the sets $x_i([r,r+1])$, $i\geq 1$, are mutually disjoint. Hence, $d(x_i(r+1), x_j(r+1))\geq 2$ 
if $i\not= j$. This contradicts the compactness of $\partial\overline B(v,r+1)$.

Therefore, write $\partial\overline B(v,r)=\{ t_i\}_{i=1}^N$ and choose $\{ x_i\}_{i=1}^N\subseteq end(T,v)$
as above (so that $x_i(r)=t_i$).
Clearly, $end(T,v)= \cup_{i=1}^N\overline B(x_i,\epsilon)$.

Conversely, assume that $end(T,v)$ is compact, let $r>0$ be given, and show that
$\overline B(v,r)$ is compact in $T$ by showing every sequence in 
$\overline B(v,r)$ has a convergent subsequence. 
Let $\{ t_i\}_{i=1}^\infty$ be a sequence in $\overline B(v,r)$ and choose
$\{ x_i\}_{i=1}^\infty\subseteq end(T,v)$ 
such that $t_i=x_i(d(v,t_i))$ for all $i\geq 1$.
By passing to a subsequence, we may assume there exists $x_0\in end(T,v)$ such that
$x_i\to x_0$ in $end(T,v)$ as $i\to\infty$.
Hence, there exists $N$ such that $d_e(x_0,x_i)\leq e^{-r}$ for all $i\geq N$.
That is, $x_0(t)=x_i(t)$ if $0\leq t\leq r$ and $i\geq N$.
In particular, $t_i=x_0(d(v,t_i))$ for all $i\geq N$.
So $t_i$ is in the compact subset $x_0([0,r])$ of $\overline B(v,r)$ for all $i\geq N$.
Thus, $\{ t_i\}_{i=1}^\infty$ has a convergent subsequence.
\end{proof}

\begin{corollary} 
\label{cor:simplicial}
Let $(T,v)$ be a geodesically complete, rooted $\br$-tree.
$T$ is a locally finite simplicial tree if and only if
$end(T,v)$ is compact and
has distance set 
$$\{ t\in\br ~|~ \hbox{ there exists }
x,y \in end(T,v) \hbox{ such that } d_e(x,y) = t\}$$
contained in $\{ 0\} \cup \{ e^{-i} ~|~ i= 0, 1, 2, \dots\}$.
\end{corollary}

\begin{proof}
Necessity follows from Proposition~\ref{prop:proper tree} and
obvious facts about the metric $d_e$ when $T$ is simplicial.
Conversely, given the distance set condition, declare
all points of $T$ of the form $x(n)$ with $x\in end(T,v)$ and
$n\in\{ 0,1,2,\dots\}$ to be vertices; likewise, sets of the 
form $x([n,n+1])$ are edges. It is easily seen that this makes $T$ into
a simplicial tree. Compactness of $end(T,v)$ guarantees local finiteness.
\end{proof}

\begin{remark}
Let $(T,v)$ be a rooted, geodesically complete, locally finite simplicial tree.
The ends of $(T,v)$ are in one-to-one correspondence with
infinite sequences of distinct edges
$e_0, e_1, e_2, \dots$ such that $v\in e_0$ and for $i\geq 1$, $e_{i-1}\cap e_i$ consists of exactly one vertex, say $v_i$, 
and the vertices $v, v_1, v_2, \dots$ are distinct.
\end{remark}

\begin{remark}
One can verify that proper $\br$-trees are equivalent to the $\br$-trees
called {\it simplicial} 
in \cite{Bes} that are additionaly required to be locally finite.
\end{remark}

We include the following definition from \cite{Hug}.

\begin{definition} A {\it cut set\/} $C$ for a geodesically 
complete, rooted $\br$-tree $(T,v)$ is a subset $C$ of $T$ such that $v\notin C$ and for
every isometric embedding $\alpha\co [0,\infty)\to T$ with $\alpha(0)=v$
there exists a unique $t_0>0$ such that $\alpha(t_0)\in C$.
For $v\not= c\in T$, let $T_c$ denote the {\it subtree of $(T,v)$
descending from $c$}; that is,
$T_c=\{x\in T ~|~ c\in [v,x]\}$.
An {\it isometry at infinity\/} 
between geodesically 
complete, rooted $\br$-trees
$(T,v)$ and $(S,w)$ is a triple
$(f, C_T, C_S)$ where $C_T$ and $C_S$ are cut sets of $T$ and $S$,
respectively, and 
$f\co \cup\{ T_c ~|~ c\in C_T\} \to \cup\{ S_c ~|~ c\in C_S\}$
is a homeomorphism such that
\begin{enumerate}
\item $f(C_T)=C_S$, and
\item for every $c\in C_T$, $f|\co T_c\to S_{f(c)}$ is an isometry.
\end{enumerate}
An isometry at infinity $(f,C_T, C_S)\co (T,v)\to (S,w)$ 
is a {\it uniform isometry at infinity\/} provided
there exist $\epsilon, \delta >0$ such that $C_T=\partial B(v,\epsilon)$ 
and $C_S=\partial B(w,\delta)$.
\end{definition}

\paragraph{The end space functor.}
Let $\mathbf{U}_1$ be the category whose objects are 
compact ultrametric spaces of diameter less than or equal to $1$ and whose
morphisms are isometries.
Let 
$\left\{ \begin{array}{c}
\mathbf{U}_2\\ \mathbf{U}_3
\end{array}\right\}$
be the category whose objects are compact ultrametric spaces and whose
morphisms are
$\left\{ \begin{array}{c}
\text{uniform local similarities} \\ \text{local
similarities}
\end{array}\right\}$.
Let 
$\left\{ \begin{array}{c}
\mathbf{T}_1 \\ \mathbf{T}_2\\ \mathbf{T}_3
\end{array}\right\}$
be the category whose objects are proper, rooted, geodesically
complete $\br$-trees and whose
morphisms are
$$\left\{ \begin{array}{c}
\text{rooted isometries} \\ \text{equivalence classes of uniform 
isometries at infinity} \\ \text{equivalence classes 
of isometries at infinity}
\end{array}\right\}.$$
The equivalence classes just referred to are germs-at-infinity (see
\cite{Hug} for precise definitions).
The following result follows from Proposition~\ref{prop:proper tree} and
\cite{Hug}.

\begin{proposition}
\label{proposition:functor}
The end space functor $\mathcal E$ restricts to equivalences of
categories $\mathcal E\co \mathbf{T}_i\to\mathbf{U}_i$ for $i=1,2,3$.
\end{proposition}

\paragraph{Balls in the ends of simplicial trees.}
It will be convenient to have the following description of the metric balls in 
$end(T,v)$, where $(T,v)$ is a rooted, geodesically complete, locally finite simplicial tree.
For each $x\in end(T,v)$ and $0< \epsilon \leq 1$,
$$B(x,\epsilon) = \{ y\in end(T,v) ~|~ d_e(x,y)<\epsilon\}
= \{ y\in end(T,v) ~|~ -\ln\epsilon< t_0\}$$
and
$$\overline B(x,\epsilon) = \{ y\in end(T,v) ~|~ d_e(x,y)\leq\epsilon\}
= \{ y\in end(T,v) ~|~ -\ln\epsilon\leq t_0\},$$
where $t_0=\sup\{ t\geq 0 ~|~ x(t) = y(t)\}$.
Let $\lceil -\ln\epsilon\rceil$
be the smallest positive integer greater than or equal to $-\ln\epsilon$.
Then $x(\lceil -\ln\epsilon\rceil)$ is a vertex of $T$ that we denote by
$v_{\{x,\epsilon\}}$.
Let $T_{\{x,\epsilon\}}$ denote the subtree of $T$ descending from
$v_{\{x,\epsilon\}}$; i.e.,
$$T_{\{x,\epsilon\}}=\bigcup\{ y(t) ~|~ y\in end(T,v), y(\lceil-\ln\epsilon\rceil)= v_{\{x,\epsilon\}}, \text{~and
$t\geq\lceil-\ln\epsilon\rceil$}\}.$$
Then $(T_{\{x,\epsilon\}}, v_{\{x,\epsilon\}})$ is itself a rooted, geodesically complete, 
locally finite simplicial tree.
We make the identification
$$\overline B(x,\epsilon) = end(T_{\{x,\epsilon\}}, v_{\{x,\epsilon\}}),$$
where $y\in \overline B(x,\epsilon)$ is identified with 
$\tilde y\in end(T_{\{x,\epsilon\}},v_{\{x,\epsilon\}})$
defined by $\tilde y(t) = y(t+\lceil-\ln\epsilon\rceil)$ for $t\geq 0$. Conversely, of course,
$\tilde y\in end(T_{\{x,\epsilon\}},v_{\{x,\epsilon\}})$ is identified
with $y\in \overline B(x,\epsilon)$ defined by
$$y(t) = \left\{ \begin{array}{ll}
x(t) & \text{for $0\leq t\leq\lceil-\ln\epsilon\rceil$}\\
\tilde y(t-\lceil-\ln\epsilon\rceil) &\text{for $t\geq\lceil-\ln\epsilon\rceil$}. \end{array}\right.
$$

Likewise, 
let $\lceil\lceil -\ln\epsilon\rceil\rceil$
be the smallest positive integer greater than $-\ln\epsilon$.
Thus, 
$\lceil -\ln\epsilon\rceil \leq \lceil\lceil -\ln\epsilon\rceil\rceil$,
with equality if and only if $-\ln\epsilon$ is an integer.
Then $x(\lceil\lceil -\ln\epsilon\rceil\rceil)$ is a vertex of $T$ that we denote by
$v_{\langle x,\epsilon\rangle}$.
Let $T_{\langle x,\epsilon\rangle}$ denote the subtree of $T$ descending from
$v_{\langle x,\epsilon\rangle}$; i.e.,
$$T_{\langle x,\epsilon\rangle}=\bigcup\{ y(t) ~|~ y\in end(T,v), y(\lceil\lceil-\ln\epsilon\rceil\rceil)= 
v_{\langle x,\epsilon\rangle}, \text{~and
$t\geq\lceil\lceil-\ln\epsilon\rceil\rceil$}\}.$$
Then $(T_{\langle x,\epsilon\rangle}, v_{\langle x,\epsilon\rangle})$ 
is itself a rooted, geodesically complete, simplicial tree.
We make the identification
$$B(x,\epsilon) = end(T_{\langle x,\epsilon\rangle}, v_{\langle x,\epsilon\rangle}),$$

\subsection{Examples}
\label{subsection:tree examples}

In this section we give examples of a few trees and their end spaces.
These examples appear again in Section~\ref{section:TheSymmetryAtInfinityGroup}.

\begin{example}{\bf The Cantor tree $C$.}
\label{ex:Cantor Tree}
The Cantor tree $C$, also called the infinite binary tree, is a locally finite
simplicial tree. It has a root $v$ of 
valency two (i.e., there exists exactly two edges containing $v$) and every other vertex is
of valency three. 
If $w$ is a vertex different from $v$, then the two edges that contain $w$ and are separated 
from $v$ by $w$ are not labelled identically.
Each edge is labelled $0$ or $1$ so that for every vertex $w$, at least one
edge containing $w$ is labelled $0$ and at least one is labelled $1$.

Let $end(C)=end(C,v)$ since the root $v$ is understood. An element of $end(C)$, being an infinite
sequence of successively adjacent edges in $C$ beginning at $v$, can be labelled uniquely by
an infinite sequence of $0$'s and $1$'s. Thus,
$$end(C) =\{ (x_0,x_1,x_2,\dots) ~|~ x_i\in\{ 0,1\} \text{ for each $i$}\}$$
and 
$$d_e((x_i),(y_i))=
\left\{ \begin{array}{ll}
0 & \mbox{if $(x_i)=(y_i)$}\\
1/e^{n} & \mbox{if $(x_i)\not= (y_i)$ and $n=\inf\{ i\geq 0 ~|~ x_i\not= y_i\}$ .}
\end{array}\right. $$ 
\end{example}

\begin{example}{\bf The Fibonacci tree $F$.}
\label{ex:Fibonacci Tree}
The Fibonacci tree $F$ is a subtree of $C$ with the same root $v$ and labelling scheme. 
In $F$, only edges labelled $0$ are allowed to follow edges labelled $1$ as one moves away
from the root. Thus,
$$end(F) =\left\{ (x_0,x_1,x_2,\dots)\in end(C) ~|~ x_i=1 \text{ implies $x_{i+1}=0$}\right\}.$$
See \cite{Hug} for some compaisons of the 
Cantor and Fibonacci trees. 
\end{example}

\begin{example}{\bf The Sturmian tree $S$.}
\label{ex:Sturmian Tree}
The Sturmian tree $S$ is also a subtree of $C$ with the same root $v$ and labelling scheme. 
In $S$, only edges labelled $1$ are allowed to follow edges labelled $1$ as one moves away
from the root. Thus,
$$end(S) =\{ (x_0,x_1,x_2,\dots)\in end(C) ~|~ x_i=1 \mbox{ implies $x_{i+1}=1$}\}.$$
In particular, $end(S)$ is countably infinite: $end(S) = \{z_\infty, z_0, z_1, z_2,\dots\}$, where
$
z_\infty  =  (0,0,0,\dots),
z_0  =  (1,1,1,\dots),
z_1  =  (0,1,1,1,\dots),
z_2  =  (0,0,1,1,\dots), \dots .
$
The metric is given by
$d(z_i,z_j)= e^{-\min\{ i,j\}} \text{ if } i\not= j.$
\end{example}
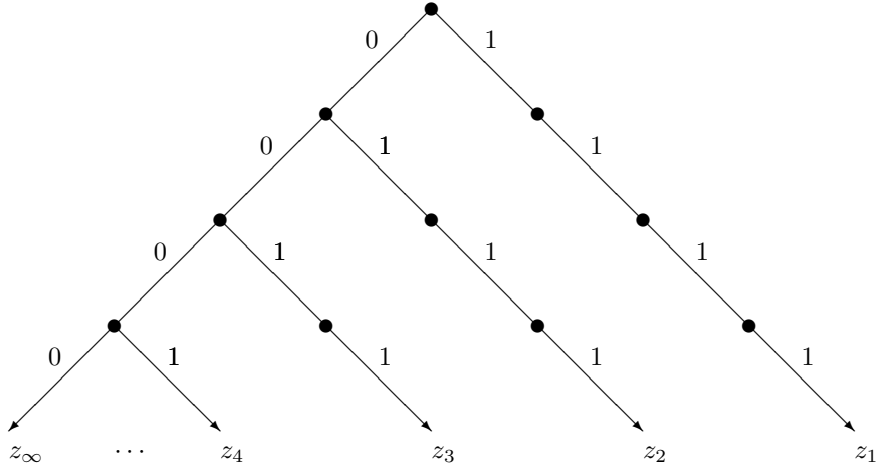
\begin{figure}[htbp]


\begin{picture}(400,160)(-20,-10)


\put(0,0){\line(1,1){160}}
\put(0,0){\vector(-1,-1){0}}
\put(160,160){\circle*{5}}

\put(40,40){\line(1,-1){40}}
\put(80,0){\vector(1,-1){0}}
\multiput(40,40)(80,0){4}{\circle*{5}}

\put(80,80){\line(1,-1){80}}
\put(160,0){\vector(1,-1){0}}
\multiput(80,80)(80,0){3}{\circle*{5}}

\put(120,120){\line(1,-1){120}}
\put(240,0){\vector(1,-1){0}}
\multiput(120,120)(80,0){2}{\circle*{5}}

\put(160,160){\line(1,-1){160}}
\put(320,0){\vector(1,-1){0}}

\put(0,-10){$z_\infty$}
\put(320,-10){$z_1$} \put(240,-10){$z_2$} \put(160,-10){$z_3$} \put(80,-10){$z_4$} \put(40,-10){$\cdots$}

\multiput(15,25)(40,40){4}{$0$}
\multiput(60,25)(40,40){4}{$1$}
\multiput(60,25)(80,0){4}{$1$}
\multiput(100,65)(80,0){3}{$1$}
\multiput(140,105)(80,0){2}{$1$}

\end{picture}
\caption{The Sturmian tree}
\label{fig:sturm-tree}
\end{figure}

\begin{example}{\bf The $n$-regular tree $R_n$.}
For $n=1,2,3,\dots$, the $n$-regular tree $R_n$ is 
the simplicial tree such that every vertex has valency $n+1$.
It is homogeneous so that a root can be chosen arbitrarily.
It is geodesically complete and locally finite.
The edges can be labelled by the integers $0,1,\dots,n$ so that 
for each vertex each label appears on exactly one edge 
containing the vertex. 
Thus,
$$end(R_n) =\{ (x_0,x_1,x_2,\dots) ~|~ x_i\in\{ 0,1,\dots,n\} \text{ and 
$x_{i+1}\not= x_i$ for each $i$}\}.$$
\end{example}


\begin{example}{\bf The infinite $n$-ary tree $A_n$.}
For $n=1,2,3,\dots$, the infinite $n$-ary tree $A_n$ is 
the rooted,
geodesically complete, locally finite simplicial tree such that
every vertex except the root has 
valency $n+1$ , and the root has valency $n$.
For example, $A_2$ is the Cantor tree.
\end{example}

\begin{example}{\bf The $n$-ended tree $E_n$.} 
For $n=1,2,3,\dots$, $(E_n,v)$ is 
the simplicial tree such that the root $v$ has valency $n$ and
all other vertices have valency $2$. Thus, $end(E_n)$ consists of $n$
points, each a distance $1$ from any other.
\end{example}

\begin{example}{\bf The irrational tree $T_\alpha$.} 
Let $\alpha$ be  positive irrational number and $\alpha=
[a_0,a_1,a_2,\dots]$ its continued fraction expansion.
Thus, $\{ a_i\}_{i=0}^\infty$ is a sequence of non-negative
integers such that $a_i\geq 1$ if $i\geq 1$ and
$$\alpha ~=~ a_0 + \cfrac{1}{a_1 +\cfrac{1}{a_2 +\cfrac{1}{\ddots}}}$$
Consider the compact metric space
$$X_\alpha = 
\left\{ (x_0,x_1,x_2,\dots) \in  
\prod_0^\infty\bz_+ ~|~ \text{$0\leq x_i\leq a_i$ and 
$x_i=a_i$ implies $x_{i+1}=0$} \right\}$$
with ultrametric
$$d_e((x_i),(y_i))=
\left\{ \begin{array}{ll}
0 & \mbox{if $(x_i)=(y_i)$}\\
1/e^{n} & \mbox{if $(x_i)\not= (y_i)$ and 
$n=\inf\{ i\geq 0 ~|~ x_i\not= y_i\}$ .}
\end{array}\right. $$ 
Let $(T_\alpha, v)$ be the rooted, geodesically complete, locally finite
simplicial tree such that $end(T_\alpha, v)$ is isometric to $X_\alpha$.
For example, the golden mean ${\frac {1+\sqrt{5}} 2} = [1,1,1,\dots]$
and $T_{\frac {1+\sqrt{5}} 2} = F$, the Fibonacci tree.
The spaces $X_\alpha$ appear in Mingo \cite{Min}.
\end{example}


\section{Local rigidity, locally rigid actions, and Hausdorffness}
\label{section:LR,LRA,H}

\subsection{Locally rigid ultrametric spaces}
\label{subsection:LocallyRigidUltrametricSpaces}

The main goal of this section is to characterize when the groupoid of local isometries on an ultrametric space is
Hausdorff. The answer is in terms of a local rigidity condition on the ultrametric space.
We also discuss the second countability of the local isometry groupoid of
locally rigid ultrametric spaces.

\begin{definition}
\label{def:locally-rigid}
A metric space $(X,d)$ is {\it locally rigid\/} 
if for every $x\in X$ there exists $\epsilon_x>0$ such that any isometry
$h\co B(x,\epsilon_x)\to B(x,\epsilon_x)$ is the identity.
\end{definition}

\begin{lemma}
\label{lemma:locally-rigid}
An ultrametric space $X$ is locally rigid
if and only if for every $x\in X$ there exists $\epsilon_x>0$ such that for any
$0<\epsilon\leq\epsilon_x$, every isometry
$h\co B(x,\epsilon)\to B(x,\epsilon)$ is the identity.
\end{lemma}

\begin{proof} This follows immediately from Lemma~\ref{lemma:IsometryExtension}.
\end{proof}

\begin{theorem}
\label{theorem:Hausdorff-groupoid}
For an ultrametric space $(X,d)$ the following are equivalent.
\begin{enumerate}
\item For every $x\in X$, $\epsilon> 0$ and isometry $g\co B(x,\epsilon)\to B(x,\epsilon)$
such that $g(x)=x$ there exists $\delta=\delta(\epsilon,x,g)>0$ 
such that $g|\co B(x,\delta)\to B(x,\delta)$ is the 
identity.
\item
\label{theorem:Hausdorff-groupoid:2}
For every $x\in X$ there exists $\epsilon_x>0$ such that if $g\co B(x,\epsilon_x)\to 
B(x,\epsilon_x)$ is an isometry with $g(x)=x$, then $g$ is the identity.
\item $X$ is locally rigid.
\item $\ligx$ is Hausdorff.
\end{enumerate}
\end{theorem}

\begin{proof} (1) implies (2). 
Suppose on the contrary that $X$ satisfies (1) but not (2).
Then there is a sequence of ``circles'' $C_i$ (in the sense of Lemma~\ref{lemma:CircularEquidistance}), 
$i\in\bn$, about 
some $x\in X$ 
of decreasing diameter and non-trivial isometries $g_i\co C_i\to C_i$. These can be pieced
together to give a non-trivial isometry $g\co B(x,\epsilon)\to B(x,\epsilon)$ which is non-trivial
on each ball about $x$.

(2) implies (3).
Suppose on the contrary that there exists $x\in X$ without the property in Definition~\ref{def:locally-rigid}.
Property (2) implies that there exists $\epsilon_0>0$ such that if $0<\epsilon\leq
\epsilon_0$ and $g\co B(x,\epsilon)\to B(x,\epsilon)$ is an isometry with $gx=x$, then $g$ is
the identity (this uses Lemma~\ref{lemma:IsometryExtension}
in a manner similar to how it is used in Lemma~\ref{lemma:locally-rigid}).
Choose $0< \epsilon_2 < \epsilon_1 < \epsilon_0$ 
so that for $i=1,2$ there exists an
isometry $h_i\co X\to X$ such that $h_iB(x,\epsilon_i)=B(x,\epsilon_i)$,
$h_i|(X\setminus B(x,\epsilon_i))=\textrm{inclusion}$ (this uses Lemma~\ref{lemma:IsometryExtension} again) and
$h_i$ is not the identity. It follows that $h_ix\not=x$ for $i=1,2$.
By choosing $\epsilon_1$ and $h_1$ before $\epsilon_2$ and $h_2$, we may assume that
$\epsilon_2<d(x,h_1x)$.
Consider the composition
$$g\co B(h_1x,\epsilon_2) \stackrel{h_1^{-1}}{\longrightarrow}
B(x,\epsilon_2) \stackrel{h_2}{\longrightarrow}
B(x,\epsilon_2) \stackrel{h_1|}{\longrightarrow} B(h_1x,\epsilon_2).$$
This isometry can be extended (by Lemma~\ref{lemma:IsometryExtension}) to an isometry $\tilde g\co X\to X$
such that $\tilde g|(X\setminus B(h_1x,\epsilon_2))$ is the inclusion. 
Now $\tilde gB(x,\epsilon_1)=B(x,\epsilon_1)$ and $\tilde gx=x$ (because
$x\notin B(h_1x,\epsilon_2)$. Thus, $\tilde g$ is the identity.
Since $\tilde gh_1x= h_1h_2x$ we have $h_1x=h_1h_2x$ and $x=h_2x$, 
a contradiction.

(3) implies (4).
Let $[g_1,x_1]\not= [g_2,x_2]$ in $\ligx$. If $d(x_1,x_2)=\epsilon >0$, then
$U(g_1,x_1,\epsilon)\cap U(g_2,x_2,\epsilon)=\emptyset$. 

If $d(g_1x_1,g_2x_2)=\epsilon >0$, then choose $\epsilon_i >0$ such that
$g_iB(x_i,\epsilon_i)\subseteq B(g_ix_i,\epsilon)$ for $i=1,2$ and observe that
$U(g_1,x_1,\epsilon_1)\cap U(g_2,x_2,\epsilon_2)=\emptyset$.
Finally suppose $x_1=x_2$ and $g_1x_1=g_2x_2$. Choose $\epsilon >0$ so that $g_i$ is defined
on $B(x_i,\epsilon)$ for $i=1,2$ and so that $\epsilon\leq\epsilon_x$ 
where $\epsilon_x > 0$ comes from Lemma~\ref{lemma:locally-rigid}.
Then $h=g_2^{-1}g_1\co B(x_1,\epsilon)\to B(x_1,\epsilon)$ is an isometry
so $h$ is the identity.
Hence $[g_1,x_1]=[g_2,x_2]$.

(4) implies (1).
Let $x\in X$, $\epsilon>0$ and $g\co B(x,\epsilon)\to B(x,\epsilon)$ be an isometry such that
$gx=x$. Suppose on the contrary that $g$ 
does not equal the identity on a sufficiently small ball
about $x$. Lemma~\ref{lemma:modification} gives another isometry $\tilde g\co B(x,\epsilon)\to B(x,\epsilon)$ such that
$\tilde gx=x$, $\tilde g$ is non-trivial arbitrarily close to $x$ and there exist points $y$
arbitrarily close to $x$ such that $\tilde g$ is the identity on sufficiently small balls 
about $y$. It follows that $U(\tilde g, x,\epsilon_1)\cap U(\id,x,\epsilon_2)\not=\emptyset$
for all $\epsilon_1,\epsilon_2 >0$, contradicting the Hausdorff property of $\ligx$.
\end{proof}

\begin{remark}\label{remark:locally-rigid-unit-space}
If $X$ is a locally rigid ultrametric space, then $X$ is an open and closed subset of
$\ligx$. To see this, recall that the embedding $\alpha\co X\to\ligx$ is given by
$\alpha(x)=[\id,x]$. If $[g,x]\in\ligx$ is not in the image of $\alpha$, then 
$[g,x]\not=[\id,x]$. Local rigidity, in particular Theorem~\ref{theorem:Hausdorff-groupoid}~(2), implies
that $gx\not= x$. If $0 <\epsilon <{\frac 1 2} d(x,gx)$, then
$U(g,x,\epsilon)\cap \alpha(X)=\emptyset$.
This shows that $\alpha(X)$ is closed in $\ligx$. It is open by Remark~\ref{remark:unit-space}.
\end{remark}

\begin{example}
\label{example:locally-rigid}
\begin{enumerate}
\item The end spaces of the following trees are not locally
rigid: the Cantor tree $C$, the $n$ regular tree $R_n$ ($n\geq 2)$,
and the $n$-ary tree $A_n$ ($n\geq 2$).
\item The end spaces of the following trees are locally
rigid:  the Fibonacci tree $F$ and the Sturmian tree $S$.
\item The end space of the irrational tree $T_\alpha$ is locally rigid if 
and only if $\alpha$ is equivalent to the golden mean $\frac {1+\sqrt{5}}
2$ under the action of $SL(2,\bz)$ by fractional linear transformations
(because this condition is equivalent to the continued fraction
expansion of $\alpha$ eventually ending in all $1$'s).
\end{enumerate}
\end{example}

\begin{lemma} 
\label{lemma:compact-locally-rigid}
If $X$ is a compact ultrametric space, then $X$ is locally rigid if and only if
there exists
$\epsilon_X>0$ such that for any $0<\epsilon\leq\epsilon_X$ and $x\in X$, every isometry
$h\co B(x,\epsilon)\to B(x,\epsilon)$ is the identity.
\end{lemma}

\begin{proof} Assume $X$ is locally rigid.
For each $x\in X$ let $\epsilon_x>0$ be given by Definition~\ref{def:locally-rigid}
and let $\epsilon_X$ be a
Lebesgue number for $\{ B(x,\epsilon_x) ~|~ x\in X\}$. Then if $0<\epsilon\leq\epsilon_X$
and $x\in X$, there exists $y\in X$ such that
$B(x,\epsilon)\subseteq B(x,\epsilon_X)\subseteq B(y,\epsilon_y)$.
Thus, $\epsilon\leq\epsilon_y$ and $B(x,\epsilon)=B(y,\epsilon)$, so any isometry
$h\co B(x,\epsilon)\to B(x,\epsilon)$ is the identity.

The converse is obvious.
\end{proof}

\begin{proposition}\label{prop:seccount}
If $X$ is a second countable, locally rigid ultrametric space, then $\ligx$ is second countable.
\end{proposition}

\begin{proof}
Let $\{ x_i\}_{i=1}^\infty$ be a countable dense subset of $X$. For each $i$, choose the least 
positive integer $i_0$ such that if $j\geq i_0$ and $g\co B(x_i, 1/j)\to B(x_i,1/j)$ is an isometry, then 
$g=\id_{B(x_i,1/j)}$.
Now suppose $j\geq i_0$ and $g,h\co B(x_i,1/j)\to X$ are two different isometric embeddings. Then 
$B(gx_i,1/j)\not= B(hx_i,1/j)$ (for otherwise 
$g^{-1}h\co B(x_i, 1/j)\to B(x_i, 1/j)$ would be a nontrivial isometry) and, hence,
$B(gx_i, 1/j)\cap B(hx_i, 1/j)=\emptyset$. It follows that for each $i$ and each $j\geq i_0$ there
are at most countably many distinct isometric embeddings, say
$g_{(i,j,k)}\co B(x_i, 1/j)\to X$, $1\leq k < N_{(i,j)}$, where $N_{(i,j)}$ is either a positive integer or $\infty$.
The proof will be complete once we show that
$$\mathcal{B} =
\{ U(g_{(i,j,k)}, x_i, 1/j) ~|~ 1\leq i <\infty, ~i_0\leq j <\infty, ~1\leq k < N_{(i,j)}\}$$
is a countable basis for $\ligx$. 
Given $U(g,x,\epsilon)$, we show that $U(g,x,\epsilon)$ is a union of elements of $\mathcal{B}$.
If $y\in B(x,\epsilon)$, let $0 < 1/n\leq \epsilon$ be chosen such that any self-isometry of
$B(y,1/n)$ is the identity. Then there exists $x_i\in B(y,1/n)$. Thus,
$B(x_i,1/n)= B(y,1/n)$ and by the choice of $i_0$, $i_0\leq n$.
It follows that
$U(g,x_i,1/n)\in\mathcal{B}$ and
$[g,y]\in U(g,x_i,1/n)\subseteq U(g,x,\epsilon)$.
Finally, note that if
$U(g_{(i,j,k)}, x_i, 1/j), U(g_{(i', j', k')}, x_{i'}, 1/j')\in\mathcal{B}$,
$j'\geq j$ and their intersection is nonempty, then
$B(x_{i'}, 1/j')\subseteq B(x_i, 1/j)$ and
$g_{(i', j', k')} = g_{(i,j,k)}|B(x_{i'}, 1/j')$
(because they must agree somewhere, hence, they agree everywhere on their common domain).
Hence, the intersection is
$U(g_{(i', j', k')}, x_{i'}, 1/j')$.
\end{proof}

The following two corollaries follow from Theorem~\ref{theorem:Hausdorff-groupoid}, Proposition~\ref{prop:seccount},
Lemmas~\ref{lemma:etale} and \ref{lemma:loccom}, and Remark~\ref{remark:opensubgrpoids}. 

\begin{corollary} If $X$ is a locally compact, second countable, locally rigid ultrametric space,
then $\ligx$ is a locally compact, second countable, Hausdorff \'etale groupoid.
\end{corollary}

\begin{corollary}
\label{corollary:groupoid-properties}
If $X$ is a compact, locally rigid ultrametric space,
then $\ligx$ is a locally compact, second countable, Hausdorff \'etale groupoid.
\end{corollary}

Finally, we establish the following two results that complete the proof
of Theorem~\ref{theorem:ls}(1) of the introduction.

\begin{lemma}
\label{lem:ls-seccount}
If $X$ is a compact, locally rigid ultrametric space, then $\lsgx$ is second countable.
\end{lemma}

\begin{proof}
Use the first part of the proof of Proposition~\ref{prop:seccount}
to find a countable basis $\{ B_i\}_{i=1}^\infty$ of $X$ by open balls,
each of which admits only countably many distinct isometric embeddings
into $X$.
Assume that for every $\epsilon>0$, there are only finitely many $i$'s
with $diam(B_i)>\epsilon$.
Because the distance set of $X$ is countable (see \cite{Ber} and the
proof of Proposition~\ref{prop:scale}), there exists a sequence
$\{\lambda_j\}_{j=1}^\infty$ of positive numbers such that if 
$g\co B_i\to g(B_i)$ is a similarity onto some open ball in $X$ and $x\in B_i$, 
then
$sim(g,x)=\lambda_j$ for some $j$; i.e., $g$ is a $\lambda_j$-similarity
(the $\lambda_j$'s are all ratios of distances in $X$).

Now if $g\co B_i\to g(B_i)$ and $h\co B_i\to h(B_i)$ are two
$\lambda_j$-similarities such that $g(B_i)\cap h(B_i)\not=\emptyset$,
then $g(B_i)=h(B_i)$ and $h^{-1}g\co B_i\to B_i$ is an isometry.
If the radius of $B_i$ is sufficiently small, then local rigidity implies
$g=h$. 
Hence, there exist only countably many distinct similarities of $B_i$
onto open balls of $X$, say $g_{i,k}$, where $1\leq k< N(i)$ and $N(i)\leq
\infty$.

Choose $x_i\in B_i$ for all $i$.
The proof will be complete once we show that
$$\mathcal{B} =
\{ U(g_{(i,k)}, x_i, diam(B_i)) ~|~ 1\leq i <\infty,  ~1\leq k < N(i)\}$$
is a countable basis for $\lsgx$. 
Given a basis element $U(g,x,\epsilon)$ of $\lsgx$,
$B(x,\epsilon)$ can be written as a union of $B_i$'s.
If $B_i\subseteq B(x,\epsilon)$, then $g|B_i = g_{i,k}$ for some $k$.
It follows that $U(g,x,\epsilon)$ 
is a union of elements of $\mathcal{B}$.
\end{proof}

\begin{lemma}
\label{lem:ls-hausdorff}
If $X$ is a compact, locally rigid ultrametric space, then $\lsgx$ is Hausdorff.
\end{lemma}

\begin{proof} This is very similar to the proof of 
Theorem~\ref{theorem:Hausdorff-groupoid} above and
Theorem~\ref{theorem:Hausdorffness} below.
Let $[g_1,x_1]\not= [g_2,x_2]$ in $\lsgx$. 
It is easy to reduce to the case that $x_1=x_2$
and $g_1x_1=g_2x_2$.
If $sim(g_1,x_1)=sim(g_2,x_2)$, then
$g_2^{-1}g_1\co B(x_1,\epsilon)\to B(x_1,\epsilon)$ is an isometry
for some $\epsilon >0$.
Local rigidity implies that $g_2^{-1}g_1$ is the identity; hence,
$[g_1,x_1]=[g_2,x_2]$.
If $sim(g_1,x_1)\not=sim(g_2,x_2)$, then 
$U(g_1,x_1,\epsilon)\cap U(g_2,x_2,\epsilon)=\emptyset$
for some sufficiently small $\epsilon>0$. 
\end{proof}

\begin{corollary} 
\label{cor:ls-grpd-props}
If $X$ is a compact, locally rigid ultrametric space, then  
$\lsgx$ is a locally compact, Hausdorff, second countable, 
\'etale groupoid.
\end{corollary}

\begin{proof}
This follows from Lemmas~\ref{lemma:etale},  \ref{lemma:loccom},
\ref{lem:ls-seccount}, and \ref{lem:ls-hausdorff}.
\end{proof}


\subsection{Locally rigid actions}
\label{subsection:LocallyRigidActions}

\begin{definition} Let $X$ be a metric space with a subgroup $\Gamma\leq LS(X)$. The action of $\Gamma$ on
$X$ is {\it locally rigid} if for every $x\in X$ and for every $g\in\Gamma_x$ such that $sim(g,x)=1$, there exists
$\epsilon > 0$ such that $g\in\Gamma_y$ for all $y\in B(x,\epsilon)$.
\end{definition}

Note that if $\Gamma$ acts locally rigidly on $X$ and $H$ is a subgroup of $\Gamma$, then $H$ also acts locally 
rigidly on $X$.

\begin{lemma}
Let $X$ be a metric space with a subgroup $\Gamma\leq LS(X)$. The following are equivalent:
\begin{enumerate}
\item The action of $\Gamma$ on $X$ is locally rigid.
\item For every $x\in X$ and for every $g,h\in\Gamma_x$ such that $sim(g,x)=sim(h,x)$, there exists
$\epsilon > 0$ such that $gy=hy$ for every $y\in B(x,\epsilon)$.
\item For every $x\in X$ and for every $g,h\in\Gamma$ such that $gx=hx$ and $sim(g,x)=sim(h,x)$, there exists
$\epsilon > 0$ such that $gy=hy$ for every $y\in B(x,\epsilon)$.
\end{enumerate}
\end{lemma}

\begin{proof}
(1) implies (2): Let $x\in X$ and $g,h\in\Gamma_x$ such that $sim(g,x)=sim(h,x)$ be given.
Then $h^{-1}g\in\Gamma_x$ and $sim(h^{-1}g,x)=1$.
Since the action of $\Gamma$ on $X$ is assumed to be locally rigid,
there exists
$\epsilon > 0$ such that $h^{-1}g\in\Gamma_y$ for all $y\in B(x,\epsilon)$; i.e., $gy=hy$ for all $y\in B(x,\epsilon)$.

(2) implies (3): Let $x\in X$ and $g,h\in\Gamma$ such that $gx=hx$ and $sim(g,x)=sim(h,x)$ be given.
Then $h^{-1}g\in\Gamma_x$ and $sim(h^{-1}g,x)=1= sim(\id_X,x)$.
Hence, there exists
$\epsilon > 0$ such that $h^{-1}gy=y$ for all $y\in B(x,\epsilon)$; i.e., $gy=hy$ for all $y\in B(x,\epsilon)$.

(3) implies (1): Let $x\in X$ and $g\in\Gamma_x$ such that $sim(g,x)=1$ be given. 
Since $sim(g,x)=sim(\id_X,x)$,
the result is obvious.
\end{proof}

\begin{theorem} Let $X$ be an ultrametric space.
\label{theorem:Hausdorffness}
\begin{enumerate}
\item $\lsgx$ is Hausdorff if and only if
for every $x\in X$ and for every $[g,x], [h,x]\in\lsgx$ such that $gx=hx$ and $sim(g,x)=sim(h,x)$, it follows that
$[g,x]=[h,x]$.
\item If $\Gamma\leq LS(X)$ and $\Gamma$ acts locally rigidly on $X$, then $\gammax$ is Hausdorff.
\label{hlr}
\end{enumerate}
\end{theorem}

\begin{proof}
(1) Assume first that $\lsgx$ is Hausdorff, and let 
$x\in X$ and $[g,x], [h,x]\in\lsgx$ such that $gx=hx$ and $sim(g,x)=sim(h,x)$ be given.
Suppose on the contrary that
$[g,x]\not=[h,x]$.
Choose $\epsilon>0$ so that $h^{-1}g$ is an isometry on $B(x,\epsilon)$.
Since $[h^{-1}g,x]\not= [\id_X,x]$, $h^{-1}g$ is
non-trivial arbitrarily close.
Hence, Lemma~\ref{lemma:modification} implies that there  
exists an isometry $\tilde g\co B(x,\epsilon)\to B(x,\epsilon)$ such that $\tilde g(x)\in\Gamma_x$,
$\tilde g$ is non-trivial arbitrarily close to $x$, and for every $\delta > 0$, $\delta\leq
\epsilon$, there exists $y\in B(x,\delta)$ and $\mu>0$ such that
$\tilde g|\co B(y,\mu)\to B(y,\mu)$ is the identity.
It follows that $[\tilde g, x]$ and $\id_X,x]$ can not be separated by open sets in $\lsgx$.

Conversely, let $[g,x], [h,y]\in\lsgx$ be given.
If $x\not= y$, choose $\epsilon > 0$ with $\epsilon\leq d(x,y)$; it is easy to see that 
$U(gx,\epsilon)\cap U(h,y,\epsilon)=\emptyset$.
If $gx\not= hy$, choose $\epsilon >0$ such that $g(B(x,\epsilon))\cap h(B(y,\epsilon))=\emptyset$;
it is easy to see that 
$U(gx,\epsilon)\cap U(h,y,\epsilon)=\emptyset$.
If $sim(g,x)\not= sim(h,y)$ choose $\epsilon > 0$ such that 
$sim(g,z)=sim(g,x)$ for all $z\in B(x,\epsilon)$ and $sim(h,z)=sim(h,y)$ for all $z\in B(y,\epsilon)$;
it is easy to see that 
$U(gx,\epsilon)\cap U(h,y,\epsilon)=\emptyset$.
[In each of these three cases, $\epsilon$ must be chosen so small that the germs are represented on 
$\epsilon$-balls.]
Thus, we are left with the case that
$x=y$, $gx=hx$ and $sim(g,x)=sim(h,x)$. Since the assumption in this case is that
$[g,x]=[h,x]$, there is nothing to separate.

(2) Let $[g,x], [h,y]\in\gammax$ be given. As in the  proof just given, it is easy to reduce to the case that
$x=y$, $gx=hx$ and $sim(g,x)=sim(h,x)$.
The assumption that $\Gamma$ is acting locally rigidly implies $h^{-1}g=\id$ near $x$; that is, $[g,x]=[h,x]$.
\end{proof}

The converse of Theorem~\ref{theorem:Hausdorffness}(\ref{hlr})
need not hold as the next example shows.

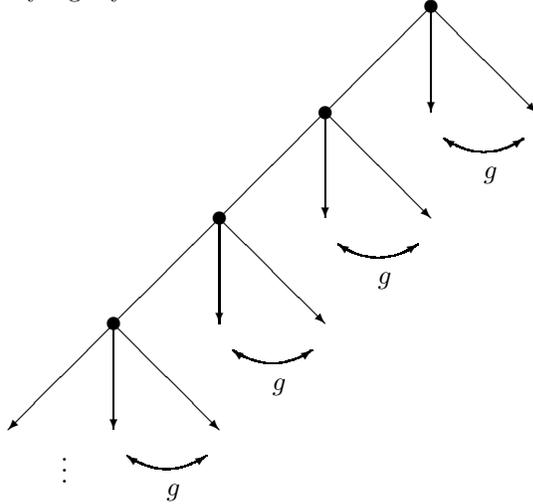
\begin{figure}[htbp]

\begin{picture}(400,160)(-50,-20)


\put(0,0){\line(1,1){160}}
\put(0,0){\vector(-1,-1){0}}

\multiput(40,40)(40,40){4}{\line(1,-1){40}}
\multiput(80,0)(40,40){4}{\vector(1,-1){0}}
\multiput(40,40)(40,40){4}{\circle*{5}}
\multiput(40,40)(40,40){4}{\line(0,-1){40}}
\multiput(40,0)(40,40){4}{\vector(0,-1){0}}

\multiput(0,0)(40,40){4}{\qbezier(45,-10)(60,-20)(75,-10)}
\multiput(45,-10)(40,40){4}{\vector(-3,2){0}}
\multiput(75,-10)(40,40){4}{\vector(3,2){0}}
\multiput(60,-25)(40,40){4}{$g$}

\put(20,-20){$\vdots$}

\end{picture}

\caption{A finite subgroup of $LS(X)$ not acting locally rigidly.
See Example~\ref{example:finite-not-lr}}
\label{fig::finite-not-lr}
\end{figure}

\begin{example} 
\label{example:finite-not-lr}
Let $X= \{ x_\infty, x_{a0}, x_{a1}, x_{a2},\dots ~|~ a\in\{0,1\}\}$ with ultrametric $d$ given by
$d(x_{ai},x_{aj})= e^{-\min\{ i,j\}}$ if $i\not= j$ and $a\in\{ 0,1\}$, and
$d(x_\infty, x_{0i}) = d(x_\infty, x_{1i}) = d(x_{0i}, x_{1i}) = e^{-i}$ for $i=0,1,2,\dots$.
The space $X$ is the end space of the tree in 
Figure~\ref{fig::finite-not-lr}.
Define $g\co X\to X$ by $gx_\infty=x_\infty$ and $gx_{ai}=x_{|a-1|i}$ for $i=0,1,2,\dots$ and $a\in\{ 0,1\}$.
Let $\Gamma$ be the subgroup of $LS(X)$ generated by $g$ (thus, $\Gamma$ is cyclic of order $2$).
Note that $\Gamma$ does not act locally rigidly on $X$ even though $\gammax$ is Hausdorff.
This example also shows that finite subgroups need not act 
locally rigidly.

\end{example}

\begin{theorem}
\label{thm:lr} If $X$ is an ultrametric space, then the following are equivalent:
\begin{enumerate}
\item
\label{thm:lr-one}
$X$ is locally rigid. 
\item
\label{thm:lr-two} 
$LS(X)$ acts locally rigidly on $X$.
\item
\label{thm:lr-three}
Every subgroup $\Gamma$ of $LS(X)$ acts locally rigidly on $X$.
\item
\label{thm:lr-four}
$LI(X)$ acts locally rigidly on $X$.
\item
\label{thm:lr-five}
There exists a group $\Gamma$ such that $Isom(X)\leq \Gamma\leq LS(X)$ and $\Gamma$ acts locally rigidly on $X$.
\item
\label{thm:lr-six}
$Isom(X)$ acts locally rigidly on $X$.
\item 
$\ligx$ is Hausdorff. 
\end{enumerate}
\end{theorem}

\begin{proof}
(1) implies (2):
Let $x\in X$ and $g\in LS(X)$ such that $gx=x$ and $sim(g,x)=1$.
Since $X$ is locally rigid, there exists $\epsilon_x>0$ such that if 
$0<\epsilon\leq\epsilon_x$ and $h\co B(x,\epsilon)\to B(x,\epsilon)$ is an isometry, then
$h=\id$.
Now choose $\epsilon >0$ such that $\epsilon\leq\epsilon_x$ and 
$g|\co B(x,\epsilon)\to B(x,\epsilon)$ is an isometry. Thus,
$g\in\Gamma_y$ for all $y\in B(x,\epsilon)$.

That (2) implies (3) implies (4) implies (5) implies (6) is obvious from the comment made above
that subgroups of groups acting locally rigidly also act locally rigidly.

(6) implies (1): Suppose on the contrary that $X$ is not locally rigid. Using
Theorem~\ref{theorem:Hausdorff-groupoid}~(\ref{theorem:Hausdorff-groupoid:2}), there exist 
$x\in X$ and a sequence $\epsilon_1 > \epsilon_2>\epsilon_3\cdots >0$ such that
$\lim_{i\to\infty}\epsilon_i=0$ together with 
isometries $h_i\co B(x,\epsilon_i)\to B(x,\epsilon_i)$
and $y_i\in B(x,\epsilon_i)\setminus B(x,\epsilon_{i+1})$ such that
$h_iy_i\not= y_i$.
Define $h\co X\to X$ by
$$h(z)= \left\{ \begin{array}{ll}
z & \textrm{if $z=x$ or $z\notin B(x,\epsilon_1)$}\\
h_iz & \textrm{if $z\in B(x,\epsilon_i)\setminus B(x,\epsilon_{i+1})$}\\
\end{array} \right.$$
Then $h\in Isom(X)$, $hx=x$ and $h$ is non-trivial arbitrarily close to $x$,
contradicting the assumption that $Isom(X)$ acts locally rigidly.

Finally, (1) and (7) are equivalent by Theorem~\ref{theorem:Hausdorff-groupoid}.
\end{proof}

We now discuss countability properties of the germ groupoid.

\begin{lemma} 
\label{lemma:sc}
If $X$ is a second countable ultrametric space with a countable subgroup
$\Gamma\leq LS(X)$, then the groupoid $\gammax$ is second countable.
\end{lemma}

\begin{proof}
Let $\Gamma= \{g_i\}_{i=1}^\infty$ and let $\{ x_j\}_{j=1}^\infty$ be a countable dense subset of $X$.
It follows that $\{ U(g_i, x_j, 1/k) ~|~ i,j,k\in\{ 1,2,3,\dots\}\}$ is a countable basis for $\gammax$.
For given any basis element $U(h,x,\epsilon)$ with $h\in\Gamma$, $x\in X$ and $\epsilon >0$, and
any germ $[h,y]\in U(h,x,\epsilon)$, simply choose $i,j,k$ such that $g_i=h$, $1/k<\epsilon$ and
$x_j\in B(y, 1/k)$.
Then $[h,y]\in U(g_i, x_j, 1/k)\subseteq U(h,x,\epsilon)$.
\end{proof}

\begin{example}
There exists a compact ultrametric space $X$ and an uncountable subgroup $\Gamma$ of $Isom(X)$ such that
$\Gamma$ acts locally rigidly on $X$ and $\gammax$ is not second countable.
Let $X$ be the end space of the Cantor tree 
(see Example~\ref{ex:Cantor Tree}). For $s\in\{ 0,1\}$, let
$\bar s= |s-1|$. For 
$x= (x_1, x_2, x_3, \dots)\in X$, 
let
$\overline{x} = (\overline{x_1}, \overline{x_2}, \overline{x_3}, \dots)$.
For 
$x= (x_1, x_2, x_3, \dots)\in X$, 
define $\alpha_x\co X\to X$ as follows. First, $\alpha_x(x)=\bar x$.
Second, if $x\not= y = (y_1, y_2, y_3,\dots)\in X$, let
$n=\min\{ i ~|~ x_i\not= y_i\}$ and
$\alpha_x(y)= (\overline{y_1}, \overline{y_2},\dots,\overline{y_n},y_{n+1},y_{n+2},\dots)$.

We show now that each $\alpha_x$ is an isometry.
For suppose $y,z\in X$ with $y\not= z\not= x\not= y$ and let
$\ell =\min\{ i ~|~ x_i\not= y_i \},
m =\min\{ i ~|~ x_i\not= z_i \},$ and
$n =\min\{ i ~|~ y_i\not= z_i \}.
$
We may write
$$y= (x_1,\dots, x_{\ell-1}, \overline{x_\ell}, y_{\ell + 1}, y_{\ell+2},\dots)$$
and
$$\alpha_x(y)= (\overline{x_1},\dots, \overline{x_{\ell-1}}, {x_\ell}, y_{\ell + 1}, y_{\ell+2},\dots).$$
It follows that 
$d(x,y)=d(\alpha_x(x),\alpha_x(y))$.
To see that
$d(y,z)=d(\alpha_x(y),\alpha_x(z))$,
assume without loss of generality that $\ell\leq m$.
We may write
$$z= (x_1,\dots, x_{m-1}, \overline{x_m}, z_{m + 1}, z_{m+2},\dots)$$
and
$$\alpha_x(z)= (\overline{x_1},\dots, \overline{x_{m-1}}, {x_m}, z_{m + 1}, z_{m+2},\dots).$$
It follows that $n\geq \ell$. If $n=\ell$, then
$$z= (x_1,\dots,x_{\ell-1},x_\ell,\dots, x_{m-1}, \overline{x_m}, z_{m + 1}, z_{m+2},\dots)$$
and
$$\alpha_x(z)= (\overline{x_1},\dots,\overline{x_\ell},\dots \overline{x_{m-1}}, {x_m}, z_{m + 1}, z_{m+2},\dots)$$
from which it follows that
$d(y,z)=d(\alpha_x(y),\alpha_x(z))$.
If $n>\ell$, then $m=\ell$ (from the ultrametric property) and we may write
$$z= (x_1,\dots,x_{\ell-1},\overline{x_\ell},y_{\ell+1},\dots, y_{n-1}, \overline{y_n}, z_{n+1}, z_{n+2},\dots)$$
and
$$\alpha_x(z)= (\overline{x_1},\dots,\overline{x_\ell-1},x_\ell, y_{\ell+1},\dots, y_{n-1}, \overline{y_{n}}, z_{n+1}, z_{n+2},\dots)$$
from which it follows that
$d(y,z)=d(\alpha_x(y),\alpha_x(z))$.
Finally, to see that $\alpha_x$ is bijective, note that
the inverse of $\alpha_x$ is given by $\alpha_x^{-1}= \alpha_{\overline x}$.

Let $\Gamma$ be the subgroup of $Isom(X)$ generated by $\{ \alpha_x ~|~ x\in X\}$.
Clearly, $\Gamma$ is uncountable.

We will now show that $\Gamma$ acts locally rigidly on $X$.
Let $x\in X$ and $\alpha\in\Gamma$ be given.
Write $\alpha=\alpha_k\circ\cdots\alpha_1$, where for each $1\leq j\leq k$
there exists $x^j\in X$ such that $\alpha_j=\alpha_{x^j}$.
Let $a^0=x$ and $a^j=\alpha_j\circ\cdots\alpha_1(x)$ for $1\leq j\leq k$.
Let $J= \{ j\in\{ 1,\dots, k\} ~|~ x^j\not= a^{j-1}\}$ and $J' =  \{ j\in\{ 1,\dots, k\} ~|~ j\notin J\}$.
For each $1\leq j\leq k$, let 
$$P_j= \left\{ \begin{array}{ll}
\min\{ i ~|~ x_i^j\not= a_i^{j-1}\} & \textrm{if $j\in J$}\\
\infty & \textrm{if $j\in J'$}.\\
\end{array} \right.$$
Let
$$P= \left\{ \begin{array}{ll}
\max\{ P_j ~|~ j\in J\} & \textrm{if $J\not=\emptyset$}\\
1 & \textrm{if $J=\emptyset$}.\\
\end{array} \right.$$
Now let $y\in X$ be any point such that $y\not=x$ and such that if 
$$N=\max\{ i ~|~ x_i\not= y_i\},$$ 
then $N >P$.
Suppose $\alpha(x)=x$.
We will prove that $\Gamma$ acts locally rigidly by showing
$\alpha(y)=y$.
We may write
$$y = (x_1,\dots,x_{N-1},\overline{x_N},y_{N+1}, y_{N+2},\dots).$$
Let $b^0=y$ and $b^j=\alpha_j\circ\cdots\circ\alpha_1(y)$ for $1\leq j\leq k$.
Since $\Gamma$ acts by isometries on $X$, we have
$a_i^j=b_i^j$, for $1\leq i\leq N-1$,
and $a_N^j\not= b_N^j$, whenever $0\leq j\leq k$.
In particular,
$(\alpha y)_i=y_i$ for $1\leq i\leq N$.
We will therefore be done once we establish the
\begin{claim}
For each $j\in\{ 0,\dots,k\}$, $b_i^j=y_i$ for $i\geq N+1$.
\end{claim}
\begin{proof}
This is true for $j=0$, so we proceed by induction, assuming $j>0$ and
$b_i^{j-1}=y_i$ for $i\geq N+1$.

\noindent
{\it Case 1.} $j\in J$.
Recall $P_j=\min\{ i ~|~ x_i^j\not= a_i^{j-1}\} \leq N-1$.
We may write
$$b^{j-1} = (a_1^{j-1},\dots,a_{N-1}^{j-1},\overline{a_N^{j-1}},y_{N+1}, y_{N+2},y_{N+3},\dots)$$
$$ = (x_1^{j},\dots,x_{P_j-1}^j,\overline{x_{P_j}^j},a_{P_j+1}^{j-1},\dots,a_{N-1}^{j-1},\overline{a_N^{j-1}},y_{N+1}, y_{N+2},y_{N+3},\dots).$$
Thus,
$$b^j=\alpha_j(b^{j-1})
= (\overline{x_1^{j}},\dots,\overline{x_{P_j-1}^j},{x_{P_j}^j},a_{P_j+1}^{j-1},\dots,a_{N-1}^{j-1},\overline{a_N^{j-1}},y_{N+1}, y_{N+2},y_{N+3},\dots),$$
and there is agreement where claimed.

\noindent
{\it Case 2.} $j\in J'$.
In this case $x^j=a^{j-1}$ and we may write
$$b^{j-1} = (a_1^{j-1},\dots,a_{N-1}^{j-1},\overline{a_N^{j-1}},y_{N+1}, y_{N+2},y_{N+3},\dots)$$
$$ = (x_1^{j},\dots,x_{N-1}^j,\overline{x_{N}^j},y_{N+1}, y_{N+2},y_{N+3},\dots).$$
Thus,
$$b^j=\alpha_j(b^{j-1})
= (\overline{x_1^{j}},\dots,\overline{x_{N-1}^j},{x_{N}^j},y_{N+1}, y_{N+2},y_{N+3},\dots),$$
and there is agreement where claimed.

This completes the proof of the claim and also the assertion that $\Gamma$ acts locally rigidly on $X$.
\end{proof}

We now show that $\gammax$ is not second countable.
Note that $\{ U(\alpha_x,x,1) ~|~ x\in X\}$ is an uncountable collection of open subsets of
$\gammax$ and the germ $[\alpha_x,x]\in U(\alpha_x,x,1)$ for every $x\in X$.
However, if $x\not= y\in X$, then $\alpha_y(y)=\overline{y}$ while $\alpha_x(y)\not=\overline{y}$.
Thus, $[\alpha_y,y]\notin U(\alpha_x,x,1)$. This implies that $\gammax$ has no countable basis.

Finally, note that $\Gamma$ does not act freely on $X$.
For example, let
$x^1= (1,0,1,0,1,0,\dots),
x^2= (1,1,0,0,0,0,\dots),
x^3= (1,0,0,0,0,0,\dots)$ and
$p= (0,0,0,0,\dots)$.
Then $\alpha=\alpha_{x^3}^2\circ\alpha_{x^2}\circ\alpha_{x^1}\in\Gamma$,
$\alpha(p)=p$ and $\alpha(x^1)\not= x^1$. Hence, $\alpha\not= 1$ but it fixes the point $p$.

This completes the discussion of the example.
\end{example}

Finally, we give a proof of Theorem~\ref{theorem:second-main}(1) from the Introduction.

\begin{theorem}
\label{theorem:second-main-one-proof}
If $X$ is a compact ultrametric space with a countable subgroup $\Gamma\leq LS(X)$ acting locally rigidly
on $X$, then
$\gammax$ is a locally compact, Hausdorff, second countable, \'etale groupoid.
\end{theorem}

\begin{proof}
This follows from
Lemma~\ref{lemma:etale},
Lemma~\ref{lemma:loccom},
Remark~\ref{remark:opensubgrpoids},
Theorem~\ref{theorem:Hausdorffness}(2), and
Lemma~\ref{lemma:sc}.
\end{proof}

\begin{remark} If in the hypothesis of Theorem~\ref{theorem:second-main-one-proof},
``compact'' is replaced by ``locally compact and second countable,''
then the conculsion still holds with the same proof.
\end{remark}

\begin{proposition}
\label{prop:ConjugateActions}
If $X$ is a metric space with a subgroup $\Gamma\leq LS(X)$ acting locally rigidly on $X$ and $h\in LS(X)$, then $h^{-1}\Gamma h$ also acts locally rigidly on $X$.
\end{proposition}

\begin{proof}
Let $x\in X$ and $g\in (h^{-1}\Gamma h)_x$ such that $sim(g,x)=1$ be given.
We need to show that there exists $\epsilon >0$ such that 
$g|\co B(x,\epsilon)\to X$ is the inclusion.
Note that $hgh^{-1}\in \Gamma_{hx}$. 
To see that $sim(hgh^{-1}, hx)=1$, choose $\delta >0$ and $\lambda >0$
such that $h|\co B(x,\delta)\to B(hx,\lambda\delta)$ is a $\lambda$-similarity. 
Then $h^{-1}|\co B(hx,\lambda\delta)\to B(x,\delta)$ is a $\lambda^{-1}$-similarity. Assume that $\delta$ is small enough that $g|\co B(x,\delta)\to B(x,\delta)$ is a $1$-similarity.
Then $hgh^{-1}|\co B(hx,\lambda\delta)\to B(hx,\lambda\delta)$ is a 
$1$-similarity.

Since $hgh^{-1}\in\Gamma$, it follows that there exists $\epsilon >0$ such that $hgh^{-1}|\co B(hx,\epsilon)\to X$  is the inclusion, from which it follows that 
$g|\co B(x,\epsilon)\to X$ is the inclusion.
\end{proof}


\section{The approximating groupoids}
\label{section:TheApproximatingGroupoids}

This section contains a proof that $\ligx$ is an AF groupoid if $X$ is a compact,
locally rigid ultrametric space.

Throughout this section, let $(X,d)$ denote an ultrametric space. 

\begin{definition} The {\it pseudogroup $\plix$ of local isometries on $X$\/}  is 
the set of all isometries between open subsets of $X$. That is, an element of $\plix$ 
consists of open subsets $U, V$ of $X$ and an isometry $g\co U\to V$.
\footnote{Of course, 
$\plix$ has the structure of a pseudogroup, but we do not explicitly use it.}
\end{definition}

\begin{definition} Let $\epsilon >0$. The {\it $\epsilon-$local isometry groupoid\/}
$\eligx$ {\it of $X$\/} is the subset of $\plix\times X$ given by
$$\eligx=\{ (g,x)\in \plix\times X ~|~ g\co B(x,\epsilon)\to B(gx,\epsilon)\}.$$
\end{definition}

Thus, $(g,x)\in\eligx$ means $g$ is an isometry from $B(x,\epsilon)$ onto
$B(gx,\epsilon)$.

The groupoid structures on $\eligx$ are the obvious ones. Thus, 
the unit space is $X$; the
domain $d\co \eligx\to X$ and range $r\co \eligx\to X$ maps are given by
$d(g,x)=x$ and $r(g,x)=gx$. If $(g_1,x_1)$ and $(g_2,x_2)$ are in $\eligx$, then the composition
is defined by $(g_2,x_2)(g_1,x_1)=(g_2g_1,x_1)$ provided $x_2=g_1x_1$.\footnote{Thus,
$[\eligx]^2=\{ ((g_2,x_2), (g_1,x_1))\in\eligx\times\eligx ~|~ x_2=g_1x_1\}.$}
The inverse is
$(g,x)^{-1}=(g^{-1},gx)$.

A basis for a topology on $\eligx$ consists of all sets $U_\delta(g,x)$ where
$(g,x)\in\eligx$, $0<\delta\leq\epsilon$ and
$$U_\delta(g,x)= \{ (h,y)\in\eligx ~|~ d(x,y)<\delta \textrm{ and $d(gz,hz)<\delta$ for all
$z\in B(x,\epsilon)$}\}.$$

\begin{proposition}
\label{prop:approx-groupoid}
If $X$ is an ultrametric space and $\epsilon >0$, then
$\eligx$ is a Hausdorff topological groupoid. Moreover, the domain and range maps are
open.
\end{proposition}

\begin{proof} 
To see that the collection of all $U_\delta(g,x)$ forms a basis, first note that the
collection certainly covers $\eligx$. And if $(g,x)\in U_{\delta_1}(g_1,x_1)\cap
U_{\delta_2}(g_2,x_2)$, let $\delta=\min\{\delta_1,\delta_2\}$ and observe that
$U_\delta(g,x)\subseteq U_{\delta_1}(g_1,x_1)\cap
U_{\delta_2}(g_2,x_2)$.

To see that the resulting topology is Hausdorff, let $(g_1,x_1)\not= (g_2,x_2)$ in $\eligx$
and choose $0<\delta\leq\epsilon$ such that 
$$\delta < \left\{ \begin{array}{ll}
d(x_1,x_2) &\textrm{if $x_1\not= x_2$}\\
\sup\{ d(g_1z,g_2z) ~|~ z\in B(x,\epsilon)\} &\textrm{if $x_1=x_2$}
\end{array} \right.$$
and observe that $U_\delta(g_1,x_1)\cap U_\delta(g_2,x_2)=\emptyset$.

To see that $d,r\co \eligx\to X$ are continuous, let $x\in X$ and $0<\delta\leq\epsilon$.
Then one can check that $d^{-1}(B(x,\delta))=\cup\{ U_\delta(h,y) ~|~ d(x,y)<\delta\}$ and
$$r^{-1}(B(x,\delta))=\cup\{ U_\delta(h,h^{-1}y) ~|~ d(x,y)<\delta\}.$$

To see that $d,r$ are open, let $(h,y)\in\eligx$ and $0<\delta\leq\epsilon$. Then one can check
that $d(U_\delta(h,y))=B(y,\delta)$ and $r(U_\delta(h,y))=B(hy,\delta)$.

That multiplication $m\co [\eligx]^2\to \eligx$ is continuous follows from the following fact:
if $(h,y)\in U_\delta(g,x)$ and $m((hk^{-1},ky),(k,y))=(h,y)$, then
$$[U_\delta(hk^{-1},ky)\times U_\delta(k,y)]
\cap[\eligx]^2\subseteq m^{-1}(U_\delta(g,x)).$$

Finally, inversion is continuous because $[U_\delta(g,x)]^{-1}=U_\delta(g^{-1},gx)$.
\end{proof}

\begin{theorem} 
\label{theorem:approx-groupoid}
If $X$ is a compact, locally rigid ultrametric space, then
there exists $\epsilon_X >0$ such that for every $0<\epsilon\leq\epsilon_X$:
\begin{enumerate}
\item  $\eligx$ is a Hausdorff, locally compact, 
\'etale groupoid,
\item $\eligx$ is an elementary groupoid in the sense of 
Renault \cite{Ren}.
\end{enumerate}\end{theorem}

\begin{proof}  Let $\epsilon_X$ be given by Lemma~\ref{lemma:compact-locally-rigid}. 

(1) For $0<\epsilon\leq\epsilon_X$,
we already know that $\eligx$ is Hausdorff (Proposition~\ref{prop:approx-groupoid}). To say that it is \'etale means
$r\co \eligx\to X$ is a local homeomorphism. To verify this it suffices to let 
$(g,x)\in\eligx$ and show that $r|\co U_\delta(g,x)\to B(gx,\delta)$ is injective whenever
$0<\delta\leq\epsilon$ (because the proof of Proposition~\ref{prop:approx-groupoid} 
shows $r|$ is continuous, open and
surjective).
To this end let $(h_i,y_i)\in U_\delta(g,x)$ for $i=1,2$ such that $r(h_1,y_1)=h_1y_1=
h_2y_2=r(h_2,y_2)$. Then $h_i\co B(y_i,\epsilon)\to B(h_iy_i,\epsilon)$ is an isometry
for $i=1,2$. Hence, $h=h_2^{-1}h_1\co B(y_1,\epsilon)\to B(y_2,\epsilon)=B(y_1,\epsilon)$ is
an isometry. The choice of $\epsilon_X$ implies $h$ is the identity, so $h_1=h_2$ from which
it also follows that $y_1=y_2$. 
Finally, note that this also implies that $\eligx$ is locally compact, being locally 
homeomorphic to the compact space $X$.

(2) Let $0<\epsilon\leq\epsilon_X$. 
According to Renault \cite[page 123]{Ren} we need to show that $\eligx$ is the disjoint union
of a sequence of elementary groupoids $G_i$ of type $n_i$ 
(the definitions will be recalled below). In
fact, we will show that the sequence is finite, say $1\leq i\leq i_\epsilon$.
Let $\mathcal{B}_\epsilon$ be the collection of all open $\epsilon$-balls in $X$. Since $X$ is compact
ultrametric, $\mathcal{B}_\epsilon$ is a finite collection and any two distinct members of
$\mathcal{B}_\epsilon$ are disjoint. By the choice of $\epsilon_X$, if $B_1,B_2\in\mathcal{B}_\epsilon$
then either there exists a unique isometry $B_1\to B_2$ or, $B_1$ and $B_2$ are not isometric.
Thus, we may express $\mathcal{B}_\epsilon$ as a finite disjoint union $\cup_{i=1}^{i_\epsilon}
\mathcal{B}_i$ such that if $1\leq i,j\leq i_\epsilon$, $B_1\in\mathcal{B}_i$, $B_2\in\mathcal{B}_j$, 
then $B_1$ and $B_2$ are isometric if and only if $i=j$; moreover, if $i=j$, then there exists
a unique isometry $B_1\to B_2$.

For $1\leq i\leq i_\epsilon$, let
$G_i=\{ (g,x)\in \eligx ~|~ B(x,\epsilon)\in\mathcal{B}_i\}$
and let $n_i$ equal the cardinality of $\mathcal{B}_i$.
Clearly, $\eligx=\cup_{i=1}^{i_\epsilon}G_i$ and the $G_i$'s are mutually disjoint subgroupoids 
of $\eligx$. It remains to show that each $G_i$ is elementary of type $n_i$.
Given $i$ choose $x_i\in X$ such that $B(x_i,\epsilon)\in\mathcal{B}_i$.
Let $\hat G_i= \{ gx_i ~|~ (g,x_i)\in G_i\}\subseteq X$. 
Note that if $(g,x_i),(h,x_i)\in G_i$ and $gx_i=hx_i$, then $g=h$.

Now $\hat G_i\times\hat G_i$ has a natural groupoid structure with set of composable pairs
$[\hat G_i\times\hat G_i]^2=
\{((g_1x_i,g_2x_i),(g_3x_i,g_4x_i))\in\hat G_i\times \hat G_i\times\hat G_i\times\hat G_i
~|~ g_2x_i=g_3x_i\}$,
unit space $(\hat G_i\times\hat G_i)^0=\hat G_i\subseteq X$,
$d\co \hat G_i\times\hat G_i\to\hat G_i$ and $r\co \hat G_i\times\hat G_i\to\hat G_i$
given by 
$d(gx_i,hx_i)=hx_i$ and $r(gx_i,hx_i)=gx_i$ and multiplication
$$(g_1x_i,g_2x_i)\cdot (g_2x_i,g_3x_i)=(g_1x_i,g_3x_i).$$
Clearly, $(r,d)\co \hat G_i\times\hat G_i\to\hat G_i\times\hat G_i$ is bijective (in fact, it is 
the identity). This is what it means to be a transitive principal groupoid on $n_i$ elements \cite[page 6]{Ren}. 

Now give $B(x_i,\epsilon)$ the trivial groupoid structure 
(that is, $B(x_i,\epsilon)$ is the unit space
and $(x,y)$ is composable if and only if $x=y$). Since $X$ is compact, $B(x_i,\epsilon)$ is
a second countable metric space.
Note that $G_i$ is isomorphic to the product
$B(x_i,\epsilon)\times (\hat G_i\times\hat G_i)$ via
$G_i\to B(x_i,\epsilon)\times (\hat G_i\times\hat G_i)$; $(g,x)\mapsto (h^{-1}x,ghx_i,hx_i),$
where $h\co B(x_i,\epsilon)\to B(x,\epsilon)$ is the unique isometry. This means that $G_i$ is
an elementary groupoid of type $n_i$.
\end{proof}

\begin{remark} 
\label{remark:approx}
Under the hypothesis and notation of  Theorem~\ref{theorem:approx-groupoid}, 
note that the topology of
$\eligx$ is second countable, being a finite union $\cup_{i=1}^{i_\epsilon}G_i$ and each
$G_i$ is homeomorphic to a product of $B(x_i,\epsilon)$ and a finite set. 
In particular, $\eligx$ is compact.
In fact, since
$B(x_i,\epsilon)$ is closed in $X$, hence compact, $\eligx$ has a countable basis of compact
open sets. In fact, there is a countable basis of compact open $\eligx$-sets in the sense
of Renault \cite[page 10]{Ren}. To see this, let $g\co B(x_i,\epsilon)\to B(gx_i,\epsilon)$ and
$h\co B(x_i,\epsilon)\to B(hx_i,\epsilon)$ be isometries, and let $A(g,h)=
\{ (gh^{-1},hy) ~|~ y\in B(x_i,\epsilon)\}$. These sets correspond to the images of
$B(x_i,\epsilon)$ under the constructions giving a basis of compact open sets for $\eligx$. 
Since $d(gh^{-1},hy)=hy$ and $r(gh^{-1},hy)=gy$, $d,r$ restricted to $A(g,h)$ are injective as
required. This property of $\eligx$ is important when applying Renault's results on AF groupoids
and AF algebras (see \cite[page 130]{Ren}).
\end{remark}

\begin{theorem}
\label{theorem:AF}
If $X$ is a compact, locally rigid ultrametric space, then:
\begin{enumerate}
\item $\ligx$ is an AF groupoid in the sense of Renault \cite{Ren},
\item The groupoid $C^*$-algebra $C^*\ligx$ is a unital AF $C^*$-algebra.
\end{enumerate}\end{theorem}

\begin{proof} (1) First note that the unit space $X$ of $\ligx$ is totally disconnected (since 
it is ultrametric). Thus, we only need to show that $\ligx$ is the inductive limit of a 
sequence of elementary groupoids (see \cite[pages 122--123]{Ren}). For this choose a sequence
$\epsilon_X > \epsilon_1> \epsilon_2>\cdots$ such that $\lim_{i\to\infty}\epsilon_i=0$ 
where $\epsilon_X$
is given by Lemma~\ref{lemma:compact-locally-rigid}.
We will observe that
$$\ligx=\lim_{\to}\eiligx .$$
First note that $\eiligx$ is an open subgroupoid of $\ligx$ via the embedding 
$(g,x)\mapsto [g,x]$. For $0<\delta\leq\epsilon_i$, the embedding takes $U_\delta(g,x)$
onto $U(g,x,\epsilon_i)$.
Likewise, $\eiligx$ is an open subgroupoid of $\eioligx$ via the embedding 
$(g,x)\mapsto (g|B(x,\epsilon_{i+1}),x)$.
Finally observe that $\ligx=\cup_{i=1}^\infty\eiligx$.

(2) That the groupoid $C^*$-algebra is AF 
follows from (1) and Renault \cite[1.15, page 134]{Ren}.
It is unital because $X$ is compact.
\end{proof}

It should be mentioned that Renault proved that every AF $C^*$-algebra
is the $C^*$-algebra of an AF groupoid and the groupoid is unique
up to isomorphism \cite[1.15, page 134]{Ren}.


\section{Trees, Bratteli diagrams and path groupoids}
\label{section:TreesBratteliDiagramsPathGroupoids}

Let $(T,v)$ be a rooted, geodesically complete,
locally finite simplicial tree. 
The purpose of this section is to define a Bratteli diagram 
$B(T,v)$ associated with $(T,v)$ and to prove that $\ligx$ is 
isomorphic to the path groupoid of
$B(T,v)$, provided $X = end(T,v)$ is locally rigid.


\subsection{Recollections on Bratteli diagrams}
\label{subsec:RecollectionOfBratteli}

The material in this section is well-known. 
See Blackadar \cite{Bla}, Bratteli \cite {Bra}, 
Davidson \cite{Dav},
Effros \cite{Eff}, 
Elliott \cite{Ell}, Exel and Renault \cite{ExR},
Giordano, Putnam, and Skau \cite{GPS}, and
Herman, Putnam and Skau \cite{HPS}
for more details. In particular, the discussion below relies
heavily on the expositions in \cite{GPS} and \cite{HPS}.

We begin with the definition of
a Bratteli diagram, which, for us, comes equipped with a distinguished
initial vertex.

\begin{definition} A directed graph
$D=(\mathcal{V},\mathcal{E})$ 
with vertex set $\mathcal{V}$, edge set $\mathcal{E}$,
initial map $s\co  \mathcal{E}\to\mathcal{V}$, and terminal map
$r\co \mathcal{E}\to\mathcal{V}$ is a {\it Bratteli diagram} if
\begin{enumerate}
\item $\mathcal{V}$ is given as the union $\mathcal{V}=\bigcup_{n=0}^\infty\mathcal{V}_n$ of
mutually disjoint, finite, nonempty sets $\mathcal{V}_n$,
\item for each edge $e\in\mathcal{E}$, if the initial vertex $s(e)\in\mathcal{V}_n$, then the terminal vertex
$r(e)\in\mathcal{V}_{n+1}$,
\item for each vertex $v\in\mathcal{V}$ there are at most finitely many edges $e\in\mathcal{E}$ with $s(e)=v$,
\item $\mathcal{V}_0$ consists of a single vertex $v_0$,
\item every vertex is the initial vertex of some edge,
\item every vertex except $v_0$ is the terminal vertex of some edge.
\end{enumerate}
\end{definition}

For example, let $(T,v)$ be a rooted, geodesically complete, locally finite, simplicial tree. By specifying
the root, $T$ is naturally a directed graph (edges are directed away from the root). Thus, $T$ is a Bratteli diagram,
where $\mathcal{V}_n$ consists of those vertices a distance $n$ from $v$ (with respect to the metric discussed in
Section\ref{section:trees}).

We now recall the construction of two invariants associated to a Bratteli
diagram $D=(\mathcal{V},\mathcal{E})$, namely, the unital dimension
group $(G(D), G_+(D), [1])$ and the unital AF $C^*$-algebra $AF(D)$.
Both of these invariants depend on a sequence of matrices, 
which we now define.

For each $i=0, 1, 2, \dots$, let $m_i= |\mathcal{V}_i|$, the cardinality
of $\mathcal{V}_i$, and write
$\mathcal{V}_i =\{ v_1^i, \dots, v_{m_i}^i\}$ (in particular, $v_0=v_1^0$).
For $i= 0, 1, 2, \dots$, $1\leq k\leq m_{i+1}$, and $1\leq \ell\leq m_i$,
let 
$$a_{k\ell}^i = |\{ e\in\mathcal{E} ~|~ \text{$s(e) = v_\ell^i$
and $r(e) = v_k^{i+1}$}\}|,$$
the number of edges from the $\ell^{th}$ vertex at the $i^{th}$
level to the $k^{th}$ vertex at the $(i+1)^{st}$ level.
Thus, $A_i\co = [a_{k\ell}^i]$ is an $(m_{i+1}\times m_i)$-matrix
with nonnegative integral entries. Moreover, no column and no row
of $A_i$ consists entirely of zeroes.

The direct limit $G(D)$ of the sequence 
$$\bz \stackrel{A_0}{\longrightarrow}
\bz^{m_1} \stackrel{A_1}{\longrightarrow} \bz^{m_2} 
\stackrel{A_2}{\longrightarrow}\bz^{m_3} 
\stackrel{A_3}{\longrightarrow} \cdots
\bz^{m_i} \stackrel{A_i}{\longrightarrow} 
\bz^{m_{i+1}}{\longrightarrow}\cdots$$
is a partially ordered abelian group with positive cone $G_+(D)$
given by the direct limit of 
$$\bz_+ \stackrel{A_0}{\longrightarrow}
\bz_+^{m_1} \stackrel{A_1}{\longrightarrow}\bz_+^{m_2} \stackrel{A_2}{\longrightarrow}
\bz_+^{m_3} \stackrel{A_3}{\longrightarrow}\cdots
\bz_+^{m_i} \stackrel{A_i}{\longrightarrow}\bz^{m_{i+1}}_+
{\longrightarrow}\cdots.$$
(Here $\bz_+=\{ 0,1,2,\dots\}$.)

\begin{definition}
The pair $(G(D), G_+(D))$ is the {\it dimension group associated
to the Bratteli diagram $D$}. The class $[1]\in G(D)$ of the unit $1\in\bz$
is an order unit\footnote{An element $u$ of the positive cone $G_+$ of a 
partially ordered abelian group $G$ is an {\it order unit}
if for every $x\in G$ there exists $n\in\bn$ such that
$x\leq nu$.}
and the triple $(G(D), G_+(D), [1])$ is
the {\it unital dimension group associated to the Bratteli diagram $D$}.
\end{definition}

The second invariant associated to a
Bratteli diagram $D$ by using the sequence of matrices $A_i$
is a direct limit of
finite dimensional $C^*$-algebras (i.e., finite direct sums of
matrix algebras over $\bc$) defined as follows.

In general, let $\bm_r$ denote the $C^*$-algebra 
of $(r\times r)$-matrices over $\bc$.
For each $v\in\mathcal{V}$, let $k(v)$ be the number of directed
paths in $D$ from $v_0$ to $v$.
Let $C_0=\bc=\bm_1$ and, for $i=1,2,3,\dots$, let
$$C_i= \bigoplus_{j=1}^{m_i} \bm_{k(v_j^i)}.$$
The matrices $A_i:= [a_{k\ell}^i]$ defined above may be considered
to be matrices of multiplicities determining
unital $C^*$-algebra homomorphisms, also denoted $A_i$,
$A_i\co C_i\to C_{i+1}$ for each $i=0, 1, 2,\dots$.
Hence, there is a direct sequence
$$\bc = C_0 \stackrel{A_0}{\longrightarrow}
C_1 \stackrel{A_1}{\longrightarrow} C_2 \stackrel{A_2}{\longrightarrow}
C_3 \stackrel{A_3}{\longrightarrow}\cdots
C_i \stackrel{A_i}{\longrightarrow} C_{i+1}
{\longrightarrow}\cdots.$$
Let $\mathrm{AF}(D)$ denote the $C^*$-direct limit of the sequence
just described. It is a unital AF algebra with unit $[1]$, the class of
$1\in\bc$. 

The two invariants of a Bratteli diagram
defined above are related via $K$-theory.
It is well-known that the $K_0$ group of the $C^*$-algebra
$\mathrm{AF}(D)$
is $G(D)$; in fact, the unital, partially ordered
abelian groups, $(K_0(\mathrm{AF}(D)), K_0(\mathrm{AF}(D))_+, [1])$
and $(G(D), G_+(D), [1])$, are isomorphic (see \cite{Dav}).

Finally, we recall the equivalence relation on Bratteli diagrams
that are classified by these invariants.

\begin{definition} A {\it telescoping} of a Bratteli diagram
$D= (\mathcal{V}, \mathcal{E})$ to a Bratteli diagram 
$D'= (\mathcal{V}', \mathcal{E}')$ consists of a subsequence
$0=m_0 < m_1 < m_2<\cdots$ of $\bz_+$
such that
\begin{enumerate}
\item $\mathcal{V}_n'=\mathcal{V}_{m_n}$ for all $n=0,1,2,\dots$, and
\item if $n\in\bz_+$, $x\in\mathcal{V}_n'$, and $y\in\mathcal{V}_{n+1}'$,
then the number of edges in $D'$ from $x$ to $y$ is exactly the number of
directed paths in $D$ from $x$ to $y$.
\end{enumerate}
\end{definition}

Two Bratteli diagrams 
$D= (\mathcal{V}, \mathcal{E})$ and 
$D'= (\mathcal{V}', \mathcal{E}')$
are {\it isomorphic} if 
there exists an isomorphism $\varphi\co  D\to D'$ of directed graphs
such that $\varphi(\mathcal{V}_n)=\mathcal{V}_n'$ for all $n=0, 1, 2,\dots$.
They are {\it equivalent} if they are equivalent 
under the equivalence relation generated by isomorphism and telescoping.

Two partially ordered abelian groups $(G, G_+)$ and $(G', G_+')$
are {\it isomorphic} is there is a group isomorphism $\varphi\co G\to G'$
such that $\varphi(G_+)=G_+'$. If in addition $u\in G$ and $u'\in G'$ 
are given order units and $\varphi(u)=u'$, then the unital
partially ordered abelian groups $(G, G_+, u)$ and $(G', G_+', u')$
are {\it isomorphic}. 

\begin{theorem}[Bratteli, Elliott] 
\label{theorem:Bratteli, Elliott}
For two Bratteli diagrams $D, D'$, the following are
equivalent:
\begin{enumerate}
\item $D$ and $D'$ are equivalent Bratteli diagrams.
\item $(G(D), G_+(D), [1])$ and $(G(D'), G_+(D'), [1])$ are isomorphic
unital partially ordered abelian groups.
\item $(\mathrm{AF}(D), [1])$ and  $(\mathrm{AF}(D'), [1])$ are isomorphic
unital $C^*$-algebras.
\end{enumerate}
Moreover, 
$(G(D), G_+(D))$ and $(G(D'), G_+(D'))$ are isomorphic
partially ordered abelian groups if and only if
$\mathrm{AF}(D)$ and  $\mathrm{AF}(D')$ are stably isomorphic
$C^*$-algebras.
\end{theorem}

The equivalence of the first two conditions is due to Bratteli \cite{Bra};
the equivalence of the second two, as well as the final statement, is
due to Elliott \cite{Ell}.


\subsection{The Bratteli diagram $B(T,v)$ associated to a tree $(T,v)$}
\label{subsec:ConstructionOfBTv}

Let $(T,v)$ be a rooted, geodesically complete, locally finite simplicial tree. 
As mentioned above, the choice of root $v$ 
gives an orientation to each edge of $T$: the edges point away from the root. Thus, $(T,v)$ is a connected,
directed graph (in fact, a Bratteli diagram). Let $s(e)$ denote the initial, and $r(e)$ the terminal, vertex of the edge $e$.

For notation, let $vert(T)$ be the set of vertices of $T$ and 
let $V_i$ be the set of vertices at level $i$. Thus,
$$V_i=\{ w\in vert(T) ~|~ \textrm{the minimal simplicial path from $v$ to $w$ has length $i$}\}.$$
For each $w\in vert(T)$, $w\not= v$, let $T_w$ denote the subtree of $T$ descending from $w$.\footnote{Thus, $T_w$ 
contains all vertices and edges of $T$ that are in directed paths beginning at $w$.}
We let $T_w$ be rooted at $w$. If $w\in V_i$, then we say {\it $(T_w,w)$ is a level $i$ 
rooted subtree of $(T,v)$.} 

In turn, a {\it level one subtree of a level $i$ rooted subtree $(T_w,w)$ of $(T,v)$\/}  is a level 
$(i+1)$ subtree $(T_u,u)$ of $(T,v)$ that is also a subtree of $(T_w,w)$ (i.e., $u$ is a vertex of $T_w$).
For each $i=0,1,2,\dots$ 
let $m_i$ be the number of rooted isometry classes of level $i$ subtrees of $(T,v)$ and
let $T_1^i, T_2^i,\dots T_{m_i}^i$ be a complete set of representatives of
the rooted isometry classes.  
Note that $m_0=1$ and $T_1^0=T$.
Thus, for each level $i$ subtree $S$ of $T$ there exists a unique integer $l$ such that $1\leq l\leq m_i$ and 
$T_l^i$ is rooted isometric to $S$.
Call these chosen subtrees the {\it admissible} ones.

For each $i\geq 1$ and for each level $i$ subtree $S$ of $T$, choose a rooted isometry
$$\alpha(T_l^i,S)\co T_l^i\to S$$
where $l$ is the unique integer such that  $1\leq l \leq m_i$ and $T_l^i$ is rooted isometric to $S$.
In choosing these isometries, insist that
$$\alpha(T_l^i,T_l^i)=\id_{T_l^i} \hbox{~for each $i\geq 1$ and $1\leq l\leq m_i$}.$$
Call these chosen rooted isometries the {\it admissible} ones. 

Define an equivalence relation $\sim$ on $T$ as follows. For two distinct vertices $w_1, w_2$ of $T$, we have 
$w_1\sim w_2$ if and only if there exists $i\geq 1$ such that $w_1, w_2\in V_i$ and $T_{w_1}$ is rooted
isometric to $T_{w_2}$.
For two distinct edges $e_1, e_2$ of $T$, we have $e_1\sim e_2$ if and only if
each of the following hold:
\begin{enumerate}
\item there exists $i\geq 1$ such that $s(e_1), s(e_2) \in V_i$,
\item $T_{s(e_1)}$ is rooted isometric to $T_{s(e_2)}$ (in particular, $s(e_1)\sim s(e_2)$),
\item if $l$ is the unique integer with $1\leq l\leq m_i$ such that $T_l^i$ is rooted isometric to 
$T_{s(e_1)}$ (which, of course, also implies $T_l^i$ is rooted isometric to $T_{s(e_2)}$), then 
$$\alpha(T_l^i, T_{s(e_1)})^{-1}(e_1) = \alpha(T_l^i, T_{s(e_2)})^{-1}(e_2)$$
as edges of $T_l^i$.
\end{enumerate}

Let $B(T,v) = T/\sim$, which has the structure of a connected directed graph. Level $i$ vertices of
$B(T,v)$ are equivalence classes of level $i$ vertices of $(T,v)$, and edges of $B(T,v)$ are equivalence
classes of edges of $(T,v)$ with the induced orientation. There is an initial vertex of $B(T,v)$, namely the class
$[v]$ (which consists only of $v$). Thus, $B(T,v)$ is a Bratteli diagram
and is called {\it the Bratteli diagram associated to $(T,v)$}.

Note that the quotient map $\kappa\co  T\to B(T,v)$ is a morphism of directed graphs.

\begin{proposition}
The  Bratteli diagram $B(T,v)$ is well-defined
up to isomorphism.
\end{proposition}

\begin{proof}
We must show that if other choices of level $i$
subtrees 
$T_1^i, T_2^i,\dots T_{m_i}^i$ and admissible isometries $\alpha(T_l^i,S)$
are made, then the resulting Bratteli diagram is isomorphic to $B(T,v)$.
The vertex set $\mathcal V$ of $B(T,v)$, and its expression as 
$\mathcal{V}=\bigcup_{n=0}^\infty\mathcal{V}_n$,
is obviously independent of the choices. 

Thus, it remains to show that if $w_1, w_2$ are vertices of $T$ such that
$w_1\in V_i$ and $w_2\in V_{i+1}$ for some $i$, then the number of edges from $[w_1]$
to $[w_2]$ in $B(T,v)$ is independent of the choices.
For this, note that the number of edges in $B(T,v)$ beginning at $[w_1]$ equals the number of
edges in $T_{w_1}$ beginning at $w_1$, and $e\mapsto \kappa(e)$, where $e$ is an edge in $T_{w_1}$ beginning at
$w_1$, gives the bijection. 
Now observe that for such an edge $e$ in $T_{w_1}$  from $w_1$ to some $w_3$, its image 
$\kappa(e)$ ends at $[w_2]$ if and  only if $T_{w_2}$ is rooted isometric to $T_{w_3}$.
\end{proof}

\begin{remark}
\label{remark:explicit}
We give here an explicit description of the vertices $\mathcal V$
and edges $\mathcal E$ of the Bratteli diagram $B(T,v)$.
For $i= 0, 1, 2, \dots$, the level $i$ vertices can be written as a set
of equivalence classes $\mathcal{V}_i =\{ [v_1^i], \dots, [v_{m_i}^i]\}$,
where $v_\ell^i$ is the root of $T_\ell^i$ $(1\leq\ell\leq m_i)$.
The number of edges in $B(T,v)$ from $[v_\ell^i]$ to
$[v_k^{i+1}]$, where $i= 0, 1, 2, \dots$, $1\leq\ell\leq m_i$, and $1\leq k\leq 
m_{i+1}$, is nonzero if and only if there exists $w\in [v_k^{i+1}]$
such that $w\in T_\ell^i$. When such a vertex $w$ exists,
the number of edges from $[v_\ell^i]$ to
$[v_k^{i+1}]$ is the number of level $1$ subtrees of $T_\ell^i$ that
are rooted isometric to $T_k^{i+1}$.
\end{remark}

\begin{figure}[htbp]

\begin{picture}(350,70)(-10,0)

\put(0,30){\circle*{5}}

\multiput(60,0)(60,0){5}{\circle*{5}}
\multiput(60,60)(60,0){5}{\circle*{5}}

\put(0,30){\line(2,-1){60}}
\put(0,30){\line(2,1){60}}

\multiput(60,60)(60,0){4}{\line(1,-1){27}}
\multiput(120,0)(60,0){4}{\line(-1,1){27}}

\multiput(60,0)(60,0){4}{\line(1,1){60}}
\multiput(60,0)(60,0){4}{\line(1,0){60}}

\multiput(310,30)(10,0){3}{\circle*{2}}

\end{picture}
\caption{Bratteli diagram of the Fibonacci tree $F$}
\label{fig:brat-fib-tree}
\end{figure}
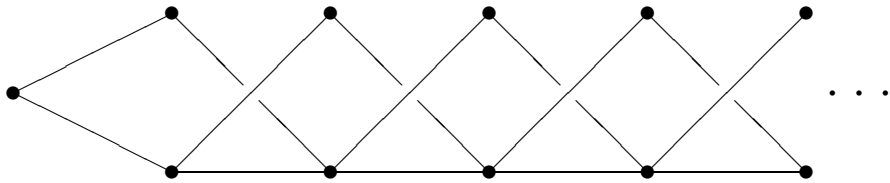

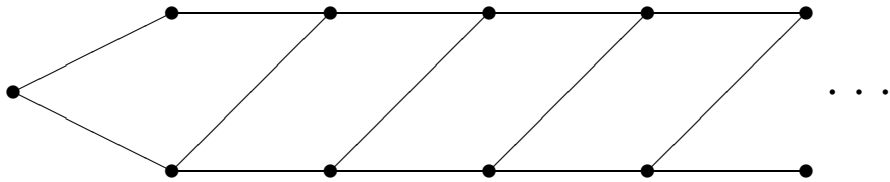
\begin{figure}[htbp]

\begin{picture}(350,70)(-10,0)

\put(0,30){\circle*{5}}

\multiput(60,0)(60,0){5}{\circle*{5}}
\multiput(60,60)(60,0){5}{\circle*{5}}

\put(0,30){\line(2,-1){60}}
\put(0,30){\line(2,1){60}}

\multiput(60,0)(60,0){4}{\line(1,1){60}}
\multiput(60,0)(60,0){4}{\line(1,0){60}}
\multiput(60,60)(60,0){4}{\line(1,0){60}}

\multiput(310,30)(10,0){3}{\circle*{2}}

\end{picture}

\caption{Bratteli diagram of the Sturmian tree $S$}
\label{fig:brat-sturm-tree}
\end{figure}

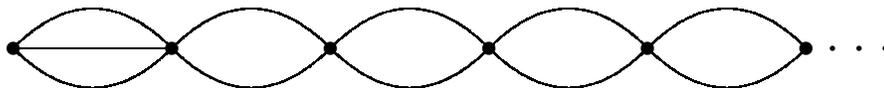
\begin{figure}[htbp]


\begin{picture}(350,70)(-10,0)

\multiput(0,30)(60,0){6}{\circle*{5}}

\multiput(0,0)(60,0){5}{\qbezier(0,30)(30,60)(60,30)}
\multiput(0,0)(60,0){5}{\qbezier(0,30)(30,0)(60,30)}

\put(0,30){\line(2,0){60}}

\multiput(310,30)(10,0){3}{\circle*{2}}

\end{picture}
\caption{Bratteli diagram of the $2$-regular tree $R_2$}
\label{fig:brat-2reg-tree}
\end{figure}

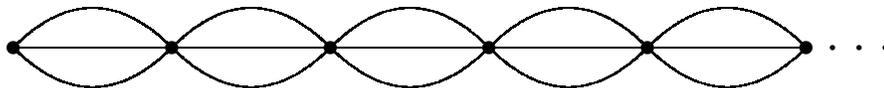
\begin{figure}[htbp]


\begin{picture}(350,70)(-10,0)

\multiput(0,30)(60,0){6}{\circle*{5}}

\multiput(0,0)(60,0){5}{\qbezier(0,30)(30,60)(60,30)}
\multiput(0,0)(60,0){5}{\qbezier(0,30)(30,0)(60,30)}

\put(0,30){\line(2,0){60}}
\multiput(0,30)(60,0){5}{\line(2,0){60}}

\multiput(310,30)(10,0){3}{\circle*{2}}

\end{picture}
\caption{Bratteli diagram of the $3$-ary tree $A_3$}
\label{fig:brat-3ary-tree}
\end{figure}

\begin{example} If $(T,v)$ denotes the Cantor tree, the Fibonacci tree,
the Sturmian tree, the $2$-regular tree $R_2$, or the $3$-ary tree
$A_3$  as defined in 
Section~\ref{subsection:tree examples}, the corresponding Bratteli
diagram $B(T,v)$ is pictured in 
Figures~\ref{fig:brat-cantor-tree} through ~\ref{fig:brat-3ary-tree}.
In each case, the initial vertex appears on the far left.
\end{example}

\begin{figure}[htbp]


\begin{picture}(350,70)(-10,0)

\multiput(0,30)(60,0){6}{\circle*{5}}

\multiput(0,0)(60,0){5}{\qbezier(0,30)(30,60)(60,30)}
\multiput(0,0)(60,0){5}{\qbezier(0,30)(30,0)(60,30)}

\multiput(310,30)(10,0){3}{\circle*{2}}

\end{picture}
\caption{Bratteli diagram of the Cantor tree $C$}
\label{fig:brat-cantor-tree}
\end{figure}
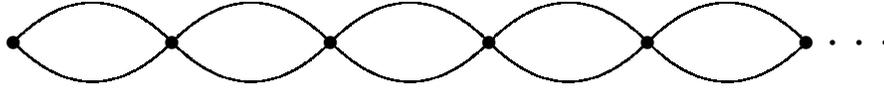

\subsection{Recollections on path groupoids}
\label{subsec:Recollection-path-groupoids}

In this section, we recall the construction of the groupoid of infinite
directed paths beginning at a distinguished vertex of
a directed graph. 
The main properties of this groupoid, due to Renault \cite{Ren},
are summarized in Theorem~\ref{theorem-Renault} below.
For more details, see
Kumjian, Pask, Raeburn, and Renault \cite{KPRR}, 
Paterson \cite{Pat}, and Renault \cite{Ren}.

We will only need the path groupoid of Bratteli diagrams, but it is just as easy to 
recall the definitions for arbitrary directed graphs.

Let $D=(\mathcal{V},\mathcal{E})$ be a directed graph 
with vertex set $\mathcal{V}$, edge set $\mathcal{E}$,
initial map $s\co  \mathcal{E}\to\mathcal{V}$, terminal map
$r\co \mathcal{E}\to\mathcal{V}$, 
and distinguished vertex $v_0$.
Assume that for each vertex $v\in\mathcal{V}$ there are at most finitely many edges $e\in\mathcal{E}$ with
initial vertex $s(e)=v$ (thus, $D$ is {\it row finite}).

A {\it path in $D$ beginning at $v_0$} is an infinite sequence 
$\alpha=(\alpha_0, \alpha_1, \alpha_2, \dots)$ of edges such that $s(\alpha_0)=v_0$ and for each $n\geq 0$, 
$r(\alpha_n)=s(\alpha_{n+1})$. Note that by convention, our paths are infinite.

The {\it path groupoid} 
$\pgd$
is the set of all pairs $(\alpha, \beta)$ of paths in $D$ beginning at $v_0$ such that $\alpha$ and
$\beta$ are {\it tail equivalent}, i.e., there exists $n\geq 0$ such that $\alpha_k=\beta_k$ for all
$k\geq n$. Tail equivalence of $\alpha$ and $\beta$ is denoted by $\alpha\sim\beta$.
The unit space is $\mathcal{P}=\mathcal{P}(D, v_0)$, the set of all paths in $D$ beginning at $v_0$.
The domain $d\co \pgd\to\mathcal{P}$ and range $r\co \pgd\to\mathcal{P}$ maps are given by
$d(\alpha,\beta)=\alpha$ and $r(\alpha,\beta)=\beta$.
Pairs $(\gamma,\delta), (\alpha, \beta)\in\pgd$ are composable if and only if $\beta=\gamma$, in which case
$(\gamma,\delta)\cdot (\alpha,\beta)=(\alpha,\delta)$.\footnote{Thus, $\pgd$ is the groupoid associated to the
equivalence relation of tail equivalence on $\mathcal{P}$.}

Observe that $\mathcal{P}$ has a natural topology; namely, consider $\mathcal{P}$ as a subspace of the countably infinite
product $\prod_0^\infty\mathcal{E}$ where $\mathcal{E}$ is given the discrete topology and the product has the
product topology.
This makes $\mathcal{P}$ a compact, totally disconnected metric space.\footnote{The topology on $\mathcal{P}$ is metrized by the ultrametric
$$d(\alpha,\beta) =  \left\{ \begin{array}{ll}
0 &\textrm{if $\alpha= \beta$}\\
e^{-n}, ~\textrm{where $n=\min\{ j ~|~ \alpha_j\not=\beta_j$\}} &\textrm{otherwise.}
\end{array} \right.$$}

For example, let $(T,v)$ be a rooted, geodesically complete, locally finite, simplicial tree
considered as a directed graph as in Section\ref{subsec:RecollectionOfBratteli}.
Then $\mathcal{P}(T,v) = end(T,v)$ as topological spaces.
In addition, we may form the Bratteli diagram $B(T,v)$ with distinguished vertex
$[v]$ associated to the tree $(T,v)$. In this case, we write $ \mathcal{P}(B(T,v),[v]) =\pb$ and
$\mathcal{PG}(B(T,v),[v])=\pgb$.

Returning to the general discussion,
we want to put a topology on $\pgd$ so that it is a locally compact groupoid with unit space $\mathcal{P}$.
It is {\it not} the subspace topology from $\mathcal{P}\times\mathcal{P}$ that we want, because that would not,
in general, be locally compact. Instead, we procede as follows. For each $n\geq 0$, define an equivalence relation
$\sim_n$ on $\mathcal{P}$ by $\alpha\sim_n\beta$ if and only if
$\alpha_k=\beta_k$ for all $k\geq n$. Let 
$$\mathcal{R}_n= \{ (\alpha,\beta)\in\mathcal{P}\times\mathcal{P} ~|~ \alpha\sim_n\beta\}.$$
Thus, $\pgd=\bigcup_{n=0}^\infty\mathcal{R}_n$. Let each $\mathcal{R}_n$ have the subspace topology
from $\mathcal{P}\times\mathcal{P}$. 
Each $\mathcal{R}_n$ is easily seen to be closed in $\mathcal{P}\times\mathcal{P}$; hence,
$\mathcal{R}_n$ is compact.
Finally, give $\pgd$ the inductive (direct) limit topology.\footnote{$U\subseteq\pgd$
is open if and only if $U\cap\mathcal{R}_n$ is open for all $n\geq 0$.}
Note that for each $n\geq 0$, the quotient space $\mathcal{P}/\sim_n$ 
is Hausdorff. In the terminology of
Exel and Lopes \cite{ExL} each $\sim_n$ is a {\it proper equivalence relation} and
tail equivalence $\sim$ is an {\it approximately proper equivalence relation}.

For a Bratteli diagram, there is the following result about the path 
groupoid. 

\begin{theorem}[Renault]
\label{theorem-Renault}
Let $D$ be a Bratteli diagram.
\begin{enumerate}
\item The path groupoid $\pgd$ is a locally compact, Hausdorff,
second countable, \'etale, AF groupoid.
\item The groupoid $C^*$-algebra
$C^*(\pgd)$ is isomorphic to $AF(D)$ as a unital $C^*$-algebra.
\end{enumerate}
\end{theorem}

These results can be found in Renault \cite{Ren};
see Exel and Renault \cite{ExR} for a recent alternative 
treatment. 
The groupoid $C^*$-algebra in the second statement 
is defined in \cite{Ren}.

The first statement in  Theorem~\ref{theorem-Renault} above, combined with
Theorem~\ref{theorem:bratteli-one} below, 
gives another way of establishing the first two
statements in Theorem~\ref{theorem:first-main}.

\subsection{Theorems on path groupoids of Bratteli diagrams}
\label{subsec:the-theorems}

The main result of this section is the following theorem.
It concerns locally rigid end spaces of trees and relates their
groupoids of local isometries to the path groupoids of the Bratteli
diagram associated to the tree.

\begin{theorem}
\label{theorem:bratteli-one}
Let $(T,v)$ be a rooted, geodesically complete, locally finite, simplicial tree.
If $end(T,v)=X$ is locally rigid, then the quotient map 
$\kappa\co T\to B(T,v)$ induces an isomorphism of groupoids
$$\kappa_*\co \ligx\to\mathcal{PG}(B(T,v)).$$
\end{theorem}

\begin{proof}
We first show that the quotient map $\kappa\co T\to B(T,v)$ induces a homeomorphism
$\kappa_\#\co  end(T,v)=X\to \pb=\mathcal{P}$ 
between unit spaces of the groupoids
$\ligx$ and $\pgb$.
This part of the proof does not use the local rigidity hypothesis.

Define $\kappa_\#$ as follows. Represent $x\in end(T,v)$ 
(which is a geodesic ray $x\co [0,\infty)\to T$ with $x(0)=v$)
by an infinite sequence
of edges
$(x_0, x_1, x_2, \dots)$ of $T$. That is,
$x_i= x([i, i+1])$ for all $i=0,1,2,\dots$.
Then set $\kappa_\#x= (\kappa x_0, \kappa x_1, \kappa x_2, \dots)\in\mathcal{P}$.

Note the following simple fact about the map $\kappa$.

\begin{fact}
\label{fact:kappa}
If $e_1,e_2$ are edges of $T$ such that $e_1\not= e_2$ and $s(e_1)=s(e_2)$, then
$\kappa(e_1)\not=\kappa(e_2)$ as edges of $B(T,v)$.
\end{fact}

This is true because otherwise $\alpha(T_l^i, T_{s(e_1)})^{-1}(e_1)=
\alpha(T_l^i, T_{s(e_1)})^{-1}(e_2)$, where $l$ is the integer such that $T_i^l$ is rooted isometric to 
$T_{s(e_1)} = T_{s(e_2)}$, contradicting the fact that $\alpha(T_l^i, T_{s(e_1)})$ is an isometry.

From this fact it follows that $\kappa$ is a local homeomorphism in the sense that for all $t\in T$ there
exists an open neighborhood $U_t$ of $t$ in $T$ such that $\kappa|\co U_t\to\kappa(U_t)$ is a
homeomorphism. (However, $\kappa(U_t)$ need not be open in $B(T,v)$ because there might be 
edges $e_1\not= e_2$ in $T$ 
which go to edges in $B(T,v)$ 
with the same terminal vertices $r(\kappa(e_1)) = r(\kappa(e_2))$.)

To see that $\kappa_\#$ is surjective, suppose $\alpha=(\alpha_0, \alpha_1, \alpha_2, \dots)\in\mathcal{P}$.
Choose an edge $e_0$ in $T$ beginning at $v$ such that $\kappa(e_0)=\alpha_0$.
Then $\kappa(r(e_0))= r(\alpha_0)=s(\alpha_1)$. By the local homeomorphism property of $\kappa$ mentioned 
above, there exists an edge $e_1$ in $T$ beginning at $r(e_0)$ such that $\kappa(e_1)=\alpha_1$.
Continue this process to construct $x=(e_0, e_1, e_2, \dots)\in end(T,v)$
such that $\kappa_\#x=\alpha$.

To see that $\kappa_\#$ is injective, suppose $x\not= y$ in $end(T,v)$, and let $t_0=\sup\{ t\geq 0 ~|~ x(t)=y(t)\}$.
Then $e_1=x([t_0,t_0+1])$ and $e_2=y([t_0, t_0+1])$ are distinct edges of $T$ with $s(e_1)= s(e_2)$.
It follows from Fact~\ref{fact:kappa} that $\kappa(e_1)\not=\kappa(e_2)$. It follows that
$\kappa_\#x\not=\kappa_\#y$.

Moreover,
\begin{eqnarray*}
\lefteqn{t_0=\min\{ j\in\{ 0,1,2,\dots\} ~|~ x(j+1)\not= y(j+1)\}}\\
& &
=\min\{ j\in\{ 0,1,2,\dots\} ~|~ \kappa(x([j,j+1]))\not= \kappa(y([j,j+1]))\}.
\end{eqnarray*}
Since $d_e(x,y)= e^{-t_0}$, it follows that $\kappa_\#$ is an isometry with respect to the natural metric on 
$\mathcal{P}$. This completes the proof that $\kappa_\#$ is a homeomorphism.

In order to define $\kappa_*\co \ligx\to\pgb$, choose $\epsilon_X>0$ by Lemma~\ref{lemma:compact-locally-rigid}.
This uses the local rigidity assumption on $X$.
Represent a given groupoid element 
$[g,x]\in\ligx$ by an isometry $g\co B(x,\epsilon)\to B(gx,\epsilon)$ with 
$0 < \epsilon\leq\epsilon_X$.

\begin{claim}
\label{claim:kappa-tail}
$\kappa_\# x$ and $\kappa_\# gx$ are tail equivalent.
\end{claim} 

\begin{proof*}\emph{of Claim.}
Since 
$$B(x,\epsilon)=end(T_{\langle x,\epsilon\rangle}, v_{\langle x,\epsilon\rangle}) \text{ and }
B(gx,\epsilon)= end(T_{\langle gx,\epsilon\rangle}, v_{\langle gx,\epsilon\rangle}),$$
(see Section\ref{section:trees}) the isometry $g$ induces a rooted isometry
$$\tilde g\co  (T_{\langle x,\epsilon\rangle}, v_{\langle x,\epsilon\rangle})\to
(T_{\langle gx,\epsilon\rangle}, v_{\langle gx,\epsilon\rangle}).$$ This map $\tilde g$ is defined on edges as follows:
if $e$ is an edge of $T_{\langle x,\epsilon\rangle}$, choose $y\in end(T,v)$ such that $e=y([i,i+1])$ for some $i$.
Then $y\in B(x,\epsilon)$ and so $gy\in B(gx,\epsilon)$. Thus,
$(gy)([i,i+1])$ is an edge of $T_{\langle gx,\epsilon\rangle}$ and we set
$\tilde g e=(gy)([i,i+1])$.
Isometries $B(x,\epsilon)\to B(gx,\epsilon)$ correspond bijectively
to rooted isometries 
$(T_{\langle x,\epsilon\rangle}, v_{\langle x,\epsilon\rangle})\to
(T_{\langle gx,\epsilon\rangle}, v_{\langle gx,\epsilon\rangle})$ (e.g., see \cite{Hug}). 
Since $\epsilon\leq\epsilon_X$, $g$ is the unique isometry from $B(x,\epsilon)$ to $B(gx,\epsilon)$.
Thus, $\tilde g$ is the unique rooted isometry.
It follows that if $T_l^i$ is the admissible level $i$ 
subtree of $T$ that is rooted isometric to $T_{\langle x,\epsilon\rangle}$,
and $\alpha_1\co T_l^i\to T_{\langle x,\epsilon\rangle}$ and 
$\alpha_2\co T_l^i\to T_{\langle gx,\epsilon\rangle}$, then $\tilde g\alpha_1=\alpha_2$.
From the definition of $\tilde g$, it follows that if $i\geq \lceil -\ln\epsilon\rceil$,
then $\tilde g(x([i,i+1]))=(gx)([i,i+1])$. Thus,
$x([i,i+1])\sim (gx)([i,i+1])$ for $i\geq\lceil -\ln \epsilon\rceil$.
Thus, $\kappa_\# x$ and $\kappa_\# gx$ are tail equivalent.
\end{proof*}

Thus, define 
$$\kappa_\ast([g,x])=(\kappa_\# x,\kappa_\# gx)\in\mathcal{P}\times\mathcal{P}.$$
The claim implies that
$(\kappa_\# x,\kappa_\# gx)\in\pgb$.

Note that $\kappa_\ast([g,x])$ is well-defined in the sense that it does not depend
on the isometry representing $[g,x]$, only on the germ of the isometry at $x$ (in fact, a 
feature of local rigidity is that $\kappa_\ast([g,x])$ only depends on $x$ and $gx$).

Note also that the diagram
$$\begin{CD}
X @>{\kappa_\#}>> \mathcal{P}\\
@V{\alpha}VV @VV{\Delta}V\\
\ligx @>{\kappa_\ast}>> \pgb
\end{CD}$$
commutes, where the vertical maps are the natural inclusions of unit spaces ($\alpha$ is given
in Remark~\ref{remark:unit-space}
and $\Delta$ is the diagonal map $\Delta(\beta,\beta)$).	

It is equally obvious that the diagrams 
$$\begin{CD}
\ligx @>{\kappa_\ast}>> \pgb\\
@V{d}VV @VV{d}V\\
X @>{\kappa_\#}>> \mathcal{P}
\end{CD}
~~~\text{and}~~~
\begin{CD}
\ligx @>{\kappa_\ast}>> \pgb\\
@V{r}VV @VV{r}V\\
X @>{\kappa_\#}>> \mathcal{P}
\end{CD}
$$
commute.

To see that $\kappa_\ast$ is multiplicative, suppose 
$[g_1,x_1], [g_2,x_2]\in\ligx$ with $x_2=g_1x_1$. Then 
$\kappa_\ast([g_2,x_2]\cdot [g_1,x_1])=
\kappa_\ast([g_2g_1,x_1])=
(\kappa_\#x_1,\kappa_\#g_2g_1x_1)=
(\kappa_\#g_1x_1, \kappa_\#g_2g_1x_1)\cdot(\kappa_\#x_1,\kappa_\#g_1x_1)=
\kappa_\ast([g_2,g_1x_1])\cdot \kappa_\ast([g_1,x_1])$.

It only remains to show that $\kappa_\ast\co \ligx\to\pgb$ is a homeomorphism.
Let $\epsilon_X>0$ be given by Lemma~\ref{lemma:compact-locally-rigid} and
choose a positive integer $N\geq -\ln\epsilon_X$. For $i=0,1,2,3,\dots$, let
$\epsilon_i=e^{-(N+i)}$.

As in the proof of Theorem~\ref{theorem:AF}, $\ligx$ is the union of approximating groupoids,
$\ligx =\cup_{i=0}^\infty\mathcal{G}_{LI}^{\epsilon_i}(X)$.
Elements of $\mathcal{G}_{LI}^{\epsilon_i}(X)$ are written $(g,x)$, but when considered
in $\ligx$, are written $[g,x]$.

\begin{claim}
\label{claim:restrict-homeo}
For every $i\geq 0$, $\kappa_\ast$ restricts to a homeomorphism
$$\kappa_\ast|\co \mathcal{G}_{LI}^{\epsilon_i}(X)\to\mathcal{R}_{N+i}.$$
\end{claim}

\begin{proof}
If $(g,x)\in\mathcal{G}_{LI}^{\epsilon_i}(X)$, then $g\co B(x,\epsilon_i)\to B(gx,\epsilon_i)$
is an isometry. We first need to observe that
$(\kappa_\#x,\kappa_\#gx)\in\mathcal{R}_{N+i}$. 
According to the proof of Claim~\ref{claim:kappa-tail},
$x([j,j+1])\sim(gx)([j,j+1])$ for
$j\geq\lceil-\ln\epsilon_i\rceil$.
Since $-\ln\epsilon_i=N+i$, we have the desired observation.
This shows $\kappa_\ast(g,x)\in\mathcal{R}_{N+i}$.

To see that $\kappa_*|$ is continuous,
let $(g,x)\in\mathcal{G}_{LI}^{\epsilon_i}(X)$ and let $\epsilon>0$ be given. We may assume $\epsilon\leq\epsilon_i$.
Recall from Section\ref{section:TheApproximatingGroupoids} that
$(g,x)$ has an open neighborhood $U_\epsilon(g,x)$ in $\mathcal{G}_{LI}^{\epsilon_i}(X)$
given by
$$U_\epsilon(g,x) = \{ (h,y)\in\mathcal{G}_{LI}^{\epsilon_i}(X) ~|~
d(x,y)<\epsilon \text{ and } d(gz,hz)<\epsilon \text{ for every } z\in B(x,\epsilon_i\}.$$ 
Let $(h,y)\in U_\epsilon(g,x)$.
Then $d(x,y)<\epsilon$ and $d(gy,hy)<\epsilon$.
Since $g$ is an isometry, $d(gx,gy)<\epsilon$.
Hence, $d(x,y)<\epsilon$ and $d(gx,hy)<2\epsilon$.
Since it was shown above that $\kappa_\#$ is an isometry, we have
$d(\kappa_\# x,\kappa_\# y)<\epsilon$ and
$d(\kappa_\# gx,\kappa_\# hy)<2\epsilon$ in $\mathcal{P}$.
Thus, the distance between
$\kappa_*(g,x) = (\kappa_\# x,\kappa_\# gx)$ and $\kappa_*(h,y)=
(\kappa_\# y,\kappa_\# hy)$ in $\mathcal{P}\times\mathcal{P}$ is small
if $\epsilon$ is small enough.
Since $\mathcal{R}_{N+i}$ is topologized as a subspace of
$\mathcal{P}\times\mathcal{P}$, this verifies
the continuity of $\kappa_*|$.

To see that $\kappa_*|$ is injective, suppose
$(g,x), (h,y) \in \mathcal{G}_{LI}^{\epsilon_i}(X)$.
Since $\kappa_*(g,x)= (\kappa_\# x,\kappa_\# gx)$,
$\kappa_*(h,y) = (\kappa_\# y,\kappa_\# hy)$, and 
$\kappa_\#$ is injective, it follows that $\kappa_*(g,x) = \kappa_*(h,y)$
implies that $x=y$ and $gx=hy$. By the choice of $\epsilon_X$, it follows that $g=h$.

To see that $\kappa_*|$ is surjective, let $(\alpha,\beta)\in\mathcal{R}_{N+i}$ 
be given.
Choose $x,y\in X$ such that $\kappa_\# x=\alpha$ and $\kappa_\# y=\beta$.
It suffices to show there exists an isometry
$g\co B(x,\epsilon_i)\to B(y,\epsilon_i)$
such that $gx=y$; for then, $(g,x)\in 
\mathcal{G}_{LI}^{\epsilon_i}(X)$
and $\kappa_*(g,x)=(\alpha,\beta)$.
Since $(\alpha,\beta)\in\mathcal{R}_{N+i}$,
it follows that $\alpha_k=\beta_k$ for all $k\geq N+i$.
Denote the sequence of edges of $x$ by
$(x_0,x_1,x_2,\dots)$ and those of $y$ by
$(y_0,y_1,y_2,\dots)$.
That is, $x_k=x[k,k+1]$ and $y_k=y[k,k+1]$,
where $x,y\co [0,\infty)\to T$.

If $x_k=y+k$ for all $k$, then we are done because we can take $g=\id$.
Therefore, assume this is not the case and let
$M=\min\{ k ~|~ x_k\neq y_k\}$.

Since $x_M\neq y_M$ and $s(x_M)=x(M) =y(M)=s(y_M)$, Fact~\ref{fact:kappa}
implies $\alpha_M=\kappa(x_M)\neq\kappa(y_M)=\beta_M$.
Thus, $M<N+i$.

Since $\alpha_{N+i}=\beta_{N+i}$, it follows that the tree
$T_{x(N+i)}$ is rooted isometric to $T_{y(N+i)}$.
Let $\ell$ be the unique integer with $1\leq \ell\leq m_{N+i}$ such that
$T^{N+i}_\ell$ is rooted isometric to
$T_{x(N+i)}$ and $T_{y(N+i)}$.
Consider the following rooted isometry defined as a composition of an admissible
isometry and the inverse of an admissible isometry:
$$\tilde g := \alpha(T^{N+i}_\ell, T_{y(N+i)})\circ \alpha(T^{N+i}_\ell, T_{x(N+i)})^{-1}\co 
T_{x(N+i)} \to T_{y(N+i)}.$$
Thus, $\tilde g(x_{N+i})=y_{N+i}$.

It follows that $\tilde g$ induces an isometry
$g\co B(x,\epsilon_i)\to B(y,\epsilon_i)$.
(According to Section\ref{section:trees},
there is an isometry between closed balls, but it
restricts to an isometry between open balls.)
To show $gx=y$, it suffices to show that
$\tilde g(x_k)=y_k$ for all $k\geq N+i$. Assume to the contrary
that this is not the case, and let
$K=\min\{ k\geq N+i ~|~ \tilde g(x+k)\neq y+k\}$.
Thus, $K > N+i$.

Proceed as above: since $\alpha_K=\beta_K$, it follows that 
$T_{x(K)}$ is rooted isometric to $T_{y(K)}$.
Moreover, there exists a rooted isometry
$\tilde g_K\co T_{x(K)} \to T_{y(K)}$ such that
$\tilde g_K(x_K)=y_K$, and 
$\tilde g_K$ induces an isometry
$g_K\co \overline B(x,e^{-K})\to\overline B(y,e^{-K})$.
It follows that $g(\overline B(x,e^{-K}))=\overline B(y,e^{-K})$
(because $g(x)\in\overline B(y,e^{-K})$).
Hence,
$g_K^{-1}\circ g|\co \overline B(x,e^{-K})\to\overline B(x,e^{-K})$
is an isometry with $g_K^{-1}\circ g(x)\neq x$. 
This contradicts the choice of $\epsilon_X$.
For this, we need to observe that
$\overline B(x, e^{-K})= B(x,\eta)$ if 
$e^{-K} <\eta < e^{-(K+1)}$,
and any such $\eta$ satisfies $\eta <e^{-(K+1)}\leq\epsilon_X$.
Hence, $gx=y$ and $\kappa_*|$ is surjective.

Since 
$\mathcal{G}_{LI}^{\epsilon_i}(X)$
is compact Hausdorff by Remark~\ref{remark:approx},
this shows that $\kappa_*|$ is a homeomorphism and
completes  the proof of Claim~\ref{claim:restrict-homeo}.
\end{proof}

Note that the following diagram commutes, where the vertical arrows
are inclusion maps:
$$\begin{CD}
\mathcal{G}_{LI}^{\epsilon_i}(X) @>{\kappa_\ast|}>> \mathcal{R}_{N+i}\\
@VVV @VVV\\
\mathcal{G}_{LI}^{\epsilon_{i+1}}(X)@>{\kappa_\ast|}>> \mathcal{R}_{N+i+1}
\end{CD}$$
Recall from the proof of Theorem~\ref{theorem:AF} that $\ligx$ is the inductive limit 
(as $i\to\infty$) of the left-hand vertical maps. By definition,
$\mathcal{R}$ is the inductive limit of the right-hand vertical maps.
Hence, $\kappa_*$ is a homeomorphism.

This completes the proof of Theorem~\ref{theorem:bratteli-one}.
\end{proof}

\subsection{Summary of 
Section~\ref{section:TreesBratteliDiagramsPathGroupoids}}
\label{subsec:summary}

Let $(T,v)$ be a rooted, geodesically complete, locally finite
simplicial tree and let $X=end(T,v)$.
By examining isomorphic subtrees of $T$ rooted at the same level
of $T$, we defined a Bratteli diagram $B(T,v)$ that is a quotient of
$T$, $\kappa\co T\to B(T,v)$.

For any Bratteli diagram (in fact, for any rooted directed graph)
there is a well-known construction of a groupoid, called
the path groupoid, based on tail equivalence
of infinite directed paths beginning at the distinguished vertex of the
diagram. In our case, we denote the path groupoid of $B(T,v)$
by $\pgb$. This groupoid satisfies sufficient conditions so that 
Renault's theory can be applied to obtain a unital AF $C^*$-algebra
$C^*\pgb$.

On the other hand, Bratteli showed how to construct a unital AF
$C^*$-algebra from any Bratteli diagram. For the Bratteli diagram
$B(T,v)$, this algebra is denoted by $\mathrm{AF}(B(T,v))$.

It is well-known that Bratteli's construction and Renault's theory
lead to isomorphic unital $C^*$-algebras. In particular, there is an
isomorphism 
$$C^*\pgb\cong\mathrm{AF}(B(T,v))$$ 
of unital $C^*$-algebras.

A unital partially ordered abelian group is obtained from $(T,v)$
in two ways. First, we take the unital, ordered $K_0$-group of
a unital $C^*$-algebra and
get
$$(K_0(C^*\pgb), K_0(C^*\pgb)_+, [1]).$$
Second, there is the unital dimension group associated to a Bratteli
diagram. In particular, we get
$(G(B(T,v)), G_+(B(T,v)), [1])$.
Since this is the unital ordered $K_0$-group of  $\mathrm{AF}(B(T,v))$,
these two constructions lead to isomorphic unital partially ordered
abelian groups. In particular,
there is an isomorphism
$$(K_0(C^*\pgb), K_0(C^*\pgb)_+, [1])\cong
(G(B(T,v)), G_+(B(T,v)), [1])$$
of unital partially ordered abelian groups.

These constructions are summarized in the following diagram:
\begin{eqnarray*}
(T,v) &\leadsto & B(T,v) \\
&\leadsto &\pgb \\
& \leadsto & C^*\pgb\cong\mathrm{AF}(B(T,v))  \\
& \leadsto & (K_0(C^*\pgb), K_0(C^*\pgb)_+, [1]) \\
& \cong & (K_0(\mathrm{AF}(B(T,v)), (K_0(\mathrm{AF}(B(T,v))_+, [1]) \\
& \cong & (G(B(T,v)), G_+(B(T,v)), [1]) \\
\end{eqnarray*}

Under the assumption that $X=end(T,v)$ is locally rigid, there
was another route that led to groupoids, unital AF $C^*$-algebras, and
unital partially ordered abelian groups.
Namely, we formed the groupoid $\ligx$ of local isometries on $X$ and
verified sufficient conditions so that Renault's theory produces
a unital AF $C^*$-algebra $C^*\ligx$. We can take the $K_0$-group
of that $C^*$-algebra and get a unital partially ordered abelian
group. This route is summarized by the following diagram:

\begin{eqnarray*}
(T,v) &\leadsto & end(T,v) = X\\
& \leadsto & \ligx \\
& \leadsto & C^*\ligx \\
& \leadsto & (K_0(C^*\ligx), K_0(C^*\ligx)_+, [1]) \\
\end{eqnarray*}

The main isomorphisms established in this section in the 
locally rigid case are summarized in the following corollary.

\begin{corollary}
\label{cor:summary}
If $(T,v)$ is a rooted, geodesically complete, locally finite
simplicial tree such that $X=end(T,v)$ 
locally rigid,
then
\begin{enumerate}
\item there is an isomorphism of 
topological groupoids $\ligx \cong \pgb$,
\item there are isomorphisms of unital $C^*$-algebras 
\begin{eqnarray*}
C^*\ligx &\cong & \mathrm{AF}(B(T,v)) \\
& \cong & C^*\pgb,
\end{eqnarray*}
\item there are isomorphisms of unital partially ordered abelian
groups
\begin{eqnarray*}
\lefteqn{(K_0(C^*\ligx), K_0(C^*\ligx)_+, [1])}\\
& \qquad\qquad\cong & (K_0(C^*\pgb), K_0(C^*\pgb)_+, [1]) \\
& \qquad\qquad\cong & 
(K_0(\mathrm{AF}(B(T,v)), (K_0(\mathrm{AF}(B(T,v))_+, [1]) \\
& \qquad\qquad\cong & (G(B(T,v)), G_+(B(T,v)), [1]). \\
\end{eqnarray*}
\end{enumerate}
In particular, $B(T,v)$ is the Bratteli diagram
for the unital AF algebra $C^*\ligx$.
\end{corollary}

Here is one last consequence of the results of this section.

\begin{corollary}
If $(T,v)$ and $(S,w)$ are rooted, geodesically complete, locally finite
simplicial trees such that $X=end(T,v)$ and $Y= end(S,w)$ 
are locally rigid,
then the Bratteli diagrams
$B(T,v)$ and $B(S,w)$ are equivalent if and only if 
$K_0C^*\ligx$ and $K_0C^*\ligy$ are isomorphic as unital
partially ordered abelian groups.
\end{corollary}

\begin{proof} By Theorem~\ref{theorem:Bratteli, Elliott}, 
$B(T,v)$ and $B(S,w)$ are equivalent if and only if 
$G(B(T,v))$ and $G(B(S,w))$ are isomorphic as unital
partially ordered abelian groups. The result now follows from the
isomorphisms above. 
\end{proof}


\section{The symmetry at infinity group}
\label{section:TheSymmetryAtInfinityGroup}

Let $(T,v)$ be a rooted, geodesically complete, locally finite simplicial 
tree. In Section~\ref{section:TreesBratteliDiagramsPathGroupoids}
we defined a Bratteli diagram $B(T,v)$ associated to $(T,v)$ by
identifying infinite subtrees of $T$ that occur at the same level.
The Bratteli diagram $B(T,v)$ leads to three isomorphic
unital partially ordered
abelian groups (as is the case for all 
Bratteli diagrams)\footnote{If $X=end(T,v)$ is locally rigid, 
then $(K_0(C^*\ligx), K_0(C^*\ligx)_+, [1])$
is a fourth group isomorphic
to these.} :
\begin{eqnarray*}
\lefteqn{(K_0(C^*\pgb), K_0(C^*\pgb)_+, [1])} \\
& \qquad\qquad\cong & 
(K_0(\mathrm{AF}(B(T,v)), (K_0(\mathrm{AF}(B(T,v))_+, [1]) \\
& \qquad\qquad\cong & (G(B(T,v)), G_+(B(T,v)), [1]) \\
\end{eqnarray*}

Of these, 
$(G(B(T,v)), G_+(B(T,v)), [1])$ is obviously the most directly defined.

We want to point out in this section that this unital partially
ordered abelian group can be defined
even more directly from the tree $(T,v)$ without passing to the 
Bratteli diagram $B(T,v)$. Of course, this is implicit in
Section~\ref{section:TreesBratteliDiagramsPathGroupoids},
but we want to make clear just how elementary the idea is.
Some examples are included at the end of this section.

Since the group is defined directly from the tree $(T,v)$ and it measures
the symmetries at infinity of $(T,v)$---that is, the number of isometric
infinite subtrees of $(T,v)$---we denote it by
$\symtv$.

We use the notation from Section~\ref{subsec:ConstructionOfBTv}.
Recall that
for each $i=0,1,2,\dots$,
$m_i$ is the number of rooted isometry classes of level $i$ subtrees of 
$(T,v)$. We choose level $i$ subtrees
$T_1^i, T_2^i,\dots T_{m_i}^i$ that form a complete
set of representatives of the rooted isometry classes.
In particular, $m_0=1$ and $T_1^0=T$.

For $0\leq i$, $1\leq k\leq m_{i+1}$ and $1\leq \ell \leq m_i$, 
let $a_{k\ell}^i$ be the number of level one subtrees of $T_\ell^i$ that are rooted isometric to
$T_k^{i+1}$. The matrix $A_i=[a_{k\ell}^i]$ is an $(m_{i+1}\times m_i)-$ matrix with nonnegative
integral entries. The indexing of the rows and columns is indicated here:

\bigskip
\begin{gather*}
\begin{matrix}
T_1^i & ~~~~\cdots & ~~~T_\ell^i & ~~\cdots & ~~~T_{m_i}^i\\
\end{matrix}\\
A_i = \quad
\begin{bmatrix} 
a_{11}^i & \cdots & a_{1\ell}^i & \cdots & a_{1m_i}^i\\
\vdots &           & \vdots     &       & \vdots \\  
a_{k1}^i & \cdots &a_{k\ell}^i & \cdots & a_{km_i}^i \\
\vdots &           & \vdots     &       & \vdots \\ 
a_{m_{i+1}1}^i & \cdots & a_{m_{i+1}\ell}^i & \cdots & a_{m_{i+1}m_i}^i\\
\end{bmatrix} 
\begin{matrix}
T_1^{i+1} \\ \vdots \\ T_k^{i+1} \\ \vdots \\ T_{m_i+1}^{i+1}
\end{matrix}
\end{gather*}
\bigskip

The simplicial cone of $\bz^n$ is $\bz_+^n=\{ (x_1,x_2,\dots, x_n)\in\bz^n ~|~ x_i\geq 0\}$.
It is a subsemigroup of $\bz^n$ and the resulting partial order on $\bz^n$ is called the
{\it simplicial ordering.\/}

\begin{definition} For a rooted, locally finite simplicial tree $(T,v)$, using 
the notation above, the {\it symmetry at infinity group $\symtv$\/} is the
unital partially
ordered abelian group given by the direct limit
$$\symtv =\lim_\to (\bz \stackrel{A_0}{\longrightarrow}
\bz^{m_1} \stackrel{A_1}{\longrightarrow} \bz^{m_2} \stackrel{A_2}{\longrightarrow}\bz^{m_3} 
\stackrel{A_3}{\longrightarrow} \cdots
\bz^{m_i} \stackrel{A_i}{\longrightarrow}\cdots).$$
The {\it positive cone of $\symtv$\/} is the subsemigroup given by the 
direct limit
$$\symtv_+ =\lim_\to (\bz_+ \stackrel{A_0}{\longrightarrow}
\bz_+^{m_1} \stackrel{A_1}{\longrightarrow}\bz_+^{m_2} \stackrel{A_2}{\longrightarrow}
\bz_+^{m_3} \stackrel{A_3}{\longrightarrow}\cdots
\bz_+^{m_i} \stackrel{A_i}{\longrightarrow}\cdots).$$
The order unit $[1]$ of $\symtv$ is the 
class of $1\in\bz$ in the direct limit.
\end{definition}

Thus, we simply use $\symtv$ to denote the unital partially ordered
abelian group given by the triple
$(\symtv,\symtv_+, [1])$.

\begin{proposition} 
\label{prop:dim gp}
If $(T,v)$ is a rooted, geodesically complete, 
locally finite simplicial tree,
then $\symtv$ is isomorphic to the unital dimension group of 
the Bratteli diagram $B(T,v)$; that is, there is an isomorphism
of unital partially ordered abelian groups
$$\symtv\cong (G(B(T,v)), G_+(B(T,v)), [1]).$$
\end{proposition}

\begin{proof}
We continue to use the notation from 
Section~\ref{subsec:ConstructionOfBTv}.
Recall that the Bratteli diagram $B(T,v) = (\mathcal{V}, \mathcal{E})$ with
$\mathcal{V} = \cup_{i=0}^\infty\mathcal{V}_i$ and
$\mathcal{V}_i=\{[v_1^i],\dots, [v_{m_i}^i]\}$ where
$v_\ell^i$ is the root of $T_\ell^i$.
It follows from Remark~\ref{remark:explicit}
that the number  $a_{k\ell}^i$ defined above is exactly the number of edges
in $\mathcal E$ from $[v_\ell^i]$ to $[v_k^{i+1}]$. 
This gives the required isomorphism of  unital partially ordered abelian
groups.
\end{proof}

The following corollary is 
Theorem\ref{theorem:first-main}\ref{first-main-four} of the introduction.

\begin{corollary}
\label{sym-and-k0}
If $X$ is a locally rigid, compact ultrametric space and
$X=end(T,v)$, where $(T,v)$ is a rooted, geodesically complete,
locally finite simplicial tree, then $\symtv$ is isomorphic to
$K_0C^*\ligx$ as a unital partially ordered abelian group.
\end{corollary}

\begin{proof}
This follows from Proposition~\ref{prop:dim gp} and
Corollary~\ref{cor:summary}.
\end{proof}

%

We now give the symmetry at infinity groups of the trees in
Section~\ref{subsection:tree examples}.
All of the calculations are elementary and are well-known (perhaps
in other contexts).
For each of the trees, the natural root is denoted $v$.

\begin{example}{\bf (The Cantor tree $C$)~} $Sym_\infty(C,v)$ is 
isomorphic to the additive group of dyadic rationals
$\bz[\frac 1 2] = \{ \frac m {2^i} ~|~ m,i\in\bz\}\subseteq\bq$ 
with positive cone $Sym_\infty(C,v)_+$ corresponding to 
the nonnegative dyadic rationals
$\bz[\frac 1 2 ]_+$ and order unit $1\in\bz[\frac 1 2]$.
The direct sequence  is
$$\bz \stackrel{2}{\longrightarrow}
\bz \stackrel{2}{\longrightarrow} \bz \stackrel{2}{\longrightarrow}
\bz 
\stackrel{2}{\longrightarrow} \cdots \stackrel{2}{\longrightarrow}
\bz \stackrel{2}{\longrightarrow}\cdots.$$
\end{example}

\begin{example}{\bf (The Fibonacci tree $F$)~}
$Sym_\infty(F,v)$ is isomorphic to the two-dimensional integral lattice 
$\bz^2$ with positive cone $Sym_\infty(F,v)_+$ corresponding to 
$\{ (x,y)\in\bz^2 ~|~
\tau x+y\geq 0\}$ where $\tau=\frac {1+\sqrt{5}} 2$ is the golden mean.
The order unit is $(1,1)\in\bz^2$.
Thus, $Sym_\infty(F,v)\cong (\bz+\tau\bz, 
(\bz+\tau\bz)\cap\br_+, \tau)$.
The direct sequence  is
$$\bz
\xrightarrow {\tiny\begin{bmatrix}  1 \\ 1 \end{bmatrix}}
\bz^2 
\xrightarrow {\tiny\begin{bmatrix} 1 ~ 1\\ 1 ~ 0\end{bmatrix}}
\bz^2
\xrightarrow {\tiny\begin{bmatrix} 1 ~ 1\\ 1 ~ 0\end{bmatrix}}
\bz^2 
\xrightarrow {\tiny\begin{bmatrix} 1 ~ 1\\ 1 ~ 0\end{bmatrix}}
\cdots 
\xrightarrow {\tiny\begin{bmatrix} 1 ~ 1\\ 1 ~ 0\end{bmatrix}}
\bz^2 
\xrightarrow {\tiny\begin{bmatrix} 1 ~ 1\\ 1 ~ 0\end{bmatrix}}
\cdots.$$
\end{example}

\begin{example}{\bf (The Sturmian tree $S$)~} 
$Sym_\infty(S,v)$ is isomorphic to the two-dimensional integral lattice 
$\bz^2$ with positive cone $Sym_\infty(S,v)_+$ corresponding to 
$\{ (x,y)\in\bz^2 ~|~
x,y\geq 0 \textrm{ or } x>0\}$. This is the lexiographic order
of $\bz^2$. The order unit is $(1,1)\in\bz^2$.
The direct sequence  is
$$\bz 
\xrightarrow {\tiny\begin{bmatrix} 1 \\ 1 \end{bmatrix}}
\bz^2 
\xrightarrow {\tiny\begin{bmatrix} 1 ~ 0\\ 1 ~ 1\end{bmatrix}}
\bz^2
\xrightarrow {\tiny\begin{bmatrix} 1 ~ 0\\ 1 ~ 1\end{bmatrix}}
\bz^2 
\xrightarrow {\tiny\begin{bmatrix} 1 ~ 0\\ 1 ~ 1\end{bmatrix}}
\cdots 
\xrightarrow {\tiny\begin{bmatrix} 1 ~ 0\\ 1 ~ 1\end{bmatrix}}
\bz^2 
\xrightarrow {\tiny\begin{bmatrix} 1 ~ 0\\ 1 ~ 1\end{bmatrix}}
\cdots.$$
\end{example}

\begin{example}{\bf (The $n$-regular tree $R_n$)~} $Sym_\infty(R_n,v)$ is 
isomorphic to the additive group
$\bz[\frac 1 n] = \{ \frac m {n^i} ~|~ m,i\in\bz\}\subseteq\bq$ 
with positive cone $Sym_\infty(R_n,v)_+$ corresponding to 
the nonnegative elements of
$\bz[\frac 1 n ]_+$ and order unit $n+1\in\bz[\frac 1 n]$.
The direct sequence  is
$$\bz 
\xrightarrow {n+1}
\bz \stackrel{n}{\longrightarrow} \bz \stackrel{n}{\longrightarrow}
\bz 
\stackrel{n}{\longrightarrow} \cdots \stackrel{n}{\longrightarrow}
\bz \stackrel{n}{\longrightarrow}\cdots.$$
\end{example}

\begin{example}{\bf (The $n$-ary tree $A_n$)~} $Sym_\infty(A_n,v)$ is 
isomorphic to 
$\bz[\frac 1 n]$ 
with positive cone $Sym_\infty(A_n,v)_+$ corresponding to 
$\bz[\frac 1 n ]_+$ and order unit $1\in\bz[\frac 1 n]$.
The direct sequence  is
$$\bz 
\stackrel{n}{\longrightarrow}
\bz \stackrel{n}{\longrightarrow} \bz \stackrel{n}{\longrightarrow}
\bz 
\stackrel{n}{\longrightarrow} \cdots \stackrel{n}{\longrightarrow}
\bz \stackrel{n}{\longrightarrow}\cdots.$$
\end{example}

\begin{example}{\bf (The $n$-ended tree $E_n$)~} 
$Sym_\infty(E_n,v)$ is 
isomorphic to $\bz$ with $Sym_\infty(E_n,v)_+$ corresponding to 
the nonnegative integers $\bz_+$ and order unit $n\in\bz$.
The direct sequence  is
$$\bz 
\stackrel{n}{\longrightarrow}
\bz 
\stackrel{1}{\longrightarrow} 
\bz 
\stackrel{1}{\longrightarrow}
\bz 
\stackrel{1}{\longrightarrow} \cdots \stackrel{1}{\longrightarrow}
\bz \stackrel{1}{\longrightarrow}\cdots.$$
\end{example}

\begin{example}{\bf (The irrational tree $T_\alpha$)~} 
Let $\alpha = [a_0, a_1, a_2, \dots]$ be the continued fraction
expansion of the positive irrational number $\alpha$.
$Sym_\infty(T_\alpha,v)$ is isomorphic to the two-dimensional integral lattice 
$\bz^2$ with positive cone $Sym_\infty(T_\alpha,v)_+$ corresponding to 
$\{ (x,y)\in\bz^2 ~|~
\alpha x+y\geq 0\}$.
The order unit is $(a_0,1)\in\bz^2$.
Thus, $Sym_\infty(T_\alpha,v)\cong (\bz+\alpha\bz, 
(\bz+\alpha\bz)\cap\br_+, a_0\alpha)$.
The direct sequence  is
$$\bz 
\xrightarrow {\tiny\begin{bmatrix} a_0 \\ 1 \end{bmatrix}}
\bz^2 
\xrightarrow {\tiny\begin{bmatrix}  a_1 &1\\ 1 & 0\end{bmatrix}}
\bz^2
\xrightarrow {\tiny\begin{bmatrix} a_2 & 1\\ 1 & 0\end{bmatrix}}
\bz^2 
\xrightarrow {\tiny\begin{bmatrix} a_3 & 1\\ 1 & 0\end{bmatrix}}
\cdots. 
$$
For this calculation, see Effros and Shen \cite{EfS}.

\end{example}


\section{Scalings and micro-scalings of ultrametrics}
\label{section:Perturbations}

In this section we associate to any compact, ultrametric space $X$ a
Bratteli diagram $B(X)$. This is accomplished by taking any 
rooted, geodesically complete, locally finite simplicial tree $(T,v)$
with $X$ scale equivalent (defined below) to $end(T,v)$ and defining $B(X)=B(T,v)$. Along with scale equivalence we also introduce the notion
of micro-scale equivalence of metric spaces and prove that the equivalence
class of $B(X)$ as a Bratteli diagram only depends on $X$ up to micro-scale
equivalence.

\begin{definition} Let $d$ and $d'$ be metrics on a set $X$.
\begin{enumerate}
\item $d'$  is a {\it scaling} of $d$ if there exists a
homeomorphism $\lambda \co [0,\infty)\to [0,\infty)$ such that $\lambda d=d'$. In this case,
$d$ and $d'$ are said to be {\it scale equivalent}.
\item $d'$  is a {\it micro-scaling} of $d$ if there exist $\epsilon>0$ 
and a
homeomorphism $\lambda \co [0,\infty)\to [0,\infty)$ such that $\lambda d(x,y)=d'(x,y)$
whenever $x,y\in X$ and $\min\{ d(x,y), d'(x,y)\}<\epsilon$. In this case,
$d$ and $d'$ are said to be {\it micro-scale equivalent}.
\end{enumerate}
\end{definition}

Note that scale equivalent or micro-scale equivalent metrics are 
topologically equivalent. In addition, scale equivalence and micro-scale
equivalence are equivalence relations on the set of all metrics on 
$X$. Moreover, if $(X,d)$ is an ultrametric space, then  part of the conditions
in the definition are unnecessary in that whenever $\lambda \co [0,\infty)\to
[0,\infty)$ is a homeomorphism, $\lambda d$ is also an ultrametric.

\begin{definition} Let $h\co X\to Y$ be a bijection between
metric spaces $(X, d_X)$ and $(Y, d_Y)$.
\begin{enumerate}
\item $h\co X\to Y$ is a {\it scale equivalence} if
$d_X$ and $h^*d_Y$ are scale equivalent.
\item $h\co X\to Y$ is a {\it micro-scale equivalence} if
$d_X$ and $h^*d_Y$ are micro-scale equivalent.
\end{enumerate}
\end{definition}

Here $h^*d_Y$ denotes the pull-back metric, $h^*d_Y(x,y) = d_Y(hx,hy)$.
Note that scale equivalences and micro-scale equivalences 
are necessarily homeomorphisms.

\begin{proposition}
\label{prop:uls-mse} If $h\co X\to Y$ is a uniform local similarity
between compact metric spaces $(X, d_X)$ and $(Y, d_Y)$, 
then $h$ is a micro-scale equivalence.
\end{proposition}

\begin{proof}
Since $X$ is compact, there exist $\epsilon>0$ and $\hat\lambda>0$
such that $h|\co B(x,\epsilon)\to B(hx,\hat\lambda\epsilon)$ is
a $\hat\lambda$-similarity for all $x\in X$. Thus, the homeomorphism
$\lambda\co [0,\infty)\to[0,\infty)$, defined by $\lambda(t)=\hat\lambda t$,
shows that $d_X$ and $h^*d_Y$ are micro-scale equivalent 
(because $d_X(x,y) <\epsilon$ implies $d_Y(hx,hy)=\hat\lambda d_X(x,y)=
(\lambda\circ d_X)(x,y)$).
\end{proof}

In particular, local isometries between compact metric spaces are
micro-scale equivalences. 

\begin{proposition}
\label{prop:lr-lr}
If $(X, d_X)$ and $(Y, d_Y)$ are micro-scale equivalent
ultrametric
spaces and $X$ is locally rigid, then $Y$ is also locally rigid.
\end{proposition}

\begin{proof}
Let $h\co X\to Y$ be a bijection, $\lambda\co [0,\infty)\to[0,\infty)$
a homeomorphism, and $\epsilon >0$ such that 
$\lambda d_X(x,y) = d_Y(hx, hy)$ whenever $\min\{ d_X(x,y), d_Y(hx,hy)\} <
\epsilon$ and $x,y \in X$. 

If $0<\delta\leq\min\{\epsilon, \lambda^{-1}(\epsilon)\}$, then
$h|\co B(x,\delta)\to B(hx, \lambda(\delta))$ is a homeomorphism.
If $0<\mu\leq\min\{\epsilon, \lambda(\epsilon),
\lambda^{-1}(\epsilon)\}$, and 
$g\co B(y,\mu)\to B(y,\mu)$ is an isometry between balls in $Y$, then
$$B(h^{-1}y, \lambda^{-1}(\mu))\xrightarrow{~h|~}
B(y,\mu)\xrightarrow{~g~}
B(y,\mu)\xrightarrow{h|^{-1}}
B(h^{-1}y, \lambda^{-1}(\mu))$$
is an isometry between balls in $X$. 

Suppose $y\in Y$ is given and let $x=h^{-1}y$.
Local rigidity of the ultrametric space 
$X$ implies there exists $\epsilon_x>0$ such that
for any $0<\nu\leq\epsilon_x$, every isometry
$B(x,\nu)\to B(x,\nu)$ is the identity 
(see Lemma~\ref{lemma:locally-rigid}).

Let $\epsilon_y=\min\{ \lambda(\epsilon_x), \epsilon, \lambda(\epsilon),
\lambda^{-1}(\epsilon) \}$
and suppose $g\co B(y,\epsilon_y)\to B(y,\epsilon_y)$ is an isometry.
It follows that the composition
$$B(x, \lambda^{-1}(\epsilon_y))\xrightarrow{~h|~}
B(y,\epsilon_y)\xrightarrow{~g~}
B(y,\epsilon_y)\xrightarrow{h|^{-1}}
B(x, \lambda^{-1}(\epsilon_y))$$
is an isometry; hence, it is the identity. Thus, $g$ is the identity
and $Y$ is locally rigid.
\end{proof}

\begin{proposition}
\label{prop:micro-scale}
A micro-scale equivalence $h\co X\to Y$ of metric spaces
induces an isomorphism $h_*\co \ligx\to\ligy$ of topological
groupoids.
\end{proposition}

\begin{proof}
Let  $\lambda\co [0,\infty)\to[0,\infty)$ be
a homeomorphism and $\epsilon >0$ such that 
$\lambda d_X(x,y) = d_Y(hx, hy)$ whenever $\min\{ d_X(x,y), d_Y(hx,hy)\} <
\epsilon$ and $x,y\in X$. 

If $0<\delta\leq\min\{\epsilon, \lambda^{-1}(\epsilon)\}$, then
$h|\co B(x,\delta)\to B(hx, \lambda(\delta))$ is a homeomorphism.
If $0<\mu\leq\min\{\epsilon, \lambda(\epsilon),
\lambda^{-1}(\epsilon)\}$, and 
$g\co B(x,\mu)\to B(y,\mu)$ is an isometry between balls in $X$, then
$$B(hx, \lambda(\mu))\xrightarrow{h|^{-1}}
B(x,\mu)\xrightarrow{~g~}
B(y,\mu)\xrightarrow{~h~}
B(hy, \lambda(\mu))$$
is an isometry between balls in $Y$. 

Thus, there is a function $h_*\co \ligx\to\ligy$ defined by
$h_*[g,x] = [ hgh^{-1}, hx]$ for each local isometry germ 
$[g,x]\in\ligx$ such that the domain of $g$ has
sufficiently small radius, which can be shown to be an isomorphism of
groupoids.
\end{proof}

Note that a micro-scale equivalence need not induce
an isomorphism between local similarity groupoids.

The following result follows immediately from the preceding proposition.
Note that we have already established that, under the hypothesis of
the corollary, that $\ligx$ satisfies the conditions required to apply
Renault's theory of groupoid $C^*$-algebras 
(see Theorem~\ref{theorem:second-main-one-proof}).

\begin{corollary}
\label{cor:inv-mse}
If $(X,d)$ is a compact, locally rigid ultrametric space, then 
the unital $C^*$-algebra $C^*\ligx$  and the unital partially
ordered abelian group $K_0C^*\ligx$ are invariants of $X$
up to micro-scale equivalence of $X$.
\end{corollary}

Note that the preceding proposition and corollary imply
Theorem~\ref{theorem:first-main}(\ref{first-main-five})
of the Introduction.

We now begin to establish 
Theorem~\ref{theorem:first-main}(\ref{first-main-three})
of the Introduction.

\begin{proposition}
\label{prop:scale}
If $(X,d)$ is a compact ultrametric metric space, then 
there exists a rooted, geodesically complete, locally finite simplicial 
tree $(T,v)$ such that $X$ is scale equivalent to $end(T,v)$.
\end{proposition}

\begin{proof} It is well-known that there exists a finite or infinite
sequence $t_0 > t_1 > t_2 >  \cdots >0$ such that 
$\{ d(x,y) ~|~ x,y\in X\} = \{ 0, t_0, t_1, t_2,\dots\}$.
Moreover, the sequence is finite if and only if $X$ is finite, and if
$X$ is infinite, then $\lim_{i\to\infty}t_i=0$. See \cite{Ber}.
Let $\lambda \co [0,\infty)\to[0,\infty)$ be a homeomorphism such that
$\lambda (t_i)= e^{-i}$ for all $i$.

We may assume that $X$ has more than one point, for otherwise the proof is
trivial; hence, $(X,\lambda d)$ has diameter $1$. In \cite{Hug}, there 
is constructed a rooted, geodesically complete $\br$-tree $(T,v)$
such that $(X,\lambda d)$ is isometric to $end(T,v)$. 
It follows that $(X,d)$ is scale equivalent to $end(T,v)$.

From
Proposition~\ref{prop:proper tree} we know that $T$ 
must be a proper $\br$-tree.
It only remains to observe from the construction in \cite{Hug},
that $T$ is in fact a locally finite simplicial tree. This is because
the set of distances in $(X,\lambda d)$ is contained in 
$\{ 0, 1, e^{-1}, e^{-2}, \dots\}$.
See Corollary~\ref{cor:simplicial}.
\end{proof}

The existence of the tree in the preceding proposition allows us to make
the following definition.

\begin{definition}
If $X$ is a compact ultrametric space, then {\it the Bratteli
diagram $B(X)$ associated to $X$} is defined to be
the Bratteli diagram $B(T,v)$ associated to 
a rooted, geodesically complete, locally finite simplicial 
tree $(T,v)$ such that $X$ is scale equivalent to $end(T,v)$.
\end{definition}

The equivalence class of the Bratteli diagram $B(X)$ in the preceding
definition is well-defined as the next result shows.

\begin{theorem}
\label{mse-ebd}
Let $(T,v)$ and $(S,w)$ be rooted, geodesically complete, locally finite 
simplicial trees.
If $end(T,v)$ and $end(S,w)$ are micro-scale equivalent,
then $B(T,v)$ and $B(S,w)$ are equivalent Bratteli diagrams.
\end{theorem}

\begin{proof}
Let $X=end(T,v)$ and $Y=end(S,w)$.
Since $X$ and $Y$ are micro-scale equivalent,
there are homeomorphisms $h\co X\to Y$ and
$\lambda\co [0,\infty)\to[0,\infty)$
and $\epsilon >0$ such that 
$\lambda d_X(x,y) = d_Y(hx, hy)$ whenever $\min\{ d_X(x,y), d_Y(hx,hy)\} <
\epsilon$ and $x,y\in X$. 
We assume $0<\epsilon<1$.
Fix $M>\max\{-\ln\epsilon, -\ln\lambda^{-1}(1)\}$.
Define 
$\hat h\co T\setminus\bar B(v,M)\to S\setminus\bar B(w,-\ln\lambda(e^{-M}))$
as follows. 
Let $x\in X$ and $M<t<\infty$; thus,
$x(t)\in T\setminus\bar B(v,M)$.
Set $\hat h(x(t)) = (h(x))(-\ln\lambda(e^{-t}))$.

\begin{claim} $\hat h$ is a homeomorphism.
\end{claim}

\begin{proof*}\emph{of Claim.}
We begin by showing that $\hat h$ is well-defined.
For $x,y\in X$ and $M<t<\infty$ such that $x(t)=y(t)$, we must show
that $(h(x))(-\ln\lambda(e^{-t})) =
(h(y))(-\ln\lambda(e^{-t}))$.
Assuming $x\neq y$,
$d_X(x,y) = e^{-t_0}$, where
$t_0=\sup\{ s\geq 0 ~|~ x(s)=y(s)\}$.
Thus,
$t_0\geq t\geq-\ln\epsilon$, which implies
$d_X(x,y)\leq e^{-t}<\epsilon$. Hence,
$\lambda(e^{-t_0})=\lambda d_X(x,y) = d_Y(hx, hy) = e^{-t_1}$,
where
$t_1=\sup\{ s\geq 0 ~|~ (hx)(s) = (hy)(s)\}$.
It follows that $t_1=-\ln\lambda(e^{-t_0})$.
Since $t\leq t_0$, we have
$-\ln\lambda(e^{-t})\leq -\ln\lambda(e^{-t_0})=t_1$.
The definition of $t_1$ and this inequality imply
that 
$(h(x))(-\ln\lambda(e^{-t}))=(h(y))(-\ln\lambda(e^{-t}))$.
Therefore, $\hat h$ is well-defined.

To see that $\hat h$ is bijective, define
$g\co S\setminus\bar B(w,-\ln\lambda(e^{-M}))\to T\setminus\bar B(v,M)$
by $g(y(s)) = (h^{-1}(y))(-\ln\lambda^{-1}(e^{-s}))$
for $y\in Y$ and $-\ln\lambda(e^{-M})<s<\infty$.
It can be checked that $g$ is well-defined and
$g=(\hat h)^{-1}$.

We now proceed to show that $\hat h$ is continuous.
Suppose first that $x\in X$ and $M<s<t<\infty$, so that in $T$,
$d(x(s), x(t)) = t-s$.
In $S$, $d(\hat h(x(s)), \hat(x(t)))=
-\ln\lambda(e^{-t}) +\ln\lambda(e^{-s})=
\ln\left({\frac {\lambda(e^{-s})} {\lambda(e^{-t})}}\right)$.
The continuity of $\hat h$ on $x((M,\infty))$ follows from this.

Now suppose $x,y\in X$, $x\neq y$ and let
$d_X(x,y) = e^{-t_0}$.
Further suppose $M<t_0<s\leq t<\infty$.
In $T$, $d(x(s),y(t))= s+t-2t_0$.
In $S$, $d(\hat h(x(s)),\hat h(y(t))
= -\ln\lambda(e^{-s})-\ln\lambda(e^{-t})+2\ln\lambda(e^{-t_0})
= \ln\left(\frac{[\lambda(e^{-t_0})]^2} {\lambda(e^{-s})\lambda(e^{-t})}\right)$.
If $d(x(s), y(t))$ is small, then $s$ and $t$ are both close to
$t_0$; hence, $d(\hat h(x(s)), \hat h(y(t))$ is small.

This now establishes that $\hat h$ is continuous on connected components
of $T\setminus\bar B(v,M)$; thus, $\hat h$ is continuous.
Likewise, $g$ is continuous and $\hat h$ is a homeomorphism.
This completes the proof of the claim.
\end{proof*}

Suppose $n>M$ is an integer and $x,y\in X$. Let $T_x$ and $T_y$ be the
rooted subtrees of $T$ descending from $x(n)$ and $y(n)$, respectively.
Let $S_x$ and  $S_y$ be the rooted subtrees of $S$ descending from
$\hat h(x(n))$ and $\hat h(y(n))$, respectively.

\begin{claim} 
\label{claim:isometric}
$(T_x, x(n))$ and $(T_y, y(n))$ are rooted isometric
if and only if
$(S_x, \hat h(x(n)))$ and $(S_y, \hat h(y(n)))$ are rooted isometric.
\end{claim}

\begin{proof*}\emph{of Claim.}
We suppress the roots of the subtrees from the notation.
Suppose $T_x$ and $T_y$ are rooted isometric.
Then $end(T_x)$ and $end(T_y)$ are isometric when these 
end spaces  are given the end space metric 
as recalled in Section~\ref{sec:ends} 
(see Proposition~\ref{proposition:functor}).
However, we want to give $end(T_x)$ and $end(T_y)$
the metrics they inherit as subspaces of $end(T,v)$---that is,
under the identifications
$end(T_x)=\bar B(x, e^{-n})\subseteq end(T,v)$, and
$end(T_y)=\bar B(y, e^{-n})\subseteq end(T,v)$.
Since the pairs of possible metrics differ by a factor of $e^{-n}$,
$end(T_x)$ and $end(T_y)$ remain isometric with the subspace metrics.
Let $j\co \bar B(x, e^{-n})\to\bar B(y, e^{-n})$ be an isometry.
Then $\hat h j (\hat h|)^{-1}\co \bar B(h(x),\lambda(e^{-n}))\to
\bar B(h(y),\lambda(e^{-n}))$
is also an isometry.
This means that $end(S_x)$ and $end(S_y)$ are isometric
as subspaces of $end(S,w)$. As above, we conclude that 
$S_x$ and $S_y$ are rooted isometric.
Similar reasoning gives the converse.
This completes the proof of the claim.
\end{proof*}

We can now complete the proof that $B(T,v)$ and $B(S,w)$ are equivalent.
We will show that $\symtv\cong\symsw$.
Since these groups are isomorphic to the unital dimension groups of $B(T,v)$
and $B(S,w)$, respectively, (by Proposition~\ref{prop:dim gp}),
it follows from Bratteli's Theorem~\ref{theorem:Bratteli, Elliott}
that $B(T,v)$ and $B(S,w)$ are isomorphic.

Let $D_X= \{ t\in\br ~|~ \hbox{ there exists } x,y\in X 
\hbox{ such that } d_X(x,y)=t\}$
and  $D_Y= \{ t\in\br ~|~ \hbox{ there exists } x,y\in Y
\hbox{ such that } d_Y(x,y)=t\}$,
the distance sets of $X$ and $Y$, respectively. 
Write $D_X = \{ 0 <\cdots <\nu_{i+1}<\nu_i<\cdots<\nu_0\leq 1\}$.
Then $D_Y = \{ 0 <\cdots <\lambda(\nu_{i+1})<\lambda(\nu_i)<\cdots<
\lambda(\nu_0)\leq 1\}$.
For each $i = 0, 1, 2,\dots$, let $L_i=-\ln\nu_i$ 
and $M_i=-\ln\lambda(\nu_i)$.
Note that 
$0\leq L_0 < L_1 < L_2 < \cdots$ and 
$0\leq M_0 < M_1 < M_2 < \cdots$.
The $L_i$'s and $M_i$'s correspond to the levels in the trees $T$ and $S$, 
respectively, where nontrivial branching occurs.

For each $i=0,1,2,\dots$, let
$\mu_i$ be the number of rooted isometry classes of level $L_i$ 
subtrees of 
$(T,v)$ and let
$\tau_1^i, \tau_2^i,\dots \tau_{\mu_i}^i$ be a complete
set of representatives of the rooted isometry classes of level $L_i$
subtrees.

According to Claim\ref{claim:isometric},
$\mu_i$ is also the number of rooted isometry classes of level $M_i$
subtrees of $(S,w)$. Moreover,
$\hat h(\tau_1^i), \hat h(\tau_2^i),\dots, \hat h(\tau_{\mu_i}^i)$ is a complete
set of representatives of the rooted isometry classes of level $M_i$
subtrees of $(S,w)$.

For $0\leq i$, $1\leq k\leq \mu_{i+1}$ and $1\leq \ell \leq \mu_i$, 
let $\alpha_{k\ell}^i$ be the number of level $(L_{i+1}-L_i)$
subtrees of $\tau_\ell^i$ that are rooted isometric to
$\tau_k^{i+1}$. The matrix $\alpha_i=[\alpha_{k\ell}^i]$ is a
$(\mu_{i+1}\times \mu_i)-$ matrix with nonnegative
integral entries.

Using Claim~\ref{claim:isometric} again, it follows that 
$a_{k\ell}^i$ is also the number of level $(M_{i+1}-M_i)$
subtrees of $\hat h(\tau_\ell^i)$ that are rooted isometric to
$\hat h(\tau_k^{i+1})$.

We claim that there is an isomorphism of unital partially ordered abelian
groups:
$$\symtv \cong\lim_\to (\bz \stackrel{\alpha_0}{\longrightarrow}
\bz^{\mu_1} \stackrel{\alpha_1}{\longrightarrow} \bz^{\mu_2} 
\stackrel{\alpha_2}{\longrightarrow}\bz^{\mu_3} 
\stackrel{\alpha_3}{\longrightarrow} \cdots
\bz^{\mu_i} \stackrel{\alpha_i}{\longrightarrow}\cdots).$$

It will follow by a similar argument that $\symsw$ is also isomorphic
to this direct limit, finishing the proof.

We use the notation of Section~\ref{section:TheSymmetryAtInfinityGroup}.
In particular, we have the matrices $A_i= [a_{k\ell}^i]$
for $i=0, 1,2, \dots$, $1\leq \ell \leq m_i$, and $1leq k\leq m_{i+1}$.
Note that $\mu_i= m_{L_i}$ and we can take $\tau_\ell^i = T_\ell^{L_i}$.
In order to show that $\symtv\cong\lim_{\to}\alpha_i$,
it suffices to show that $\alpha_i = A_{(L_{i+1}-1)} \cdots 
A_{(L_i+1)}A_{L_i}$ for $i, j = 0,1,2, \dots$.
Hence, the following claim completes the proof.

\begin{claim}
For $1\leq\ell\leq m_i$, $1\leq k\leq m_{i+j+1}$, and
$i, j = 0,1,2, \dots$, the $k\ell$-entry of the product
$A_{i+j}\cdots A_i$ is the number of level $j$
subtrees of $T_\ell^i$ that are rooted isometric to $T_k^{i+j+1}$.
\end{claim}

\begin{proof*}\emph{of Claim.}
The proof is by induction on $j$. The statement is obviously 
true for $j=0$; so assume $j>0$ and the statement is true for
$j-1$. Let $B= A_{i+j-1}\cdots A_i$ and denote its entries
by $B=[b_{pq}]$.
By the inductive assumption, $b_{p\ell}$ is the number of level
$j-1$ subtrees of $T_\ell^i$ that are rooted isometric to 
$T_p^{i+j}$.
The entries of the matrix $A_{i+j}= [a_{kp}^{i+j}]$ have the following
interpretation by definition:
$a_{kp}^{i+j}$ is the number of level $1$ subtrees of $T_p^{i+j}$
that are rooted isometric to $T_k^{i+j+1}$.
Hence, the number of level $j$ subtrees of $T_\ell^i$ that are rooted
isometric to $T_k^{i+j+1}$ is given by
$\sum_{p=1}^{m_{i+j}} a_{kp}b_{p\ell}$; that is, the $k\ell$-entry of
$A_{i+j}\cdots A_i$.
\end{proof*}

This completes the proof of the theorem.
\end{proof}

The converse of the preceding theorem is not true, as the following 
example shows.

\begin{figure}[htpb]

\begin{picture}(370,100)

\put(85,0){$T$}
\put(285,0){$S$}

\multiput(0,20)(30,0){6}{\vector(-1,-2){0}}
\multiput(20,20)(30,0){6}{\vector(1,-2){0}}

\multiput(200,20)(30,0){6}{\vector(-1,-2){0}}
\multiput(220,20)(30,0){6}{\vector(1,-2){0}}

\multiput(210,40)(30,0){6}{\vector(0,-1){20}}

\multiput(0,20)(30,0){6}{\line(1,2){10}}
\multiput(20,20)(30,0){6}{\line(-1,2){10}}

\multiput(200,20)(30,0){6}{\line(1,2){10}}
\multiput(220,20)(30,0){6}{\line(-1,2){10}}

\multiput(10,40)(90,0){2}{\line(3,2){30}}
\multiput(40,60)(90,0){2}{\line(3,-2){30}}
\multiput(40,60)(90,0){2}{\line(0,-1){20}}
\put(40,60){\line(2,1){45}}
\put(130,60){\line(-2,1){45}}

\multiput(210,40)(60,0){3}{\line(3,4){15}}
\multiput(240,40)(60,0){3}{\line(-3,4){15}}
\put(225,60){\line(3,1){60}}
\put(345,60){\line(-3,1){60}}
\put(285,60){\line(0,1){20}}

\multiput(10,40)(30,0){6}{\circle*{5}}
\multiput(210,40)(30,0){6}{\circle*{5}}

\multiput(40,60)(90,0){2}{\circle*{5}}
\multiput(225,60)(60,0){3}{\circle*{5}}

\put(85,82.5){\circle*{5}}
\put(285,80){\circle*{5}}

\end{picture}
\caption{The trees of Example~\ref{example:2-3-trees}
with $B(T,v)$ and $B(S,w)$ equivalent.}
\label{fig:2-3-trees}
\end{figure}
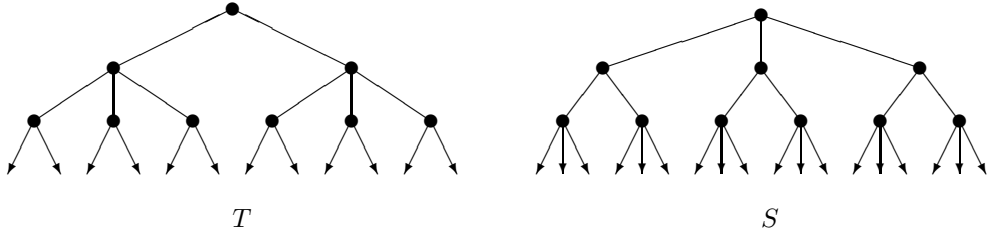

\begin{example}
\label{example:2-3-trees}
There are two rooted, geodesically complete,
locally finite simplicial trees, $(T,v)$ and $(S,w)$, such that 
$B(T,v)$ and $B(S,w)$ are equivalent Bratteli diagrams,
but $X=end(T,v)$ and $Y=end(S,w)$ are not 
micro-scale equivalent.
The trees are pictured in 
Figure~\ref{fig:2-3-trees}.
Elements $x\in X$ are sequences $ x= (x_0, x_1, x_2,\dots)$
such that
$x_i\in\left\{
\begin{array}{ll}
\{ 0,1\} & \text{if $i$ is even,}\\
\{ 0, 1, 2\} & \text{if $i$ is odd.}\\
\end{array}\right.$
Elements $y\in Y$ are sequences $ y= (y_0, y_1, y_2,\dots)$
such that
$y_i\in\left\{
\begin{array}{ll}
\{ 0,1\} & \text{if $i$ is odd,}\\
\{ 0, 1, 2\} & \text{if $i$ is even.}\\
\end{array}\right.$
Suppose $X$ and $Y$ are micro-scale equivalent.
Then there are homeomorphisms $h\co X\to Y$ and
$\lambda\co [0,\infty)\to[0,\infty)$
and $0<\epsilon\leq 1$ such that 
$\lambda d_X(x,y) = d_Y(hx, hy)$ whenever $\min\{ d_X(x,y), d_Y(hx,hy)\} <
\epsilon$ and $x,y\in X$. 
There exists integers $i_0>-\ln\epsilon$ and $c\leq i_0$ such that
$\lambda(e^{-i})= e^{c-i}$ for all $i\geq i_0$.
For each $i=0,1,2,\dots$, let
$\alpha_i =\left\{\begin{array}{ll}
2\cdot(3\cdot 2)^{\frac i 2} &\text{if $i$ is even,}\\
(2\cdot 3)^{\frac {i+1} 2} &\text{ if $i$ is odd}\\
\end{array}\right.$
and 
$\beta_i =\left\{\begin{array}{ll}
3\cdot(2\cdot 3)^{\frac i 2} &\text{if $i$ is even,}\\
(3\cdot 2)^{\frac {i+1} 2} &\text{ if $i$ is odd.}\\
\end{array}\right.$
The number $\left\{\begin{array}{l}\alpha_i\\ \beta_i\end{array}\right.$
is the maximum number of distinct points of
$\left\{\begin{array}{l}X\\ Y\end{array}\right.$
whose distances from each other are $e^{-i}$.
Clearly, $\alpha_i=\beta_{i-c}$ for all $i\geq i_0$.
In particular, $\alpha_{i_0}=\beta_{i_0-c}$
and
$\alpha_{i_0+1}=\beta_{i_0+1-c}$.
This implies $c=0$ and both $i_0$ and $i_0+1$ are odd---a contradiction;
hence, $X$ and $Y$ are not micro-scale equivalent.
The Bratteli diagrams $B(T,v)$ and $B(S,w)$ both
telescope to the Bratteli diagram
$D=(\mathcal{V,E})$ with $\mathcal{V}_i=\{ v_i\}$, a single vertex and
six edges from $v_i$ to $v_{i+1}$ for each $i=0,1,2,\dots$.
\end{example}

The following result is a restatement of
Theorem~\ref{theorem:first-main}(\ref{first-main-three})
in the Introduction.

\begin{theorem}
If $X$ is a compact, locally rigid ultrametric space, then there exists
a Bratteli diagram $B(X)$ such that $\ligx$ is isomorphic to the
path groupoid of $B(X)$.
\end{theorem}

\begin{proof} 
Proposition~\ref{prop:scale} gives a 
rooted, geodesically complete, locally finite simplicial 
tree $(T,v)$ such that $X$ is scale equivalent to $Y:=end(T,v)$.
We have defined $B(X):= B(T,v)$; hence,
we have equality of the path groupoids
$\pg(B(X)) \cong \pgb$.
Now Corollary~\ref{cor:summary} implies there is an isomorphism
of topological groupoids
$\pgb\cong\ligy$. Since $X$ and $Y$ are scale equivalent,
Proposition~\ref{prop:micro-scale} implies there is an isomorphism
of topological groupoids $\ligx\cong\ligy$.
Thus, $\ligx\cong\pg(B(X))$ as required.
\end{proof}

We now reinterpret our results on invariants for ultrametric spaces
as invariants for trees.

\begin{corollary}
Let $(T,v)$ and $(S,w)$ be rooted, geodesically complete,
locally finite simplicial trees. If $(T,v)$ and $(S,w)$ are uniformly
isometric at infinity, then
\begin{enumerate}
\item $B(T,v)$ and $B(S,w)$ are equivalent Bratteli diagrams, and 
\item $\symtv$ and $\symsw$ are isomorphic partially ordered
abelian groups.
\end{enumerate}
If, in addition, either $X:=end(T,v)$ or $Y:=end(S,w)$ is locally rigid, then so is
the other, and
\begin{enumerate}
\item[3.] 
$(K_0C^*\ligx, K_0C^*\ligx_+, [1])$ and 
$(K_0C^*\ligy, K_0C^*\ligy_+, [1])$ are isomorphic unital partially
ordered abelian groups.
\end{enumerate}
\end{corollary}

\begin{proof}
The first statement follows from Propositions~\ref{proposition:functor}
and \ref{prop:uls-mse} and Theorem~\ref{mse-ebd}.
The second statement follows from the first, Proposition~\ref{prop:dim gp},
and Theorem~\ref{theorem:Bratteli, Elliott}.
The final statement follows from the first,
Proposition~\ref{prop:lr-lr}, and Corollary~\ref{cor:summary}.
\end{proof}

\begin{corollary}
\label{unif-iso-inf-inv}
Let $(T,v)$ and $(S,w)$ be rooted, geodesically complete,
proper $\br$-trees that are uniformly
isometric at infinity.
If either $X:=end(T,v)$ or $Y:=end(S,w)$ is locally rigid, then so is
the other, and
$$\text{$(K_0C^*\ligx, K_0C^*\ligx_+, [1])$ and
$(K_0C^*\ligy, K_0C^*\ligy_+, [1])$}$$ are isomorphic unital partially
ordered abelian groups.
\end{corollary}

\begin{proof}
This follows from Propositions~\ref{proposition:functor}, \ref{prop:uls-mse},
\ref{prop:lr-lr}, and \ref{prop:scale} and Corollary~\ref{cor:inv-mse}.
\end{proof}


\section{Faithful unitary representations}
\label{section:FUR}

The goal of this section is to prove the following theorem, which is the third part of
Theorem~\ref{theorem:second-main}.

\begin{theorem}
\label{theorem:FUR}
If $X$ is a compact ultrametric space with a countable subgroup $\Gamma\leq LS(X)$ acting locally rigidly
on $X$, then
there is a faithful unitary representation of $\Gamma$ into  $C^*\gammax$.
\end{theorem}

It follows from Theorem~\ref{theorem:second-main}(i)
(Corollary~\ref{cor:ls-grpd-props})
that under the hypothesis of Theorem~\ref{theorem:FUR}
Renault's theory \cite{Ren} can be applied so that
$C^*\gammax$ is defined. 

Of course, by a faithful unitary representation of $\Gamma$ into $C^*\gammax$, we mean an injective homomorphism
of $\Gamma$ into the multiplicative group of unitary elements of $C^*\gammax$. This is proved by establishing, 
in Corollary~\ref{corollary:rep} below, an injective homomorphism $\rho\co \Gamma\to C_c(\gammax)$ into the unitary group
of the convolution algebra of $\gammax$. Since the $C^*$-algebra of $\gammax$ is a completion of the convolution
algebra of $\gammax$, Theorem~\ref{theorem:FUR} follows.

We begin by fixing notation for the convolution $*$-algebras of groups and groupoids.
For more details, see
Muhly \cite{Muh}, Paterson \cite{Pat} and Renault \cite{Ren}.

If $\Gamma$ is a discrete group, then $C_c(\Gamma)$ denotes the convolution $*$-algebra of $\Gamma$, otherwise known
as the complex algebra $\bc\Gamma$:
$$C_c(\Gamma) := \{ f\co \Gamma\to\bc ~|~ \text{$f$ has finite support}\}.$$
Multiplication and involution on this complex vector space are given by
$$(f*g)(\gamma) :=\sum_{\beta\in\Gamma}f(\beta)g(\beta^{-1}\gamma) \text{ and
$f^*(\gamma) :=\overline{f(\gamma^{-1})}$},$$
where $\overline{~\cdot~}$ denotes complex conjugation.

Now if $\g$ is a locally compact, Hausdorff \'etale groupoid, then
$$C_c(\g) := \{ f\co \g\to\bc ~|~ \text{$f$ is continuous and has compact support}\}.$$
For each $u$ in the unit space of $\g$, $r^{-1}(u):=\g^u$ is discrete. 
Thus, each element of $C_c(\g)$ restricts to a function on each $\g^u$ with finite support. Therefore, multiplication
and involution on the complex vector space $C_c(\g)$ may be defined by
$$(f*g)(y) :=\sum_{x\in\g^{r(y)}}f(x)g(x^{-1}y) \text{ and
$f^*(x) :=\overline{f(x^{-1})}$}.$$

Thus, $C_c(\g)$ is a topological $*$-algebra.

If the unit space $X=\{ xx^{-1} ~|~ x\in\g\} =\{ r(x) ~|~ x\in\g\}$ is compact, then the algebra 
$C_c(\g)$ has a unit $1$, namely, the characteristic function $\chi_X$.\footnote{Do not confuse $1$ with the function
on $\g$ that is identically $1$, i.e., $\chi_\g$. Since $\g$ need not be compact,
$\chi_\g$ need not be in $C_c(\g)$. On the other hand, $0\in C_c(\g)$ is the function that is identically $0$.}
In that case, the {\it unitary group of $C_c(\g)$} is the multiplicative group
$$\{ f\in C_c(\g)~|~ f^*f=1=ff^*\}.$$

{\it For the remainder of this section, let $X$ be a (nonempty) compact ultrametric space with a subgroup
$\Gamma\leq LS(X)$ acting locally rigidly on $X$.} Even though $\Gamma$ need not be a discrete subgroup of $LS(X)$, 
we will endow $\Gamma$ with the discrete topology. Denote the identity of $\Gamma$ by $e$; that is, $e=\id_X$.

For each $\gamma\in\Gamma$, let
$$A_\gamma :=\{ [\gamma,x] ~|~ x\in X\}\subseteq\gammax.$$

\begin{lemma}For each $\gamma\in\Gamma$, 
$A_\gamma$ is compact and open in $\gammax$.
\end{lemma}

\begin{proof}
Clearly,
$A_\gamma= \cup_{x\in X} U(\gamma,x,1)$. This shows $A_\gamma$ is open.
To see that $A_\gamma$ is compact, note that $E\co X\to A_\gamma$, defined by 
$E(x)=[\gamma,x]$, is continuous. This is because $E^{-1}U(\gamma,x,\epsilon)=B(x,\epsilon)$ for each
$x\in X$ and $\epsilon >0$.
\end{proof}

Thus, for each $\gamma\in\Gamma$, the characteristic function $\chi_{A_\gamma}$ of $A_\gamma$ is in $C_c(\gammax)$.

\begin{proposition}
\label{proposition:characteristic}
If $\gamma, \gamma_1, \gamma_2\in\Gamma$, then
\begin{enumerate}
\item \label{item:adjoint}
$\chi_{A_\gamma}^*= \chi_{A_{\gamma^{-1}}}$
\item \label{item:product}
$\chi_{A_{{\gamma}_1{\gamma}_2}}=\chi_{A_{\gamma_1}}*\chi_{A_{\gamma_2}}$
\item \label{item:unit-space} 
$A_e=X$, the unit space
\item \label{item:unit}
$\chi_{A_\gamma}=1$ if and only if $\gamma=e$
\item $A_\gamma\neq\emptyset$ (so that $\chi_{A_\gamma}\neq 0$).
\end{enumerate}
\end{proposition}

\begin{proof}
(1) First note that for $x\in X$ and $\beta, \gamma\in\Gamma$,
$[\beta,x]=[\gamma^{-1},x]$ if and only if $[\beta^{-1},\beta x]=[\gamma, \beta x]$.
Thus,
$$\chi_{A_\gamma}^*([\beta, x])= 
\overline{\chi_{A_\gamma}([\beta, x]^{-1})}= 
\chi_{A_\gamma}([\beta, x]^{-1})= 
\chi_{A_\gamma}([\beta^{-1}, \beta x])$$
$$= \left\{ \begin{array}{ll}
1 & \text{if $[\beta^{-1},\beta x]=[\gamma,\beta x]$}\\
0 &\text{otherwise} \end{array}\right.
= \left\{ \begin{array}{ll}
1 & \text{if $[\beta,x]=[\gamma^{-1}, x]$}\\
0 &\text{otherwise} \end{array}\right.$$
$$=\chi_{A_{\gamma^{-1}}}([\beta, x]).$$ 

(2) Let $[\gamma,x]\in\gammax$ be given. 
Then 
$$\chi_{A_{{\gamma}_1{\gamma}_2}}[\gamma,x]
= \left\{ \begin{array}{ll}
1 & \text{if $[\gamma,x]=[\gamma_1\gamma_2, x]$}\\
0 &\text{otherwise.} \end{array}\right.
$$
On the other hand,
$$\chi_{A_{\gamma_1}}*\chi_{A_{\gamma_2}}[\gamma,x]=
$$
$$
\sum_{[\beta,y]\in r^{-1}(r[\gamma,x])}
\chi_{A_{\gamma_1}}[\beta,y]\cdot\chi_{A_{\gamma_2}}([\beta,y]^{-1}[\gamma,x])
=
\sum_{\beta y=\gamma x}
\chi_{A_{\gamma_1}}[\beta,y]\cdot\chi_{A_{\gamma_2}}[\beta^{-1}\gamma,x].$$
Since
$$\chi_{A_{{\gamma}_1}}[\beta,y]
= \left\{ \begin{array}{ll}
1 & \text{if $[\beta,y]=[\gamma_1, y]$}\\
0 &\text{otherwise,} \end{array}\right.
$$
we have 
$$\chi_{A_{\gamma_1}}*\chi_{A_{\gamma_2}}[\gamma,x]=
\sum_{S}
\chi_{A_{\gamma_2}}[\beta^{-1}\gamma,x],
$$
where
$S=\{ [\beta,y]\in\gammax ~|~ \text{$\beta y=\gamma x$ and $[\beta, y]=[\gamma_1,y]$}\}$.
If $[\beta, y]\in S$, then $\gamma x=\beta y=\gamma_1 y$; thus, 
$y=\gamma_1^{-1}\gamma x$ and $[\beta, y] =[\gamma_1, \gamma_1^{-1}\gamma x]$.
It follows that $S=\{ [\gamma_1, \gamma_1^{-1}\gamma x]\}$; in other words,
the sum has only one term and
$$\chi_{A_{\gamma_1}}*\chi_{A_{\gamma_2}}[\gamma,x]=
\chi_{A_{\gamma_2}}[\gamma_1^{-1}\gamma,x]
=
\left\{ \begin{array}{ll}
1 & \text{if $[\gamma_1^{-1}\gamma, x]=[\gamma_2, x]$}\\
0 &\text{otherwise} \end{array}\right.
$$
$$=
\left\{ \begin{array}{ll}
1 & \text{if $[\gamma, x]=[\gamma_1\gamma_2, x]$}\\
0 &\text{otherwise} \end{array}\right.
=
\chi_{A_{{\gamma}_1{\gamma}_2}}[\gamma,x].
$$

(3) This follows from the description of the unit space in Remark~\ref{remark:unit-space}.

(4) If $\gamma=e$, then $A_e$ is the unit space by \ref{item:unit-space}.
Thus, $\chi_{A_\gamma}$ is the unit of $C_c(\gammax)$ by the general remarks made above.
Conversely, if $\chi_{A_\gamma}=1$, then $\chi_{A_\gamma}=\chi_{A_e}$ and $A_\gamma= A_e$.
Thus, for each $x\in X$, $[\gamma,x]=[\id_X,x]$; hence, $\gamma=e$.

(5) is obvious.
\end{proof}

\begin{corollary}
\label{corollary:rep}
$\rho\co \Gamma\to C_c(\gammax)$, defined by $\rho(\gamma)=\chi_{A_\gamma}$,
is an injective homomorphism into the unitary group of $C_c(\gammax)$.
\end{corollary}

\begin{proof}
That $\rho$ is a homomorphism follows from \ref{proposition:characteristic}~(\ref{item:product}).
The image of $\rho$ lies in the unitary group by 
\ref{proposition:characteristic}~(\ref{item:adjoint}),~(\ref{item:product}),~(\ref{item:unit}).
The injectivity of $\rho$ follows from \ref{proposition:characteristic}~(\ref{item:unit}).
\end{proof}

This completes the proof of Theorem~\ref{theorem:FUR}.

\begin{example}
Let $X=\{x_1,\dots, x_n\}$ be the finite ultrametric space
with $d(x_i,x_j) =1 $ if $i\neq j$.
Then $\Gamma:= LS(X)= Isom(X)$ is the symmetric group $S_n$ on $n$ elements and $\Gamma$ acts locally rigidly on $X$. The groupoid $\gammax$ is the transitive principle groupoid on $n$ elements
(that is, the trivial groupoid $X\times X$) and $C_c(\gammax)=\bm_n(\bc)$.
The homomorphism
$\rho\co \Gamma\to \bm_n(\bc)$ defined above is the representation of 
$S_n$ by permutation matrices. 
\end{example}

\section{Miscellanea}
\subsection{Isometries of trees}

In this section we show that a group acting by isometries on a tree 
sometimes leads to a group of local similarities
acting locally rigidly on the end space of the tree. For the general theory of isometries on trees, see
for example Alperin and Bass \cite{AlB}, Bestvina \cite{Bes}, Chiswell \cite{Chi}, Morgan and Shalen \cite{MoS},
and Serre \cite{Ser}.

Let $(T,v)$ be a geodesically complete, rooted, locally finite, simplicial tree, let $X=end(T,v)$, and
let $Isom(T)$ denote the group of isometries on $T$.\footnote{Some of the facts and constructions 
in this section hold in the more general context of $\br$-trees.} In particular, $X$ is compact ultrametric.

\paragraph{Description of a group homomorphism $\epsilon\co Isom(T)\to LS(X)$.}

We will use the notation and terminology from \cite{Hug}.
Let $\gamma\co T\to T$ be an isometry and let $r=1+d(v,\gamma v)$.
Then $\partial B(v,r)$ and $\partial B(\gamma v,r)$ are cut sets for $(T,v)$
(cf. \cite[Example 3.2]{Hug}). Let $x\in X$; thus, $x\co [0,\infty)\to T$ is
an isometric embedding with $x(0)=v$.
Let $\hat x\co [0,|| \gamma(x(r))||]\to X$ be the unique isometric embedding such that
$\hat x(0)=v$ and $\hat x(||\gamma(x(r)||) = \gamma(x(r))$ (here $||y|| := d(v,y)$ for all $y\in T$).

Define $\gamma_*(x)\co [0,\infty)\to T$ by
$$\gamma_*(x)(t) = 
\left\{ \begin{array}{ll}
\hat x(t) & \text{if $0\leq t\leq ||\gamma(x(r))||$}\\
\gamma\circ x(t-||\gamma(x(r))||+r) &\text{if $||\gamma(x(r))||\leq t$.} \end{array}\right.
$$
Then $\gamma_*(x)\in X$, $\gamma_*\co X\to X$ is in $LS(X)$, and 
$\epsilon\co Isom(T)\to LS(X)$ defined by $\epsilon(\gamma)= \gamma_*$ is a 
group homomorphism (cf. \cite[Section 5]{Hug}).\footnote{Of course, it is quite well-known that
$\epsilon$ is a homomorphism from $Isom(T)$ into the group of homeomorphisms of $X$; we are
just pointing out here that when the end space $X$ is given the natural metric described herein, that the image
of $\epsilon$ lies in $LS(X)$.}

It follows that for $x\in X$ and $\gamma\in Isom(T)$, $\epsilon(\gamma)(x)=x$ if and only if there exists
$t_1,t_2\geq 0$ such that
$x([t_1,\infty))=\gamma x([t_2,\infty))$. From this the following {\it key property of $\epsilon$} can be verified:
if $x\in X$, $t_0\geq 0$, and $\gamma\in Isom(T)$ such that $\gamma\in\Gamma_{x(t_0)}$ and $\epsilon(\gamma)\in
\epsilon(\Gamma)_x$, then $\gamma\in\Gamma_{x(t)}$ for all $t\geq t_0$.

\begin{remark}
\label{cd:isom}
There exists a commuting diagram of groups and group homomorphisms:
$$\begin{CD}
Isom(T,v) @>{\epsilon|}>{\cong}> Isom(X)\\
@V{}VV @VV{}V\\
Isom(T) @>{\epsilon}>> LS(X)
\end{CD}$$
The vertical arrows are inclusions of subgroups. The top horizontal map $\epsilon|$ is
easily seen to be an isomorphism (cf. Proposition~\ref{proposition:functor}
and \cite[Corollary 8.7]{Hug}).
\end{remark}

\begin{example} $\epsilon\co Isom(T)\to LS(X)$ need not be surjective.
There are many local similarities of the end space of the Sturmian tree $T$
(see Example~\ref{ex:Sturmian Tree}), but $T$ has no non-trivial isometries.
In general, $\epsilon$ is rarely surjective.
\end{example}

\begin{example} 
Let $X$ be the space in Example~\ref{example:finite-not-lr} and let $(T,v)$ be the geodesically
complete, rooted tree with $X= end(T,v)$. The subgroup $\Gamma\cong\bz/2$ of $LS(X)$ defined in 
\ref{example:finite-not-lr}
is the isomorphic image of a subgroup $\widehat\Gamma$ of $Isom(T,v)$ under $\epsilon$.
Note that ${\widehat\Gamma}_x$ is finite for each $x\in T$ even though $\Gamma$ does not
act locally rigidly on $X$.
\end{example}

\begin{theorem} 
\label{theorem:tree-isom}
Let $\Gamma$ be a subgroup of $Isom(T)$. The group $\epsilon(\Gamma)$ acts locally rigidly on
$X$ if and only if 
for every $x\in X$ and for every $\gamma\in\Gamma$ such that there exists $t_0\geq 0$ with
$\gamma\in\Gamma_{x(t)}$ for all $t\geq t_0$, there exists $t_1\geq t_0$ so that if 
$y\in X$ and $y(t_1)=x(t_1)$, then $\gamma\in\Gamma_{y(t)}$ for all $t\geq t_1$.
\end{theorem}

\begin{proof}
Assume first that $\epsilon(\Gamma)$ acts locally rigidly. Let $x\in X$, $t_0\geq 0$,
and $\gamma\in\cap\{\Gamma_{x(t)} ~|~ t\geq t_0\}$ be given.
Clearly, $\epsilon(\gamma)\in\epsilon(\Gamma)_x$ and $sim(\epsilon(\gamma),x)=1$. 
By the definition of a locally rigid action, there exists $\delta >0$ such that $\epsilon(\gamma)\in\epsilon(\Gamma)_y$
for all $y\in B(x,\delta)$. Choose $t_1\geq t_0$ such that if $y\in X$ and $y(t_1)=x(t_1)$, then $y\in B(x,\delta)$.
For such a $y$, $\gamma\in\Gamma_{y(t_1)}$ and $\epsilon(\gamma)\in \epsilon(\Gamma)_y$.
Hence, it follows from the key property of $\epsilon$ mentioned above that $\gamma\in\Gamma_{y(t)}$ for all $t\geq t_1$.

For the converse, first observe that there is a bijection from $X$ to the set of {\it open ends of $T$ in the sense of}
\cite{AlB} given by $x\mapsto [x([0,\infty))]$, the open end determined by
the image of $x$. Moreover, for $\gamma\in Isom(T)$ and $x\in X$, $\epsilon(\gamma)\in\epsilon(\Gamma)_x$ if and only if
$\gamma$ fixes the open end of $T$ determined by $x$. Also recall (e.g. from \cite{AlB}) that for
$\gamma\in Isom(T)$, $\ell(\gamma):=\min\{ d(t,\gamma(t)) ~|~ t\in T\}$ and
$A_{\gamma} := \{ t\in T ~|~ d(t,\gamma(t)) =\ell(\gamma)\}$.

Now let $x\in X$ and $\gamma\in\Gamma$ be given such that $\epsilon(\gamma)\in\epsilon(\Gamma)_x$ and
$sim(\epsilon(\gamma),x)=1$. We must show that $\epsilon(\gamma)$ fixes all points of $X$ sufficiently close to $x$.
Since $\gamma$ fixes the open end of $T$ determined by $x$, it follows that $x(t)\in A_{\gamma}$ for sufficiently
large $t$ (see \cite[Corollary 6.17]{AlB}). 

If $\gamma$ is elliptic (i.e., fixes some point of $T$), then $A_{\gamma}$ is the fixed point set of $\gamma$.
Thus, $\gamma\in\Gamma_{x(t)}$ for all sufficiently large $t$, say for $t\geq t_0$.
Let $t_1\geq t_0$ be given by the hypothesis. Then $\gamma\in\Gamma_{y(t)}$ whenever $y\in X$, $t\geq t_1$, and
$y(t_1)=x(t_1)$. That is, $\gamma$ fixes the open end determined by such a $y$ and hence, $\epsilon(\gamma)$
fixes $y$ for $y$ sufficiently close to $x$. 

On the other hand, if $\gamma$ is hyperbolic (i.e., has no fixed point), it follows that $A_{\gamma}$ is isometric to 
$\br$ (see \cite{AlB}). The condition $sim(\epsilon(\gamma),x)=1$ implies $A_{\gamma}=T$ and $x$ is isolated in $X$.
\end{proof}

\begin{remark}
The homomorphism $\epsilon:Isom(T)\to LS(X)$ is an injection if and only if
$T$ is not isometric to $\br$.
If $T$ is isometric to $\br$, one sees that $Isom(T)\cong \br\rtimes\bz/2$ and $LS(X)\cong \bz/2$. In
particular, $\epsilon$ is not injective.
Conversely, if $\epsilon$ is not injective, then there exists a non-trivial $\gamma\in\Gamma$ such that
$\epsilon(\gamma)$ is the identity on $X$. Using the notation in the proof of Theorem~\ref{theorem:tree-isom},
it follows that every point of $X$ determines an open end of $A_{\gamma}$ (see \cite[Corollary 6.17]{AlB}).
If $\gamma$ is hyperbolic, then $A_{\gamma}$ is isometric to $\br$ and $\gamma|A_{\gamma}$ is a translation;
it follows that $A_{\gamma}=T$. If $\gamma$ is elliptic, then $A_{\gamma}$ is the fixed point set of 
$\gamma$; it follows that $A_{\gamma}=T$, a contradiction.
\end{remark}

The following two corollaries follow immediately from Theorem~\ref{theorem:tree-isom}. 

\begin{corollary} Let $\Gamma$ be a subgroup of $Isom(T)$. If
$\{ t\in T ~|~ \Gamma_t\neq \{ 1\}\}$ is a bounded subset of $T$, then
$\epsilon(\Gamma)$ acts locally rigidly on $X$.
\end{corollary}

\begin{corollary} If $\Gamma\leq Isom(T)$ acts freely on $T$, then
$\epsilon(\Gamma)$ acts locally rigidly on $X$.
\end{corollary}

As an example, let $\Gamma$ be a finitely generated free group 
with a free set of generators $S=\{ s_1, s_2, s_3, \dots, s_n\}$
and identity element $e$.
Recall that the Cayley graph $\Cay(\Gamma,S)$ is a 
geodesically complete, rooted (at $e$), locally finite, simplicial tree
. Moreover, $\Gamma$ acts freely by isometries on $\Cay(\Gamma,S)$.

The following result follows from the previous corollary.

\begin{corollary} $\epsilon(\Gamma)$ acts locally rigidly on $end(\Cay(\Gamma,S),e)$.
\end{corollary}

\subsection{Local rigidity and tree lattices}
\label{subsec:tree-lattices}

We now characterize local rigidity of an ultrametric space $X$ in terms of the isometry group of $X$. This will be used 
below
when we discuss a connection with the Bass-Lubotzky theory of tree lattices. We will use the majorant topology on the isometry 
group.\footnote{Let $\mathcal{U}$ be an open cover of the metric space $X$, $h\in Isom(X)$ and define
$$N(h,\mathcal{U})=\{ g\in Isom(X) ~|~ \text{for each $x\in X$ there exists $U\in\mathcal{U}$ 
such that $g(x), h(x)\in U$}\}.$$
For fixed $h$ the sets $N(h,\mathcal{U})$ form a neighborhood basis for $h$ and the resulting topology is 
called the {\it majorant topology}.
Alternatively, consider continuous functions $\epsilon\co X\to (0,\infty]$. Then the sets
$$N_\epsilon(h) =\{ g\in Isom(X) ~|~ d(g(x), h(x) <\epsilon(h(x)) \text{~for all $x\in X$}\}$$
also form a neighborhood basis for $h$ in the majorant topology.}
This coincides with what is sometimes called the fine Whitney topology and makes the isometry group into a 
topological group (see \cite[Essay I, appendix C]{KS}). 
Of course, when $X$ is compact the majorant topology agrees with the compact-open topology and
the uniform topology (the usual sup-metric).

\begin{theorem} An ultrametric space $(X,d)$ is locally rigid if and only if the topological group
$Isom(X)$ is discrete in the majorant topology.
\end{theorem}

\begin{proof}
Assume first that 
$Isom(X)$ is discrete in the majorant topology.  
If $X$
is not locally rigid, then there exists $x\in X$ and a decreasing sequence $\{ \epsilon_i\}_{i=1}^\infty$ of
positive numbers converging to $0$ and non-trivial isometries $h_i\co B(x,\epsilon_i)\to B(x,\epsilon_i)$, $i=1,2,3,\dots$.
By Lemma~\ref{lemma:IsometryExtension} we can extend each $h_i$ to a non-trivial isometry $\tilde h_i\co X\to X$ such that
$\tilde h_i$ is the inclusion on the complement of $B(x,\epsilon_i)$.
It follows that $\tilde h_i$ converges to the identity on $X$ in $Isom(X)$ as $i\to\infty$ in the majorant topology, 
contradicting the discreteness of $Isom(X)$. 

Conversely, assume that $X$ is locally rigid. For each $x\in X$, let $\epsilon_x>0$ be given by 
Definition~\ref{def:locally-rigid}.
Let $\mathcal{U} = \{ B(x,\epsilon_x) ~|~ x\in X\}$ and suppose $g,h\co X\to X$ are two isometries that are $\mathcal{U}$-close.
If $y\in X$, then there exists $x\in X$ such that $g(y), h(y)\in B(x,\epsilon_x)$. 
It follows that $y\in g^{-1}B(x,\epsilon_x)\cap h^{-1}B(x,\epsilon_x)$.
Moreover, $g^{-1}B(x,\epsilon_x)=B(g^{-1}x,\epsilon_x)$ and
$h^{-1}B(x,\epsilon_x)= B(h^{-1}x,\epsilon_x)$. Hence, 
$g^{-1}B(x,\epsilon_x)=h^{-1}B(x,\epsilon_x)$ (by Proposition~\ref{prop:ElemProp}(1)).
Thus, $hg^{-1}|\co B(x,\epsilon_x)\to B(x,\epsilon_x)$ is an isometry and it follows that $hg^{-1}|$ is the identity
(by the choice of $\epsilon_x$); in particular, $g(y)=h(y)$. Thus, any two $\mathcal{U}$-close isometries on $X$ are equal and
$Isom(X)$ is discrete in the majorant topology.
\end{proof}

\begin{corollary} 
\label{cor:isom-finite}
A compact ultrametric space $X$ is locally rigid if and only if $Isom(X)$ is finite.
\end{corollary}

\begin{proof} If $X$ is compact, then so is $Isom(X)$. Thus, $Isom(X)$ is finite if and only if it is discrete.
\end{proof}

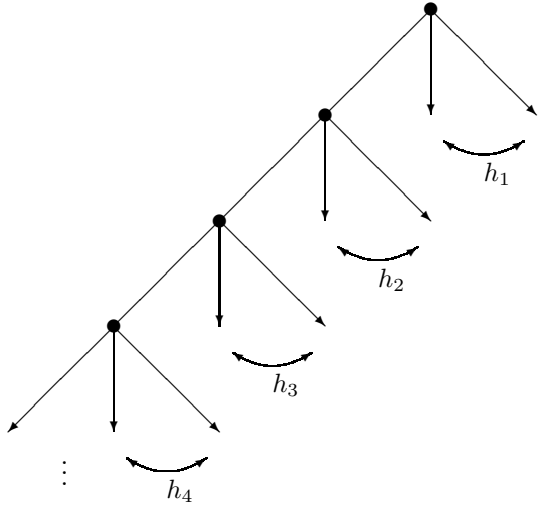
\begin{figure}

\begin{picture}(400,160)(-50,-20)


\put(0,0){\line(1,1){160}}
\put(0,0){\vector(-1,-1){0}}

\multiput(40,40)(40,40){4}{\line(1,-1){40}}
\multiput(80,0)(40,40){4}{\vector(1,-1){0}}
\multiput(40,40)(40,40){4}{\circle*{5}}
\multiput(40,40)(40,40){4}{\line(0,-1){40}}
\multiput(40,0)(40,40){4}{\vector(0,-1){0}}

\multiput(0,0)(40,40){4}{\qbezier(45,-10)(60,-20)(75,-10)}
\multiput(45,-10)(40,40){4}{\vector(-3,2){0}}
\multiput(75,-10)(40,40){4}{\vector(3,2){0}}
\put(60,-25){$h_4$}\put(100,15){$h_3$}\put(140,55){$h_2$}
\put(180,95){$h_1$}

\put(20,-20){$\vdots$}

\end{picture}

\caption{A non-discrete isometry group.
See Example~\ref{example:non-discrete}}
\label{figure:non-discrete}
\end{figure}

\begin{example}
\label{example:non-discrete} We give an example of a noncompact
ultrametric space 
$Y$ that is locally rigid, but $Isom(Y)$ is not
discrete in the compact-open topology (or, the uniform topology). 
Let $X$ be the space of Example~\ref{example:finite-not-lr}
and let $Y=X\setminus\{ x_\infty\}$.
For each $i\in\bn$, define an isometry $h_i\co Y\to Y$ by
$$h_i(z)= \left\{ \begin{array}{ll}
x_{|a-1|i} &\textrm{if $z=x_{ai}$ for $i\in\{ 0,1\}$}\\
z &\textrm{otherwise}.\\
\end{array} \right.$$
Then $h_i$ converges to the identity in the compact-open topology 
(hence, also in the uniform topology).
See Figure~\ref{figure:non-discrete}.
\end{example}

We now turn to some connections with some concepts encountered in the theory of tree lattices as developed by Bass and
Lubotzky \cite{BaL}.
Let $T$ be a locally finite simplicial tree.
The locally compact group $Aut(T)$ of simplicial auto-homeomorphisms of $T$ is a subgroup of the group $Isom(T)$ of
isometries of $T$ onto $T$.\footnote{In fact, $Aut(T) = Isom(T)$ if and only if 
$T$ is not isometric to $\br$.}${}^{,}$\footnote{The topology on $Aut(T)$ is the compact-open topology, 
so that two simplicial auto-homeomorphisms are close if
they agree on a large finite subtree. What is important about this topology is that discrete subgroups of
$Aut(T)$ are precisely the subgroups whose vertex stabilizers are finite.}
If $v\in T$ is a vertex and $Aut(T,v)\subseteq Aut(T)$ and $Isom(T,v)\subseteq Isom(T)$ are the subgroups of 
automorphisms fixing $v$, then $Aut(T,v)=Isom(T,v)$.

Fix a vertex $v\in T$, the root, and assume now that $(T,v)$ is geodesically complete.
The end space $end(T,v)$ of $(T,v)$ is a compact ultrametric space and the natural
function $Aut(T,v)=Isom(T,v)\to Isom(end(T,v))$ is an isomorphism 
(see Remark~\ref{cd:isom}). 

Recall the following definitions from \cite{BaL} and \cite{BaT}:

\begin{definition} 
\begin{enumerate}
\item A locally finite simplicial tree $T$ is {\it rigid} if $Aut(T)$ is discrete.
\item A finite, connected simplicial graph $K$ is {\it $\pi$-rigid} if $\pi_1(K) = Aut(\tilde K)$,
where $\tilde K$ is the locally finite simplicial tree that is the universal cover of $K$.
\end{enumerate}
\end{definition}

It follows that $T$ is rigid if and only if $Aut(T,v)$ is finite for all vertices $v\in T$,
if and only if $Aut(T,v)$ is finite for some vertex $v\in T$.

\begin{proposition} If $(T,v)$ is a geodesically complete, rooted locally finite simplicial tree, then $T$
is rigid if and only if the compact ultrametric space $end(T,v)$ is locally rigid.
\end{proposition}

\begin{proof}  By the remarks above, $T$ is rigid if and only if  the group 
$Isom(end(T,v))$ is finite. Hence, the result follows from Corollary~\ref{cor:isom-finite}.
\end{proof}

It follows that there is a rich source of examples of compact, locally rigid ultrametric spaces. In fact, Bass and Tits
\cite[page 185]{BaT} assert that a randomly constructed locally finite tree will have no non-trivial automorphisms. 
One therefore expects that almost all compact ultrametric spaces are locally rigid.

Bass and Kulkarni
\cite{BaK} and Bass and Tits
\cite{BaT} have 
provided examples of
$\pi$-rigid graphs $K$ without terminal vertices.
The universal covering trees $\tilde{K}$
of these graphs are rigid (because in this case $Aut(\tilde{K})$ acts freely on $\tilde{K}$) and
geodesically complete (with respect to any vertex).
Hence, these trees have end spaces that are compact, locally rigid ultrametric spaces. 

\subsection{R. J. Thompson's groups and their descendants}
\label{subsec:Thompson}

In this section we see how groups defined, generalized, and
developed by
Brown, Higman, Thompson, 
Neretin, R{\"o}ver, and others can be interpreted as groups of local similarities on compact
ultrametric spaces. 

For the groups $F$, $T$, and $V$ of Thompson \cite{Tho}, the connection
with the current paper arises from 
their description via {\it reduced tree diagrams}
in Cannon, Floyd, and Parry \cite{Can} based on work of
K. Brown \cite{Bro}. 
Higman \cite{Hig}  generalized $V$ to a family of infinite finitely presented
groups $G_{n,r}$ ($n=2,3,4,\dots$, $r=1,2,3,\dots$) and Brown \cite{Bro}
extended this to families
$F_{n,r}\leq T_{n,r}\leq G_{n,r}$ with $F_{2,1}=F$, $T_{2,1}=T$, and
$G_{2,1}= V$.
Brown, based on earlier work of J\'onsson and Tarski \cite{JonT},
used trees to describe these groups.
In particular, it is clear from Brown's work that each of these groups can be realized as subgroups
of groups of local similarities on end spaces of trees.

R{\"o}ver \cite{Rov}, with his notion of 
{\it almost automorphisms of trees},
further developed the viewpoint of Brown and described the groups $G_{n,1}$
as subgroups of homeomorphism  groups
of end spaces of trees (these homeomorphisms are local
similarities).
See also Greenberg and Serigiescu \cite{GreS} for an instance of tree
diagrams inducing groups of homeomorphisms on end spaces.

Neretin \cite{Ner84}, \cite{Ner92}, \cite{Ner03} introduced $p$-adic analogues
of the diffeomorphism group ${\rm Diff}(S^1)$ of the circle, called
groups of {\it spheromorphisms} and, later, {\it hierarchomorphisms},
that are also groups of homeomorphisms of end spaces of trees.
As is the case with R{\"o}ver's groups,
Neretin's groups are subgroups
of local similarities on the end space of a tree.
For more on Neretin's groups, see Kapoudjian \cite{Kap} and
Lavrenyuk and Sushchansky \cite{LavS}.

To indicate in a bit more detail how these groups are related to the current
paper, we need to introduce some more terminology.

Let $(T,v)$ be a rooted, geodesically complete, locally finite 
simplicial tree. 
For $i\in\bz_+$, $V_i$ denotes the
set of vertices of $T$ a distance $i$ from $v$ (the {\it vertices at
level $i$}) as in Section~\ref{subsec:ConstructionOfBTv}, and 
$\mathcal E_i$ denotes the set of edges of $T$ with one vertex in 
$V_i$ and the other in $V_{i+1}$ 
(the {\it edges at level $i$}).
For a vertex $w\in V_i$, let $\mathcal E_w\co =\{ E\in\mathcal E_i 
~|~ w\in E\}$.

An {\it order} of $(T,v)$ consists of a total order on $\mathcal E_w$
for each vertex $w\in T$, and $(T,v)$ is {\it ordered} if it comes with
an order.

Note that an order of $(T,v)$ induces a total order on $V_i$
for each $i=0,1,2,\dots$ defined inductively as follows.
If $v_1$ and $v_2$ are distinct vertices in 
$V_{i+1}$, let $E_1, E_2$ be the unique edges in $\mathcal E_i$
with $v_1\in E_1$, $v_2\in E_2$, and let $w_1\in E_1$, $w_2\in E_2$
be the vertices of these edges in $V_i$.
If $w_1=w_2$, then define $v_1<v_2$ if and only if $E_1<E_2$; if 
$w_1\neq w_2$, define $v_1<v_2$ if and only if $w_1<w_2$.

Note that every $(T,v)$ can be ordered.
Furthermore, an order on the infinite $n$-ary tree $(A_n,v)$ is often implicitly 
used: the set of
vertices immediately below a given vertex is identified with the 
ordered set $\{0,1,2,\dots n-1\}$.

If $(T,v)$ is ordered, then there is an induced (dictionary) total order
on $X=end(T,v)$. This is because $X$ is identified with the set of all
sequences $(v_0,v_1,v_2,\dots)$ such that $v_i\in V_i$ for each $i$.
Thus, we can speak of {\it order preserving} maps $X\to X$.
Moreover, a map $h\co X\to X$ is {\it locally order preserving} if
for each $x\in X$ there exists $\epsilon > 0$ such that
$h|\co B(x,\epsilon)\to X$ is order preserving.

There also is an induced total order on any collection of disjoint balls
in $X$ (this uses the description of balls in Section~\ref{sec:ends}).
In particular, such collections have a cyclic order.
Therefore, we say that a local similarity $h\co X\to X$ is
{\it cyclic order preserving} if there exists $\epsilon >0$ such that
for every $x\in X$ there exists $\lambda_x>0$ so that
$h|\co B(x,\epsilon)\to B(hx,\lambda_x\epsilon)$ is a similarity, and the
induced function
$\{ B(x,\epsilon) ~|~ x\in X\}\to\{ B(x,\lambda_x\epsilon) ~|~ x\in X\}$,
$B(x,\epsilon)\mapsto B(x,\lambda_x\epsilon)$, preserves the cyclic
order.

For the rest of this section, let $(T,v)$ be an ordered, rooted, geodesically complete,
locally finite simplicial tree, and let $X= end(T,v)$.
We denote various subgroups of $LS(X)$ as follows:
\begin{eqnarray*}
LS^{o.p.}(X) & = & \{ h\in LS(X) ~|~ \text{$h$ is order preserving}\}\\
LS_{l.o.p.}(X) & = & \{ h\in LS(X) ~|~ \text{$h$ is locally order preserving}\}\\
LS^{o.p.}_{l.o.p.}(X) & = & LS^{o.p.}(X)\cap LS_{l.o.p.}(X)\\
LS^{c.o.p.}(X) & = & \{ h\in LS(X) ~|~ \text{$h$ is cyclic order preserving}\}\\
LS^{c.o.p.}_{l.o.p.}(X) & = & LS^{c.o.p.}(X)\cap LS_{l.o.p.}(X)
\end{eqnarray*}

Although R{\"o}ver \cite{Rov} focused on spherically
homogeneous trees, it is clear that his definition 
of the almost automorphism group AAut$(T,v)$ can be made for any $(T,v)$
as above and that AAut$(T,v) = LS_{l.o.p}(X)$.

Likewise, Neretin's \cite{Ner03}  group
$Hier(T,\Gamma)$ of {\it hierarchomorphisms}, where $\Gamma\leq Isom(T)$,
is a subgroup of $LS(X)$.

\begin{example}
Let $X_n=end(A_n,v)$, the end space of the infinite $n$-ary tree.
Then,
$LS_{l.o.p}(X_n) = G_{n,1}$, the Higman--Thompson group.
In particular, there are the following
interpretations of Thompson's groups:
\begin{eqnarray*}
LS_{l.o.p}(X_2) &=& G_{2,1} ~=~ V\\
LS^{o.p.}_{l.o.p.}(X_2) &=& F\\
LS^{c.o.p.}_{l.o.p.}(X_2) &=& T
\end{eqnarray*}
\end{example}

\begin{proposition} 
\label{prop:ordered-loc-rig}
If $(T,v)$ is an ordered, rooted, geodesically complete, 
locally finite simplicial tree and $X=end(T,v)$, then
$LS_{l.o.p.}(X)$ acts locally rigidly on $X$.
\end{proposition}

\begin{proof}
Let $h\in LS_{l.o.p}(X)$ and $x\in X$ such that
$hx=x$ and $sim(h,x)=1$. Thus, there exists $\epsilon > 0$
such that $h|\co B(x,\epsilon)\to B(x,\epsilon)$ is an order preserving
isometry. This implies that $h|$ is the identity
(clearly, using the notation of Section~\ref{sec:ends}, $h$
fixes the vertex $\langle x,\epsilon\rangle\in T$ and each edge of
$\mathcal E_{\langle x,\epsilon\rangle}$; continuing by induction, $h$
fixes each vertex of $T_{\langle x,\epsilon\rangle}$).
\end{proof}

\begin{corollary}[Birget \cite{Bir}, Nekrashevych \cite{Nek}] 
For each $n=2,3,4,\dots$, there exists a 
faithful unitary representation of the Higman--Thompson group
$G_{n,1}$ into the Cuntz algebra $\mathcal O_n$.
\end{corollary}

\begin{proof}
By Proposition~\ref{prop:ordered-loc-rig},
$\Gamma=LS_{l.o.p}(X_n)$ acts locally rigidly on
$X_n = end(A_n,v)$. The groupoid
${\mathcal{G}_{\Gamma}(X_n)} = O_n$, the Cuntz groupoid, and
$C^*{\mathcal{G}_{\Gamma}(X_n)}=\mathcal{O}_n$, the Cuntz algebra, by
Renault \cite{Ren}. Hence, the corollary follows from
Theorem~\ref{theorem:second-main}(2)
(i.e., Theorem~\ref{theorem:FUR}).
\end{proof}

\subsection{Symbolic dynamics}
\label{section:SymbolicDynamics}

This brief section 
contains a preliminary comparison between the concept
of local isometry in ultrametric spaces
and the concept of tail equivalence studied in 
symbolic dynamics.

Let $A$ be an $n\times n$ matrix with entries $A_{ij}\in\{ 0,1\}$, $1\leq i,j\leq n$.
Let $X_A$ be the one-sided subshift of finite type (i.e., the one-sided topological Markov
shift) with transition matrix $A$. 
Thus,
$$X_A= \{ (x_i)_{i=1}^\infty ~|~ \textrm{for each $i$, 
$x_i\in\{ 0,1,\dots,n-1\}$ and 
$A_{x_i,x_{i+1}}=1$}\}.$$
Define a metric on $X_A$ by
$$d(x,y)= \left\{ \begin{array}{ll}
e^{-k} & \textrm{if $x_i=y_i$ for $1\leq i< k$ and $x_k\neq y_k$}\\
0 &\textrm{if $x=y$.}\\
\end{array} \right.$$
As is well known, $(X_A,d)$ is a compact ultrametric space.

\begin{proposition}
\label{prop:tail-equivalent}
If $x,y\in X_A$ and $x$ and $y$ are tail equivalent, then
there exist $\epsilon >0$ and an isometry $h\co B(x,\epsilon)\to B(y,\epsilon)$ such that
$hx=y$.
\end{proposition}

\begin{proof} Suppose $k\in\bn$ and $x_i=y_i$ for all $i\geq k$.
Define $h\co B(x,e^{-(k+1)})\to B(y,e^{-(k+1)})$ by 
$$h(z)_i=\left\{ \begin{array}{ll}
y_i & \textrm{if $1\leq i\leq k-1$}\\
z_i & \textrm{if $k\leq i$,}\\
\end{array} \right. ~~~~\textrm{for each $z\in B(x,e^{-(k+1)})$}.
$$
It is easy to check that $h$ is the desired isometry.
\end{proof}

\begin{example} The converse of Proposition~\ref{prop:tail-equivalent}
does not hold in general. For example,
if 
$A = \left[ \begin{array}{cc} 1 & 1\\ 1&1 \end{array} \right]$,
then $X_A$ is the end space of the Cantor tree and is
isometrically homogeneous, i.e., the isometry group of $X_A$ acts 
transitively.
Another example is provided by the matrix
$A=\begin{bmatrix}
1 & 1 & 1\\
1 & 1 & 1\\
1 & 0 & 0\\
\end{bmatrix}$.
In this case, $X_A$ is not isometrically homogeneous, but equal local 
isometry type of points need not imply tail equivalence.
\end{example}



{\footnotesize
\bibliographystyle{plain}
\addcontentsline{toc}{section}{\refname}
\bibliography{biblio}
}

\end{document}